\documentclass[11pt,leqno]{article}
\usepackage[russian, english]{babel}
\usepackage{graphicx}
\usepackage{epstopdf}
\usepackage{amsfonts, amsthm, amsmath}
\usepackage{amsxtra}
\usepackage{tikz}
\usetikzlibrary {positioning}
\usetikzlibrary {positioning}
\usepackage {xcolor}
\definecolor {processblue}{cmyk}{0.96,0,0,0}
\usetikzlibrary{calc}
\tikzstyle{stuff_nofill}=[circle,draw]
\usepackage{fontenc}
\usepackage[dvips]{epsfig}
\newtheorem{theo}{Theorem}[section]
\newtheorem{defi}[theo]{Definition}
\newtheorem{lemm}[theo]{Lemma}
\newtheorem{prop}[theo]{Proposition}

\newtheorem{cor}[theo]{Corollary}

\def\part{P} %partition

\newcommand{\lineann}[5][0.5]{%
    \begin{scope}[ blue,inner sep=2pt]
        \draw[dashed, blue!40] (#2,0) -- +(0,#1)
            node [coordinate, near end] (a) {};
        \draw[dashed, blue!40] (#3,0) -- +(0,#1)
            node [coordinate, near end] (b) {};
        \draw[<->] (a) -- node[] {#4} (b);
    \end{scope}
}
\newcommand{\Anm}[2][0.5]{\draw[|-|,red, very thick] (#2-#1,0)--(#2+#1,0);}

\begin{document}
\begin{center}
\textbf{\large The thermodynamic formalism and central limit theorem for stochastic perturbations of
circle maps with a break} \\
\vspace{.25in}\large{
Akhtam Dzhalilov\footnote{Natural-Mathematical Science Department, Turin Polytechnic University, Kichik Halqa Yoli 17,  Tashkent 100095, Uzbekistan. adzhalilov21@gmail.com}, Dieter  Mayer\footnote{Institut f\"ur Theoretische Physik, TU
Clausthal, Leibnizstrasse 10,  D-38678 Clausthal-Zellerfeld, Germany.  dieter.mayer@tu-clausthal.de}, Abdurahmon Aliyev\footnote{ V.I.Romanovsky Institute of Mathematics, Academy of Sciences,  Beruniy street 369, Tashkent 100170, Uzbekistan. aliyev95@mail.ru}
}
\end{center}

{\small \textbf{Abstract.}

Let
 $T\in C^{2+\varepsilon}(S^{1}\setminus\{x_{b}\}),\,\,\varepsilon>0,$ be an orientation preserving circle
 homeomorphism  with rotation number  $\rho_T=[k_{1},k_{2},..,k_{m},1,1,...],\,\,m\geq1$,
 and
 a single break point $x_{b}$, at which the derivative of $T$ has a jump.
We consider the stochastic sequence
 $ \overline{z}_{n+1} = T(\overline{z}_{n}) + \sigma \xi_{n+1},\,\overline{z}_{0}:=z\in S^1,$
where $\{\xi_{n},\,n=1,2,...\}$ is a sequence of real valued independent mean zero random
variables of comparable sizes, and $\sigma > 0$ is a small parameter -called
noise level- which controls the size of the noise.
Such stochastic perturbations of one dimensional interval maps, among them critical circle maps,  have been studied some time ago by Diaz-Espinoza and de la Llave, who showed for the resulting sum of random variables  a central limit theorem and determined the rate of convergence  to the Gaussian distribution. Their approach used the renormalization group technique well known from Feigenbaum's universality theory for the transition to chaos in critical interval maps. On the other hand Sinai and coworkers developped a thermodynamic formalism approach to such maps with a weak noise at the accumulation of period doubling to study their ergodic properties like the existence of invariant measures. We will use this formalism in the present paper to extend the result of Diaz-Espinoza and de la Llave to circle homemorphisms with a break point.   A thermodynamic formalism for circle homeomorphisms with a break point was recently constructed by  A. Dzhalilov et al.. This formalism and the sequence of dynamical partitions $P_n(T,x_b)$, determined by the trajactory  of the break point $x_b$ under $T$, allows us, following earlier work of Vul at al., to establish a symbolic dynamics for any point $z\in S^1$ and to define a transfer operator whose leading eigenvalue can be used to bound the Lyapunov function. For a special sequence $n_m, m\to\infty$ we can show, that the barycentric coefficient of any $z_k$ in the orbit of a point $z_0\in S^1$, which does not intersect  the orbit of the break point $x_b$, is universally bounded in the interval of the partition $P^{n_m}(T,x_b)$ to which it belongs. To prove the central limit theorem and to determine the convergence rate to the Gaussian distribution  we expand the stochastic sequence in a Taylor expansion in the variables $\xi_i$ which leads to the decomposition
$ \overline{z}_{n} (z_0,\sigma)= T^n({z}_{0}) + \sigma L_n(z_0) + \sigma^2 Q_n(z_0,\sigma),$ where $L_n(z)= \xi_n+\sum\limits_{k=1}^{n-1}\xi_k\prod\limits_{j=k}^{n-1} T'(z_j), z\in S^1$ is the linearized effective noise and $Q$ describes the higher order terms. This can be done however only in a neighbourhood of the point $z_n=T^n(z_0)$ which does not contain the break points of $T^n$. For this we use the sequence $\{n_m\}$ and construct a series of neighbourhoods $A_k^{n_m}$ of the points $z_k$ composed of elements of the partitions $P_{n_m}(T,x_b)$ which do not contain any break point of the map $T^{q_{n_m}}$, where the $q_n$ are the so called first return times of $T$, determined by the continued fraction expansion of the rotation number $\rho_T$ of the circle map $T$. Following  Diaz-Espinosa and de la Llave we prove the central limit theorem for the lienearized process which leads finally also to the proof of our extension of their results for circle homeomorphisms with a break point.

%\textbf{Keywords:} circle map, rotation number, break point, stochastic perterbution, central limit theorem, thermodynamic formalism.

%\textbf{Mathematics Subject Classification (2010):} 37C05; 37C15; 37E05; 37E10; 37E20; 37B10.

\makeatletter
\renewcommand{\@evenhead}{\vbox{\thepage \hfil {\it Dzhalilov A., Aliev A.}   \hrule }}
\renewcommand{\@oddhead}{\vbox{\hfill
        {\it The thermodynamic formalism and central limit theorem for stochastic perturbations of
circle maps with a break  }\hfill \thepage \hrule}} \makeatother

\label{firstpage}

\section{Introduction}

 Dynamical systems theory is mostly interested in describing the typical
behaviour of orbits as time goes to infinity, and to understand how this behaviour is
modified under small perturbations of the system. In the present work we study  stochastic perturbations
of circle maps with one break point, using as the  main tool  the  thermodynamic
formalism.
Ya.G. Sinai  constructed in [19] the first example of a thermodynamic
formalism  for Anosov's flows, which was generalized later
 in the works of D. Ruelle [18], R. Bowen [2] and others for Smale's Axiom A systems.
 E.B. Vul, Ya.G. Sinai  and K.M. Khanin finally succeeded in establishing  in [20] a thermodynamic formalism approach to Feigenbaum universality in families of critical interval maps, the first example of such an approach to  non-hyperbolic systems, whereas the standard approch by Feigenbaum and others has been the renormalization group well known from statistical mechanics.

A natural generalization of smooth interval maps and circle diffeomorphisms are piecewise smooth circle homeomorphisms
with break points (see [14]). Contrary to diffeomorphisms the invariant measure of these circle homeomorphisms $T\in C^{2+\varepsilon}(S^{1}\setminus\{x_{b}\}),\,\,\varepsilon>0$,
with break point $x_{b}$ and irrational rotation number is singular w.r.t. Lebesque measure \cite{DK1998}. The
renormalizations of such maps are exponentially approximated by fractional-linear maps \cite{KV1991}.
Consider two homeomorphisms $ T_{1} $ and $ T_{2} $ with the same irrational rotation number $ \rho = \rho (T_{1}) = \rho (T_{2}), $ and with identical breakpoint $ x_{b}=x_1=x_2. $ The question of the regularity of the conjugation  $ \Phi $ between $ T_{1} $ and $ T_{2} $ is called the  rigidity  problem. It has been intensively studied in the works of \cite{KhnKhm}, \cite{KhKo} and others.

\begin{theo}[see \cite{KhnKhm}]  \label{teor-1}
 Let $T_{1},\,\,T_{2}\in C^{2+\varepsilon}(S^{1}\setminus \{x_{b}\}),\,\,\,\varepsilon>0$, be circle homeomorphisms with break point $x_{b}$. Suppose that

{\rm 1)} they have the same rotation number $\rho (T_{1})=\rho (T_{2})=\rho $;

{\rm 2)} $ \rho $ is irrational and has a periodic continued fraction expansion of the form
$$\rho=[k_{1},k_{2},\dots,k_{s},k_{1},k_{2},\dots,k_{s},\dots],\,\,\,s\geq 1.$$

Then the conjugating homeomorphism $ \Phi $ between $ {{T}_{1}} $ and $ {{T}_{2}} $ belongs to the class $ C^{1+ \theta} (S^{1}), $ where $ \theta> 0 $ depends only on the rotation number $ \rho $.
\end{theo}

 J. Crutchfield et al. and B. Shraiman et al. considered in
\cite{CNR81} respectively \cite{SWM81} heuristically a renormalization group - respectively a field theoretic path-integral approach for weak Gaussian
noise perturbing one dimensional maps with period doubling at the onset of chaos. The main
result in those papers was that after appropriately rescaling space and time, the Lyapunov exponent at the transition  satisfies some scaling relations. Vul et al.
developed in \cite{VSK84} a rigorous thermodynamic formalism approach for critical maps with period doubling.
Among many other results these authors  studied the effect of noise on the ergodic properties of
these maps and showed that for systems with weak
noise at the accumulation of period doubling there is a stationary measure, depending on the magnitude of the noise, which converges for vanishing noise to
the invariant measure of the attractor.

O. Diaz-Espinosa and R. de la Llave studied in \cite{Llave}  stochastic perturbations of several systems using the   renormalization group technique.
Among others they proved a central limit theorem for critical circle maps with golden mean rotation number and some mild conditions on the stochastic noise.

Before turning  to the formulation of the main results of our work we recall
the general setup and more details of the two main results of  O. Diaz-Espinosa and R. de la Llave in \cite{Llave}.

Let $(\Omega,\mathcal{F}, P)$ be a probability space and $T:S^1\to S^1$ a homeomorphism of the circle $S^1 \to S^1$.
Let the stochastic sequence be defined as
\begin{equation}\label{Peterb1}
 \overline{x}_{n+1} = T(\overline{x}_{n}) + \sigma \xi_{n+1},\,\overline{x}_{0}:=x\in S^1
 \end{equation}
where $(\xi_n)$ is a sequence of independent random variables  with $p>2$ finite moments satisfying the following conditions:
\begin{equation}\label{Cond1}
 E\xi_n=0;
 \end{equation}
\begin{equation}\label{Cond2}
  const\leq(E|\xi_n|^2)^{1/2}\leq(E|\xi_n|^p)^{1/p}\leq Const.
 \end{equation}
 In the following we set
 $$x_n=T^n(x_0),\,n\geq 1,$$
 where $T^n$ denotes the $n$-th iteration of $T$.
%We define
%$$x_n=T(x_{n-1}),\,\,x_0=x. $$
The linearized effective noise is
defined as
\begin{equation}\label{Ln}L_n(x)=\xi_n+\sum\limits_{k=1}^{n-1}\xi_k\prod\limits_{j=k}^{n-1}T'(x_j),\,x\in S^1.\end{equation}
  Let $\omega_{n}(x,\sigma)$ be the stochastic process defined by
  \begin{equation}\label{proc}\omega_n(x,\sigma)=\frac{\bar{x}_n-x_n}{\sigma \sqrt{\text{var}(L_n(x)) }}.\end{equation}\\
  For an arbitrary  $z_0\in S^{1}\setminus \{T^{i}(x_b),\,i=0, -1,-2,...\}$ and
 each $s\geq 0,\,n\geq 1$ the Lyapunov functions $\Lambda_s(z_0,n)$ and $\hat{\Lambda}(z_0,n)$ are defined as follows
\begin{equation}\label{Lambda}\Lambda_s(z_0,n)=1+\sum\limits_{k=1}^{n-1}\prod\limits_{j=k}^{n-1}\left|T'(z_j)\right|^s\end{equation}
\begin{equation}\hat{\Lambda}(z_0,n)=\max\limits_{1\leq i\leq n-1}\sum\limits_{k=1}^{i}\prod\limits_{j=k}^{i}\left|T'(z_j)\right|\end{equation}\\
Using the renormalization group approach O. Diaz-Espinosa and R. de la Llave   established in \cite{Llave} a sufficient condition for the following CLT to hold for the sequence of random variables $\omega_n(x,\sigma_n) $ defined by certain one-dimensional maps

\begin{theo}\label{Llave1}(see \cite{Llave}). Let $T:M\to M$ be a $C^{2}$ map for $M=R^{1},\,I=[-1,1]$ or $S^{1}$
and let $\{\xi_{n},\,n=1,2,...\}$ be a sequence of independent random variables with $p>2$ finite moments.
%satisfying conditions (\ref{Cond1}) and (\ref{Cond2}).
Suppose that for some $x\in M$ there is an increasing sequence of positive integers $n_k$ such that
\begin{equation}\label{maincon}
\lim\limits_{k\to\infty}\frac{\Lambda_p(x,n_k)}{(\Lambda_2(x,n_k))^{p/2}}=0.
\end{equation}
Let $\sigma_k$ be a sequence of positive numbers. Assume furthermore
either of the following two conditions\\
$\bullet$ $[H1]$
 the noise satisfies conditions (\ref{Cond1}) and (\ref{Cond2}) with $p > 2$
and the sequence $\sigma_k$ satisfies
\begin{equation}\label{2con}
\lim\limits_{k\to\infty}\frac{\sup\limits_{x\in S^1}|T''(x)| \|\max\limits_{1\leq j\leq n_k}\xi_j\|_p^2(\hat{\Lambda}(x,n_k))^6\sigma_k}{\sqrt{\Lambda_2(x,n_k)}}=0;
\end{equation}
$\bullet$ $[H2]$
 the noise satisfies conditions (\ref{Cond1}) and (\ref{Cond2}) with $p\geq 4$
and the sequence $\sigma_k$  satisfies
\begin{equation}\label{2cond}
\lim\limits_{k\to\infty}\frac{\sup\limits_{x\in S^1}|T''(x)| \|\max\limits_{1\leq j\leq n_k}\xi_j\|_p^2(\hat{\Lambda}(x,n_k))^3\sigma_k}{\sqrt{\Lambda_2(x,n_k)}}=0.
\end{equation}
Then there exists a sequence of events $B_k\in \mathcal{F}$ such that
\begin{itemize}
  \item M1.\,\,\, $\lim\limits_{k\to\infty} P(B_k)=1$
   \item M2.\,\,\, The two processes defined by
     \begin{equation}\label{proce}\omega_{n_k}(x,\sigma_{k})=\frac{\bar{x}_{n_k}-x_{n_k}}{\sigma_{k}\sqrt{\text{var}(L_{n_k}(x)) }}.\end{equation}
      \begin{equation}\label{proc1}  \tilde{\omega}_{n_k}(x,\sigma_{k})=\frac{(\bar{x}_{n_k}-x_{n_k})\mathbf{1}_{B_k}}{\sqrt{\text{var}((\bar{x}_{n_k}-x_{n_k})\mathbf{1}_{B_k}) }}.\end{equation}
      converge in distribution to a standard Gaussian as $k\to \infty$.
\end{itemize}

      If furthermore  the sequence $\xi_{n}$ is supported on a compact set then
we can choose $B_{k} =\Omega$
 for all $k$.
\end{theo}
For the rate of convergence to the Gaussian in this Theorem  these authors got the following result, where $\Phi(z)$ denotes the distribution of the standard Gaussian on the real line:
\begin{theo}\label{Llave2}(see \cite{Llave}).Let $T, \xi_n$  be as in Theorem \ref{Llave1} and let $s = \min(p, 3)$.
Assume that condition (\ref{maincon}) holds at some $x\in M$. If $\sigma_k$ is a sequence
of positive numbers such that
\begin{equation}
\frac{(\hat{\Lambda}(x,n_k))^3}{\sqrt{\text{var}(L_{n_k}(x))}}\sup\limits_{x\in S^1}|T''(x)| \|\max\limits_{1\leq j\leq n_k}\xi_j\|_s^2\sigma_k\leq \left(\frac{\Lambda_s(x,n_k)}{(\Lambda_s(x,n_k))^{s/2}}\right)^2,
\end{equation}
then, we have
\begin{equation}
\sup\limits_{z\in R}|P(\omega_{n_k}(x,\sigma_k) \mathbf{1}_{B_k}\leq z)-\Phi(z)|\leq A\frac{\Lambda_s(x,n_k)}{(\Lambda_s(x,n_k))^{s/2}},
\end{equation}
where the constant $A>0$ depends only on $x$.
\end{theo}
Our main result of the present paper is
\begin{theo}\label{main}Let $T\in C^{2+\varepsilon}(S^{1}\setminus \{x_{b}\}),\,\,\,\varepsilon>0$, be a circle homeomorphism with  break point $x_{b}$, $T'(x)\geq Const >0,\,x \in[x_{b},x_{b}+1] $  and rotation number $\rho_T=[k_{1},k_{2},..,k_{m},1,1,...],\,\,m\geq1$.
  Consider  a sequence of independent random variables $(\xi_n)$ with $p>2$ finite moments satisfying the  conditions (\ref{Cond1}) and (\ref{Cond2}) for some  $x\in S^{1}\setminus \{T^{i}(x_{b}),\,i=0, -1,-2,...\}$.
Then
 \begin{enumerate}
   \item there exists a constant $\gamma>0$ such that if
   \begin{equation}\label{sigman}
   \lim\limits_{n\to\infty}\sigma_{n} n^{\gamma}=0,
   \end{equation}
the process $\omega_{q_n}(x,\sigma_{q_n})$ defined by (\ref{proce}) converges in distribution
to the standard Gaussian<
   \item furthermore, there are constants  $\tau>0$ and $\kappa>0$ depending
on $p$ and  a constant $C_1>0$ such that if $\sigma_n\leq C_1 n^{-\tau}$, then
$$\sup\limits_{z\in \mathbb{R}}\left|P(\omega_{q_n}(x,\sigma_{q_n})\leq z)-\Phi(z)\right|\leq Cq_n^{-\kappa},$$
where  $q_n$ is the first return time of $T$, the constant $C>0$  depends only on $x$ and $\Phi(z)$ is the distribution function of the standard Gaussian on $\mathbb{R}$.
 \end{enumerate}
 \end{theo}

Remarks.  Our paper is strongly influenced by the work of O. Diaz-Espinosa and R. de la Llave in \cite{Llave}.  Following  their ideas  we prove  the analog of their Theorem \ref{Llave1} for circle maps with a break point with small modifications. For estimating the  Lyapunov function however  we are using the thermodynamic formalism for such maps.\newline
We restrict our discussion to the simplest case of circle maps with one break point and eventually golden mean rotation number. Since the thermodynamic formalism in our approach  can be extended to  rotation numbers with arbitrary eventually periodic continued fraction expansions, our constructions can be in principle generalized to such cases. Technically these become however more involved and would lead to a longer paper.\newline
We emphasize that the limit theorems in this paper are considered in the setup of sums of random variables. 
\section{Preliminaries and notations}

Let $T$ be an orientation preserving circle homeomorphism with  irrational rotation number $\rho_{T}$.
Then $\rho_{T}$ can be uniquely expanded as a
continued fraction i.e.
$\rho_{T}=1/(k_{1}+1/(k_{2}+...)):=[k_{1}, k_{2},...,k_{n},...)$.
Denote by  $p_{n}/q_{n}=[k_{1}, k_{2},...,k_{n}],$ $n\geq 1,$  its  $n$-th convergent.
  The numbers $q_{n}, n\geq1$ are  called the  \textbf{first return times}  of $T$
  and satisfy the recurrence relations
  $q_{n+1}=k_{n+1}q_{n}+q_{n-1}\; n\geq 1,$
  where  $q_{0}=1$ and  $q_{1}=k_{1}.$
  Fix  an arbitrary point $z_{0}\in S^{1}.$  Its forward orbit
$ O^{+}_{T}(z_{0})=\{z_{i}=T^{i}(z_{0}),\,\,i=0,1,2...\}$ defines a sequence of natural partitions of the circle.
Indeed, denote by $I_{0}^{(n)}(z_{0})$ the closed interval in $S^{1}$ with endpoints $z_{0}$ and
$z_{q_{n}}=T^{q_{n}}(z_{0})$. In the clockwise orientation of the circle the point
$z_{q_{n}}$ is then for $n$ odd to the left of $z_0$, and for $n$ even to its
right. If $I_{i}^{(n)}(z_{0})=T^{i}(I_{0}^{(n)}(z_{0})),
i\geq 1$, denote the iterates of the interval $I_{0}^{(n)}(z_{0})$ under
$T$, it is well known, that the set $P_{n}(z_{0})$ of intervals
with mutually disjoint interiors, defined as
$$
P_{n}(z_0)=\{I_{i}^{(n)}(z_{0}),\,\,0\leq i <
q_{n+1}\}\cup\{I_{j}^{(n+1)}(z_{0}),\,\,0\leq j < q_{n}\}
$$
determines a partition of the circle for any $n$. The partition
$P_{n}(z_{0})$ is called the $n$-th \textbf{dynamical partition} of $S^1$ determined by
 the point $z_{0}$ and the map $T$.

  Proceeding from partition
$ P_{n}(z_{0})$ to $P_{n+1}(z_{0})$ all the
intervals $I_{j}^{(n+1)}(z_{0}),\,\,0\leq j\leq
q_{n}-1$, are preserved, whereas each of the intervals
$I_{i}^{(n)}(z_{0}),\; 0\leq i \leq q_{n+1}-1$, is partitioned into $k_{n+2}+1$
subintervals belonging to $P_{n+1}(z_{0})$, such
that %(see Fig 2.1)
$$ I_{i}^{(n)}(z_{0})=
I_{i}^{(n+2)}(z_{0})\cup\bigcup\limits_{s=0}^{k_{n+2}-1}I_{i+q_{n}+sq_{n+1}}^{(n+1)}(z_{0}).
$$
Obviously one has $P_{1}(z_{0})\prec P_{2}(z_{0})\prec ...\prec
P_{n}(z_{0})\prec ....$

The intervals $ I_{0}^{(n)}(z_{0}),\, I_{0}^{(n+1)}(z_{0})$ are called \textbf{generators}
of the partition  $P_{n}(z_{0}).$}

Later we will use also the so called \textbf{renormalization intervals} $J^{(n)}_i (z_0)= T^i (J^{(n)}_0(z_0))= J^{(n)}_0(z_i),\; i=0,1,2,\dots$, where $J^{(n)}_0 (z_0)=I^{(n)}_0 (z_0)\cup I^{(n+1)}_0 (z_0)$ and  $z_i=T^i(z_0)$ . Define the $Poincar\acute{e}\, map$ $\pi_n: J_0^{(n)}(z_0)\rightarrow J_0^{(n)}(z_0)$ by
$$\pi_n(x)=\left\{
               \begin{array}{ll}
                 T^{q_{n+1}}x, & \hbox{if  } x\in I_0^{(n)}(z_0), \\
                 T^{q_{n}}x, & \hbox{if  } x\in I^{(n+1)}_0(z_0).
               \end{array}
             \right.
$$

The following lemma plays a key role for
studying the metrical properties of the homeomorphism $T.$

\begin{lemm}\label{DEN}(see \cite{DK1998}) Let $T$ be a  circle homeomorphism with one break point $x_b$ with jump ratio $c_{T}(x_{b})=\sqrt{\frac{T(x_{b}-)}{T(x_{b}+)}}\neq 1$ and irrational rotation number. Suppose $T\in C^1 ([x_b, x_b+1])$ and $\underset{z\in[x_b,x_b+1]}{\text{var}} \ln T'=\bar{v}<\infty.$ Put $v=\bar{v}+2|\ln c_{T}(x_{b})|.$ If $y_0\in S^1$ and $T^i(y_0)\neq x_b, \, 0\leq i<q_n,$ then
\begin{equation}\label{DenInq}
e^{-v}\le \prod\limits_{s=0}^{q_{n}-1} DT(y_{s})\le e^{v}.
\end{equation}
holds.
\end{lemm}

Inequality (\ref{DenInq}) is called \textbf{Denjoy's inequality}.

It follows
from Lemma \ref{DEN}, that the intervals of the dynamical partition
$P_{n}(z_{0})$ have exponentially small lengths.
Indeed one finds
\begin{cor}\label{LDB} Suppose the circle map $T$ satisfies the conditions of Lemma \ref{DEN}.
Then for  an arbitrary element $I^{(n)}$  of
the dynamical partition $P_{n}(z_0)$ the following bounds hold
\begin{equation}\label{EDL}
\ell(I^{(n)})\leq const \,\theta^{n},
\end{equation}
 where   $\theta=(1+e^{-v})^{-\frac{1}{2}}<1$ and $l$ denotes Lebesque measure.
\end{cor}

Corollary \ref{LDB} implies that the trajectory
of every point $x\in S^1$ is dense in $S^1.$ This together with monotonicity of $T$
implies that the homeomorphism $T$ is topologically conjugate to the linear rotation
$T_{\rho}(x) = x + \rho \, \mod 1$.

\begin{defi}
Let $K>1$ be a constant. We call two intervals $I_1$ nd $I_2$ of the circle $S^1$  $K$-comparable if the inequality $K^{-1}|I_2|\leq |I_1|\leq K|I_2|$ holds.
\end{defi}
\begin{lemm}\label{compare} Suppose the circle homeomorphism  $T$ satisfies the conditions of Lemma \ref{DEN} and $z_{0}\in S^{1}.$
Then for an arbitrary interval  $I^{(n)}$ of the dynamical partition $P_{n}(z_{0})$  at least $n-1$ elements of $P_n(T,z_0)$ are $e^v-$comparable with $I^{(n)}$.
\end{lemm}
\proof
Let   $I^{(n)}_{j_0}\in P_n(T,z_0),\,\,\,0\leq j_0<q_{n+1}.$
First we assume $j_0=0,$ i.e. $I^{(n)}_{j_0}= I^{(n)}_{0}$. Applying Denjoy inequality  (\ref{DenInq})
we see that the intervals $T^{q_0}(I^{(n)}_0),\, T^{q_1}(I^{(n)}_0),\,...,\, T^{q_{n}}(I^{(n)}_0)$ of the partition  $P_{n}$
 are $e^v$-comparable with $I^{(n)}_{0}.$\\
If $0< j_0<q_{n+1}$ then   $q_{i_0}\leq j_0\leq q_{i_0+1},  $   for some $0\leq i_{0}\leq { n}.$\\
In this case, the intervals
 \begin{equation} \label{Comp.Int} T^{-q_{i_0}}I^{(n)}_{j_0},\,T^{-q_{i_0-1}}I^{(n)}_{j_0},\,...,\,
T^{-q_{0}}I^{(n)}_{j_0},\,T^{q_{0}}I^{(n)}_{j_0},\,...\,T^{q_{n-i_0-1}}I^{(n)}_{j_0}
 \end{equation}
 are elements of the partition  $P_n(T,z_0)$.  Applying  again  Denjoy's inequality (\ref{DenInq}) we can see
 that each interval in  (\ref{Comp.Int}) is  $e^v-$comparable with $I^{(n)}_{j_0}.$\\
For intervals  $I^{(n+1)}_{j_0},\,0\leq j_0<q_{n}$, the statement can be proved  analogously.
\endproof

We recall the following definition introduced in \cite{KO1989}.

\begin{defi}\label{SMALL} An interval $I=[\tau,t]\subset S^{1}$ is said to
be $q_{n}$-small, and its endpoints $q_{n}$-close, if the intervals
$T^{i}(I),\; 0\le i\leq q_{n}-1$, are, except for the endpoints, pairwise disjoint.
\end{defi}
It follows from the structure of the dynamical partition, that an
interval $I=[\tau,t]$ is $q_{n}$-small if and only if either
$\tau\prec t \preceq T^{q_{n-1}}(\tau)$ or $T^{q_{n-1}}(t)\preceq
\tau\prec t$.
\begin{lemm}\label{FIN}(see \cite{Fi1950}). Suppose, the homeomorphism  $T$ with irrational rotation
number $\rho_T$  satisfies the conditions of Lemma \ref{DEN} and
 the interval $I=(x,y)\subset S^{1}$ is $q_{n}$-
small. Then for any $0\leq k <q_{n}$ Finzi's inequality holds

\begin{equation}\label{FLK}
e^{-v}\leq \frac{ DT^{k}(x)}{DT^{k}(y)}\leq e^{v},
\end{equation}
 where $v$ is the total variation of $\log DT$ on $S^{1}$.
\end{lemm}

Next we formulate the thermodynamic formalism for circle maps with a break point.

Let $X_{br}$ be the set of strictly increasing pairs of functions ($f(x),\,\, x\in [-1,0],\,\, g(x),\,\, x\in[0,\alpha]$ for some $\alpha>0$) satisfying the following conditions:
\begin{itemize}
  \item $f(0)=\alpha,\,\, g(0)=-1$;
  \item $f(-1)=g(\alpha)$;
  \item  $f(g(0))=f(-1)<0$;
  \item $f^{(2)}(g(0))\geq 0$;
  \item $f(x)\in C^{2+\varepsilon}([-1,0])$, $g(x)\in C^{2+\varepsilon}([0,\alpha])$ for all $\varepsilon>0$.
  \item $f'_+(0)\not=g'_-(0)$.
\end{itemize}
These conditions allow us to construct a circle homeomorphism $G_{f,g}$ on $[-1,\alpha)$ from a pair $(f,g)\in X_{br}$  with break point $x_b =0$ (and possibly with a second break point  $x_b=-1$ if $f'(-1) \not= g'(\alpha)$) as follows
$$G_{f,g}(x)=\left\{\begin{array}{ll}
f(x)\,\,\,\,\ \textrm{if} \,\,\,\,\,x\in [-1,0),\\
g(x)\,\,\,\,\ \textrm{if} \,\,\,\,\,x\in [0,\alpha).\\
\end{array}\right.$$

Using the map $l:[-1,\alpha] \to S^1$ with $l(x)=\frac{x+1}{\alpha+1}$ we get a circle homeomorphism $l\circ G_{f,g}\circ l^{-1}$ on $S^1=\mathbb{R}\mod 1$, which we denote   for simplicity also by $G_{f,g}$ whenever its domain of definition is clear. We define the rotation number $\rho(G_{f,g})$ of $G_{f,g}$ by the rotation number of this circle homeomorphism when acting on $S^1$  .
Denote by $X_{br}(\omega)$ the subset of $(f,g)\in X_{br}$ with  $\rho(G_{f,g})=\omega =\frac{\sqrt{5}-1}{2}$ the golden mean.
Recall the jump ration $c=\sqrt{\frac{DG_{f,g}(0_-)}{DG_{f,g}(0+)}}=\sqrt{\frac{Df(0_-)}{Dg(0+)}}$  of $G_{f,g}$ at its break point $x_b=0$. \\
It is clear that in case $ c=1 $ the homeomorphism $G_{f,g}$ is  a smooth map as long as $0$ is its only break point. We assume that $ c \neq 1 $.
Define a renormalization operator $R_{br}:X_{br}(\omega)\rightarrow X_{br}(\omega)$ as follows:
$$R_{br}(f(x),g(x))=(\tilde{f}(x),\,\,\, x\in[-1,0]; \,\,\, \tilde{g}(x),\,\,\, x\in[0,\alpha']),$$
where
$$\tilde{f}(x)=-\alpha^{-1}f(g(-\alpha x)), \,\,\, \tilde{g}(x)=-\alpha^{-1}f(-\alpha x),\,\,\, \alpha'=-\alpha^{-1}f(-1).$$
Since $D\tilde{f}(0_-)=Df(-1_+)Dg(0_+)$ and $D\tilde{g}(0_+)=Df(0_-)$ we find $\tilde{c}(0)=\sqrt{\frac{DG_{\tilde{f},\tilde{g}}(0_-)}{DG_{\tilde{f},\tilde{g}}(0_+)}}=c^{-1}\sqrt{Df(-1_+)}$. On the other hand $D\tilde{f}(-1_+)=Df(f(-1)) Dg(\alpha_-)$ and $D\tilde{g}(\alpha'_-)=Df(f(-1))$
which leads to $\tilde{c}(-1)=\sqrt{\frac{DG_{\tilde{f},\tilde{g}}(\alpha'_-)} {DG_{\tilde{f},\tilde{g}}(-1_+)}}=\frac{1}{\sqrt{Dg(\alpha_-)}}$ which in general is also different from $1$. This shows that $G_{\tilde{f},\tilde{g}}$ has in general two break points , one at $x=x_b=0$ and one at $x=G_{\tilde{f},\tilde{g}}x_b$. For the product of these two jump ratios one finds $\tilde{c}(0) \tilde{c}(-1)=\sqrt\frac{{Df(-1_+)}}{{Dg(\alpha_-)}} c^{-1}$.
From the work of Khanin and Vul  in \cite{KV1991} it is known,  that $R_{br}:X_{br}(\omega)\rightarrow X_{br}(\omega)$ has an unique periodic orbit $\{f_i(x,c_i), g_i(x,c_i), i=1,2\}$ of period two, that means
$$R_{br}(f_1(x,c_1), g_1(x,c_1))=(f_2(x,c_2), g_2(x,c_2)),$$
$$R_{br}(f_2(x,c_2), g_2(x,c_2))=(f_1(x,c_1), g_1(x,c_1)).$$

where the functions $ f_i (x, c_i) $ and $ g_i (x, c_i) $, $ i = 1,2 $ have the following explicit form:
\begin{gather} \label{4}
f_i(x,c_i)=\frac{(\alpha_i+c_ix)\beta_i}{\beta_i+(\beta_i+\alpha_i-c_i)x},
\end{gather}
\begin{gather} \label{4'}
g_i(x,c_i)=\frac{\alpha_i\beta_i(x_i-c_i)}{\alpha_i\beta_ic_i+(c_i-\alpha_i-c_i\beta_i)x},
\end{gather}
with
$$\alpha_1=\frac{c-\beta_{0}^{2}}{1+\beta_{0}},\,\, \alpha_2=\frac{c^{-1}-\beta_{0}^{2}}{1+\beta_{0}},\,\, c_1=c,\,\, c_2=c^{-1},\,\, \beta_1=\beta_2=\beta_0,$$
 $\beta_0$  the unique  root of the equation
$$\beta^4-\beta^3-\beta^2\frac{(c+1)^2}{c}-\beta+1=0,$$
belonging to the interval $(0,1)$.

By using the pairs of functions $(f_i, g_i),\, i=1,2$ we define  circle homeomorphisms $G_i: [-1,\alpha_i]\to [-1,\alpha_i] ,\,i=1,2$, as

$$G_{i}(x)=\left\{\begin{array}{ll}
f_i(x,c_i)\,\,\,\,\ \textrm{if} \,\,\,\,\,x\in [-1,0),\\
g_i(x,c_i)\,\,\,\,\ \textrm{if} \,\,\,\,\,x\in [0,\alpha_i). \\
\end{array}\right.$$
where we used the fact that $f_i(0)=\alpha_i$ and $g_i(0)=-1$ respectively $f_i(-1)=g_i(\alpha_i)=-\beta_i$. Since $Df_1(0_-)=\frac{c(\alpha_1+\beta_0)-\alpha_1(\alpha_{1}+\beta_0)}{\beta_0}$and $Dg_1(0_+)=\frac{(1-\beta_0)(c-\alpha_1)}{\alpha_1\beta_0 c}$ one finds $c_{G_1}(0)=\sqrt{\frac{Df_1(0_-)}{Dg_1(0_+)}}=\sqrt{\frac{(c+\beta_0)(c-\beta_0^2)}{(1+\beta_0)(1-\beta_0^2)}}$. But this equals $1$ only iff  $c=1$. Hence $G_1$ has a break point at $x_b=0$. On the other hand one finds $Df_1(-1_+)=\frac{\beta_0(\alpha_{1}+\beta_0)}{c-\alpha_1}$
respectively $Dg_1(\alpha_1)=\frac{c\beta_0(1-\beta_0)}{\alpha_1(c-\alpha_1)}$ which leads to $c_{G_1}(-1)=\sqrt{\frac{Dg_1(\alpha_{1_-)}}{Df_1(-1_+)}}=\sqrt{\frac{(1+\beta_0)(1-\beta_0^2)c}{(c+\beta_0)(c-\beta_0^2)}}$. This equals $1$ iff $c=1$. This  shows that the map $G_1$ has indeed two break points, namely at $x=x_b=0$ and at $x=G_1(x_b)=-1$. For the product of the jump ratios one finds $c_{G_1}(0) c_{G_1}(G_1(0)) = \sqrt{c}$

Using the function $l_i(x)=\frac{1+x}{1+\alpha_i}$ the map $l\circ G_i\circ l^{-1}$ defines then a homeomorphism of the circle $S^1$ with breakpoints $x_b=\frac{1}{\alpha_i+1}$ and $x'_b=0$.

We rename the homeomorphism of the circle $S^1$ corresponding to $G_1$ by $G_{br}$ and denote by $B(G_{br})$ the set of all circle homeomorphisms which are $C^{1+\epsilon}$ conjugate to $G_{br}$.

The sequence of dynamical partitions $P_{n}(x_b), n\geq 1$ allows us to introduce a symbolic dynamics for the map $G_{br}$.
For this take an arbitrary point $x\in S^1\setminus O^+_{G_{br}}(x_b)$ where $O^+_{G_{br}}(x_b)$ denotes the forward orbit of the break point $x_b$ of $G_{br}$.
For $n\geq 0$ put $a_{n+1}=: a_{n+1}(x) = a$ if $x\in I_j^{(n+1)}(x_b), 0\leq j<q_n$.  If however $x\in I_i^{(n)}(x_b), 0\leq i< q_{n+1}$ we
know from the construction of the partition $P_{n+1}$ from $P_{n}$ in case  $\rho(G_{br})=\omega=\frac{\sqrt{5}-1}{2}$ that either $x\in
I_i^{(n+2)}(x_b), 0\leq i<q_{n+1}$ or $x\in I_{i+q_n}^{(n+1)}(x_b), 0\leq i<q_{n+1}$. In the first case we put $a_{n+1}=0$ and in the second $a_{n+1}=1$. This way we get a one-
to-one correspondence \newline
$$\varphi: S^1\setminus O^+_{G_{br}}(x_b)\leftrightarrow \{(a_1,\cdots ,a_n,\cdots), a_n\in \{a,0,1\}; \, a_{n+1}=a \iff a_n=0,\, n\geq 1\}:=\Theta_{\mathbb{A}},$$
where $\Theta_{\mathbb{A}}$ denotes the space of allowed infinite one-sided sequences of symbols from the alphabet $A=\{1,0,a\}$ with transition matrix    $$\mathbb{A}=\begin{pmatrix} 1&1&0\\ 0&0&1\\1&1&0\end{pmatrix},$$
that is
    $$\Theta_{\mathbb{A}}:=\{\underline{a}=(a_{1},\ldots,a_{n},\ldots),\,\,\,a_i\in A, \mathbb{A}_{a_i,a_{i+1}}=1\;\text{for}\; i\in \mathbb{Z}_+\}.$$
The triple $(\theta_{\mathbb{A}}, {A},\sigma)$ with $\sigma:\theta_{\mathbb{A}}\to\theta_{\mathbb{A}}$ the shift map $(\sigma(\underline{a}))_i=a_{i+1}$ is called a subshift of finite type over the alphabet $A$ with transition matrix $\mathbb{A}$.
Notice that every interval $I^{(n)}$ of the dynamical partition $P_n$ corresponds to the unique  finite word $(a_1,\,a_2,\,...,a_n)$ of  length $n$. In particular, for $n$ even  the words $(a,\,0,\,a,\,0,\,...,\,a,\,0)$ and $(0,\,a,\,0,\,a,\,...,\,0,\,a)$  correspond to the atoms $I^{(n+1)}_0(x_b)$ respectively  $I^{(n)}_0(x_b)$ of  $P_n$. \\

In \cite{DjK} A. Dzhalilov and J. Karimov  constructed  the thermodynamic formalism for maps in $B(G_{br})$ by using a closely related subshift of finite type as follows:\newline
Denote by $\mathbb{A}^t=(\mathbb{A}^t_{i,j}), i,j\in A$ the transposed matrix of $\mathbb{A}$ and by $
(\Theta_{\mathbb{A}^t}, {A}, \sigma)$ the subshift of finite type with alphabet $A$ and transition matrix $\mathbb{A}
^t$. Obviously $\underline{b}=(b_1,\cdots, b_n,b_{n+1, \cdots})\in \Theta_{\mathbb{A}^t}$ if and only if  $b_{n+1}=0$ iff $b_n=a$. Therefore a finite sequence  $(a_1,\cdots,a_n)$ is a subsequence of some $\underline{a}\in
\Theta_{\mathbb{A}}$ iff $(a_n, a_{n-1},\cdots,a_1) $ is a subsequnce of some $\underline{b}\in \Theta_{\mathbb{A}
^t}$. Since $\mathbb{A}^3_{i,j}={(\mathbb{A}^t)}^3_{j,i}>0$ for all $i,j\in A$, the map $\sigma$ is topological mixing
in both $\Theta_{\mathbb{A}} $ and $\Theta_{\mathbb{A}^t}$.

In the following we denote by $\vec{a}$ the vector $(a_1,\ldots,a_n)\in A^n: \mathbb{A}_{a_i,a_{i+1}}=1, 1\leq i
\leq n-1$, and by $\overleftarrow{b}=(b_1,\cdots, b_n): \mathbb{A}^t_{b_i,b_{i+1}}=1, 1\leq i\leq n-1$.

Define for $x\in A$
$$\underline{\gamma}(x):=\left\{\begin{array}{ll}
(a,0,a,0,\ldots,a,0,\ldots),\,\,\,\,\,\,\,\,\,x=0,1,\\
(0,a,0,a,\ldots,0,a,\ldots),\,\,\,\,\,\,\,\,\,x=a. \\
\end{array}\right.$$

\begin{theo}\label{thermo} (see \cite{DjK})  For any $G\in B(G_{br})$, there exists a function  $U_{br}:\Theta_{\mathbb{A}^t}\rightarrow (-\infty, 0)$ continuous in the product topology, such that the following properties hold:\\
{\rm 1)} For any $\underline{b}=(b_1,\ldots,b_{k},b_{k+1},\ldots,b_{n},\ldots)$,
$\underline{c}=(b_1,\ldots,b_{k},c_{k+1},\ldots,c_{n},\ldots)\in
\Theta_{\mathbb{A}^t}$ there exist  constants  $C_{1}>0$ and   $q\in(0,1),$   not depending on $\underline{b}$, $\underline{c}$ and $k$, such that
$$|U_{br}(\underline{b})-U_{br}(\underline{c})|\leq C_{1}\cdot q^k.$$\\
{\rm 2)} Let $\Delta_{s_n}^{(n)}\subset\Delta_{s_r}^{(r)}$, $1\leq r<n,\,0 \leq s_{r}\leq q_{r+1}-1,\,\,0 \leq s_{n}\leq q_{n+1}-1   $ and $\varphi(\Delta_{s_n}^{(n)})=(a_{1},\ldots ,a_{r},\ldots, a_{n})$, $\varphi(\Delta_{s_r}^{(r)})=(a_{1}\ldots,a_{r}),$ then
$$|\Delta_{s_n}^{(n)}|=(1+\psi_{r}{(a_1,\ldots,a_{n}}))|
\Delta_{s_r}^{(r)}|\exp\{\sum\limits_{s=r}^{n}U_{br}(a_s,a_{s-1},\ldots,a_r,\ldots,a_1,\underline{\gamma}(a_{1}))\},$$
where $|\psi_{r}{(a_1,\ldots,a_{n}})|\leq C_{2}\cdot q^{r}$, with the constant $C_{2}>0$   not depending on
$ r,n$ and  $(a_1,\ldots,a_{n}). $
\end{theo}
Here $\Delta^{(n)}$ are elements of the dynamical partition  $P_n(G,x_b):=\{  \Delta_{i}^{(n)}= G^i\Delta_{0}^{(n)},\,0\leq i\leq q_{n+1}-1\} \bigcup \{\Delta_{j}^{(n+1)}= G^j\Delta_{0}^{(n+1)},\,0\leq j\leq q_{n}-1 \}$ with $\Delta_{0}^{(n)}=[x_{q_n}, x_b)$ respectively $\Delta_{0}^{(n+1)}=[ x_b, x_{q_{n-1}}) $ and $x_i=G^i x_b$, where $x_b=0$  denotes the break point of $G$.
The function $U_{br}$ is defined as follows: for $\underline{b}=(b_1,\cdots,b_k,b_{k+1,}\cdots)\in \Theta_{\mathbb{A}^t}$ and  any $k\geq 1$ denote for $i=1,2$ by $\vec{a}^i_k$ the finite sequence  $\vec{a}^i_k=(b_k,\cdots, b_i)$. Choose $I(\vec{a}^i_k)\in P_{k-i+1}$ such that $\varphi(I(\vec{a}^i_k))=\vec{a}^i_k$. Obviously one has $I(\vec{a}^1_k)\subset I(\vec{a}^2_k)$. Define for $\overleftarrow{b_k}=(b_1,\cdots,b_k)$ the function $U_k=U_k(\overleftarrow{b}_k)$ as
$$U_k(\overleftarrow{b_k}):=\log \frac{|I(\vec{a}^1_k)|}{|I(\vec{a}^2_k)|}.$$
One can show that the thermodynamic limit $k\to \infty$ of the function $U_k$ exists, namely for $\underline{b}\in \Theta_{\mathbb{A}^t}$
$$U_{br}(\underline{b})=\lim\limits_{k\to \infty} U_k(\vec{b}_k)$$
is a well defined function on $\Theta_{\mathbb{A}^t}$ with values in $(-\infty,0)$.

According to a well known theorem in \cite{Ruelle} (Proposition 5.13), the following limit exists
\begin{equation}\label{ruile}\ln \lambda_{{\beta}}=\lim\limits_{n\to\infty}\frac{1}{n}\ln\left[\sum
\limits_{\vec{\epsilon}=(\epsilon_1,\cdots,\epsilon_n)}\exp\left\{
{\beta} \sum\limits_{s=1}^{n}
U_{br}(\varepsilon_s,\varepsilon_{s-1},\ldots,\varepsilon_1,\underline{\gamma}(\varepsilon_1))\right\}\right],
\end{equation}
where  $\lambda_{\beta}$ is the leading eigenvalue of the transfer operators  $D_{\beta}: C(\Theta_{\mathbb{A}^t})\rightarrow C(\Theta_{\mathbb{A}^t})$ respectively its dual $D^{*}_{\beta}:M(\Theta_{\mathbb{A}}^t)\rightarrow M(\Theta_{\mathbb{A}^t})$
\begin{equation}\label{transfer1}(D_{\beta}f)(\underline{b})=\sum\limits_{\underline{x}\in \sigma^{-1}(\underline{b})\cap \Theta_{\mathbb{A}^t}}e^{\beta U_{br}(\underline{x})}f(\underline{x})\end{equation}
\begin{equation}\label{transfer2}(D^{*}_{\beta}\mu)(c ,\underline{b})=\mathbb{A}^t_{{c,b_1}}e^{\beta U_{br}({c},\underline{b})}\mu(\underline{b}).\end{equation}
$\sigma: \Theta_{\mathcal{A}^t}\rightarrow \Theta_{\mathcal{A}^t}$ is the shift map $\sigma(\underline{b})_i=b_{i+1}, i\geq 1$, $C(\Theta_{\mathcal{A}^t})$ denotes the space of continuous functions and $M(\Theta_{\mathcal{A}^t})$  the  space of Borel measures on $\Theta_{\mathbb{A}^t}$.  We used also the formal expression $d\mu(x)=\mu(x) dx$ for the measure $d\mu$.\\

\section{Estimates for the Lyapunov function}
 In this section we estimate  the Lyapunov's functions.
Let $T\in C^{2+\varepsilon}(S^{1}\setminus \{x_{b}\}),,\,\,\,\varepsilon>0$, be a circle
homeomorphism with  break point $x_{b}$, $T'(x)\geq Const >0,\,x \in S_1\equiv
[x_{b},x_{b}+1) $  and rotation number $\rho_T=[k_{1},k_{2},..,k_{m},1,1,...]\,\text{for
some} \, m\geq 1$.
In this case obiously $q_{n+1}=q_{n}+q_{n-1},\,n\geq m$.
Consider the $n$-th dynamical partition
$P_{n}(x_b)=\{\, I^{(n)}_{j}(x_{b}),\,\,0 \leq j<q_{n+1};\,\,\,I^{(n+1)}_{i}(x_{b}),\,0\leq i<q_{n}\}$
 determined by the break point $x_{b}$ under the map $T$. We denote by $I^{(n)}({x_{b}};z_{0})  $  the interval of
  the partition $P_{n}(x_{b})$ containing the point $z_{0}$.
\\
\begin{theo}\label{Lyap1} For any $\beta>0$ and $z_0\in S^{1}\setminus \{T^{i}(x_b),\,i=0, -1,-2,...,-q_n+1\},$
there exists $n_0:=n_0(\beta)> m$, such that for every $ n>n_0$  the  following inequalities hold
\begin{equation}\label{ln}
c_{1}|I^{(n)}({x_{b}};z_{0})|^\beta \lambda_{-\beta}^n\leq \Lambda_\beta(z_0,q_{n})\leq C_{1} |I^{(n)}({x_{b}};z_{0})|^\beta \lambda_{-\beta}^n,
\end{equation}
where  $\lambda_{-\beta}$ is the leading eigenvalue of the transfer operator $D_{-\beta}$ in (\ref{transfer1}) and the positive  constants $c_{1},\,C_{1}$ depend on the map $T$.
\end{theo}
\proof
Using the elements of the partition $P_{n}({x_{b}}),$
we introduce the following sum
$$S_{n,\beta}(z_0)=
|I^{(n)}_{0}({x_{b}};z_{0})|^\beta \left(\sum\limits_{k=0}^{q_{n+1}-1}|I^{(n)}_{k}({x_{b}})|^{-\beta}+
\sum\limits_{s=0}^{q_{n}-1}|I^{(n+1)}_{s}({x_{b}})|^{-\beta}\right).$$
Then one has the following
\begin{lemm}\label{3.1} For any $z_0\in S^{1}\setminus \{T^{i}(x_b),\,i=0, -1,-2,...,-q_{n}-q_{n+1}+1\},$ $\beta>0$ and $n\geq1$
\begin{equation}\label{3}c_2\cdot S_{n,\beta}(z_0)\leq\Lambda_\beta (z_0,q_{n}+q_{n+1})\leq C_2\cdot S_{n,\beta}(z_0),\end{equation}
where the positive constants $c_2$ and $C_2$ depend only on  the total variation of  $\log T'$ and $\beta$.
\end{lemm}
\proof
For any $z_0\in S^1\setminus \{T^{-i}x_b,\, 0\leq i< q_{n}+q_{n+1}\}$ and $\beta>0,$ we consider
\\
\begin{equation}
\Lambda_{\beta}(z_0,q_{n}+q_{n+1})=1+\sum\limits_{k=1}^{q_{n}+q_{n+1}-1}\prod\limits_{j=k}^{q_{n}+q_{n+1}-1}\left|DT(z_j)\right|^{\beta}.
\end{equation}
We rewrite $\Lambda_{\beta}(z_0,q_{n}+q_{n+1})$ in the following form:
\begin{equation}\label{1}
\Lambda_{\beta}(z_0,q_{n}+q_{n+1})=\left|DT^{q_{n}+q_{n+1}-1}(z_1)\right|^{\beta}\left(1+\sum\limits_{k=1}^{q_{n}+q_{n+1}-1}\left|DT^{k}(z_1)\right|^{-\beta}\right).
\end{equation}
Since
$$DT^{q_{n}+q_{n+1}-1}(z_1)=\frac{DT^{q_{n}+q_{n+1}}(z_0)}{DT(z_0)},$$
and using the bounds
$$K_1^{-1}\leq DT(x)\leq K_1,\,\,x\in[x_b,x_b+1],\, K_1=K_1(T)>0,$$
respectively Denjoy's inequality (\ref{FLK}) we obtain
\begin{equation}\label{bda}
K_1^{-1}e^{-2v}\leq DT^{q_{n}+q_{n+1}-1}(z_1)\leq K_1 e^{2v},
\end{equation}
where $v=\underset{x\in S^1}{\mathrm{var}} DT(x)+2|\ln c_T(x_b)|$.
Next we estimate the sum in (\ref{1}). Let $P_n(T,x_b)$ and $P_n(T,z_0)$ be the $n$-th dynamical partitions determined by   $x_b$ respectively $z_0$ under the map $T$.
\\
Here two cases are possible
\\
\textbf{Case I.} $z_0\in I_{i_0}^{(n)}(x_b),\,0\leq i_0<q_{n+1}$,\\
\\
\textbf{Case II.} $z_0\in I_{j_0}^{(n+1)}(x_b),\,0\leq j_0<q_{n}$.\\
\\
We consider only case I., case II. can be treated  similarly.
\\
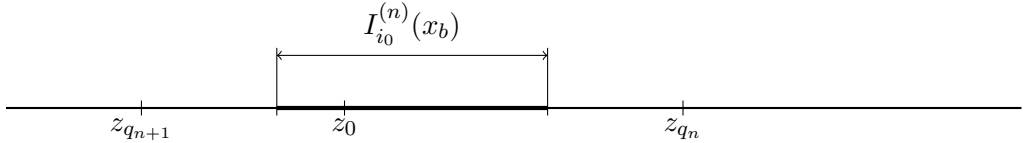
\begin {tikzpicture} [xscale=0.9]
\\
\draw [thick](0,0) -- (15,0);
\draw (2,0.1) -- (2,-0.1);
\draw (5,0.1) -- (5,-0.1);
\draw (10,0.1) -- (10,-0.1);
\draw (4,0.8) -- (4,-0.1);
\draw (8,0.8) -- (8,-0.1);
\draw[<->] (4,0.7) -- (8,0.7);
\draw [ultra thick](4,0) -- (8,-0);
\node[below] at (5,0){$z_{0}$};
\node[below] at (2,0){$z_{q_{n+1}}$};
\node[below] at (10,0){$z_{q_{n}}$};
\node[above] at (6,0.7) {$I^{(n)}_{i_0}(x_b)$};
\end{tikzpicture}
\\
Obviously
\begin{equation}
\begin{split}
1+\sum\limits_{k=1}^{q_{n+1}+q_n-1}\left|DT^{k}(z_1)\right|^{-\beta}&=|DT(z_0)|^{\beta}\left(\sum\limits_{k=1}^{q_{n+1}}|DT^k(z_0)|^{-\beta}+\sum\limits_{k=q_{n+1}+1}^{q_{n+1}+q_n}|DT^k(z_0)|^{-\beta}\right)\\
&=|DT(z_0)|^{\beta}(S_n^{(1)}(\beta)+S_n^{(2)}(\beta)).
\end{split}
\end{equation}
First we estimate the sum $S_n^{(1)}(\beta)$.
\\
By assumption $z_0\in I_{i_0}^{(n)}(x_b),\,0\leq i_0<q_{n+1}$.
In case $0\leq i_0\leq q_{n+1}-2$ we split the sum $S^{(1)}_n(\beta)$ into two parts
\begin{equation}
S^{(1)}_n(\beta)=\sum\limits_{k=1}^{q_{n+1}-i_0-1}|DT^k(z_0)|^{-\beta}+\sum\limits_{k=q_{n+1}-i_0}^{q_{n+1}}|DT^k(z_0)|^{-\beta}.
\end{equation}\newline
In case $i_0=q_{n+1}-1$ we put
\begin{equation} \label{estDT32}
S_n^{(1)}(\beta)=\sum\limits_{k=1}^{q_{n+1}}|DT^k(z_0)|^{-\beta}.
\end{equation}
If $0\leq i_0\leq q_{n+1}-2$ we estimate $DT^k(z_0)$ for $1\leq k\leq q_{n+1}-i_0-1$ and $q_{n+1}-i_0\leq k\leq q_{n+1}$ separately.
\\
If $1\leq k\leq q_{n+1}-i_0-1$ we have
$$I_{k+i_0}^{(n)}(x_b)=\int\limits_{I_{i_0}^{(n)}(x_b)}DT^k(t)dt.$$
It is clear that
$$\frac{I_{k+i_0}^{(n)}(x_b)}{DT^k(z_0)}=\int\limits_{I_{i_0}^{(n)}(x_b)}\frac{DT^k(t)}{DT^k(z_0)}dt.$$
Applying  Finzi's inequality (\ref{FLK}) we obtain:
\begin{equation}\label{estDT}
e^{-v}\frac{|I_{i_0}^{(n)}(x_b)|}{|I_{k+i_0}^{(n)}(x_b)|}\leq (DT^k(z_0))^{-1}\leq e^{v}\frac{|I_{i_0}^{(n)}(x_b)|}{|I_{k+i_0}^{(n)}(x_b)|}.
\end{equation}
These inequalities imply
\begin{equation}\label{estDT1a}
e^{-\beta v}|I_{i_0}^{(n)}(x_b)|^{\beta}\sum\limits_{k=1}^{q_{n+1}-i_0-1}|I_{k+i_0}^{(n)}(x_b)|^{-\beta}\leq
\sum\limits_{k=1}^{q_{n+1}-i_0-1}|DT^k(z_0)|^{-\beta}\leq
e^{\beta v}|I_{i_0}^{(n)}(x_b)|^{\beta}\sum\limits_{k=1}^{q_{n+1}-i_0-1}|I_{k+i_0}^{(n)}(x_b)|^{-\beta}.
\end{equation}
Next consider the case  $q_{n+1}-i_0< k\leq q_{n+1}$. Then obviously
\begin{equation}\label{estDT1}
DT^k(z_0)=DT^{q_{n+1}-i_0}(z_0)DT^{k-(q_{n+1}-i_0)}(z_{q_{n+1}-i_0}).
\end{equation}
As for the estimate (\ref{estDT}) we can show that
$$
e^{-v}\frac{|I_{q_{n+1}}^{(n)}(x_b)|}{|I_{k+i_0}^{(n)}(x_b)|}\leq (DT^{k-(q_{n+1}-i_0)}(z_{q_{n+1}-i_0}))^{-1}\leq e^{v}\frac{|I_{q_{n+1}}^{(n)}(x_b)|}{|I_{k+i_0}^{(n)}(x_b)|}.
$$

Since
$$e^{-v}|I_{k-(q_{n+1}-i_0)}^{(n)}(x_b)|\leq |I_{k+i_0}^{(n)}(x_b)|=\int\limits_{I_{k-(q_{n+1}-i_0)}^{(n)}(x_b)}DT^{q_{n+1}}(t)dt\leq e^v|I_{k-(q_{n+1}-i_0)}^{(n)}(x_b)|$$
one finds
$$
e^{-2v}\frac{|I_{q_{n+1}}^{(n)}(x_b)|}{|I_{k-(q_{n+1}-i_0)}^{(n)}(x_b)|}\leq (DT^{k-(q_{n+1}-i_0)}(z_{q_{n+1}-i_0}))^{-1}\leq e^{2v}\frac{|I_{q_{n+1}}^{(n)}(x_b)|}{|I_{k-(q_{n+1}-i_0)}^{(n)}(x_b)|}.
$$

Inserting $k=q_{n+1}-i_0$ into (\ref{estDT})  we get
$$
e^{-v}\frac{|I_{i_0}^{(n)}(x_b)|}{|I_{q_{n+1}}^{(n)}(x_b)|}\leq (DT^{q_{n+1}-i_0}(z_0))^{-1}\leq e^{v}\frac{|I_{i_0}^{(n)}(x_b)|}{|I_{q_{n+1}}^{(n)}(x_b)|}.
$$
The last   inequalities together with (\ref{estDT1}) imply for $q_{n+1}-i_0\leq k\leq q_{n+1}$
\begin{equation}\label{estDT2a}
e^{-3v}\frac{|I_{i_0}^{(n)}(x_b)|}{|I_{k-(q_{n+1}-i_0)}^{(n)}(x_b)|}\leq (DT^k(z_0))^{-1}\leq e^{3v}\frac{|I_{i_0}^{(n)}(x_b)|}{|I_{k-(q_{n+1}-i_0)}^{(n)}(x_b)|}.
\end{equation}
\newline
In  case $i_0=q_{n+1}-1$ one shows for $1\leq k\leq q_{n+1}$  as above the bounds
$$
e^{-2v}\frac{|I_{q_{n+1}-1}^{(n)}(x_b)|}{|I_{k-1}^{(n)}(x_b)|}\leq (DT^k(z_0))^{-1}\leq e^{2v}\frac{|I_{q_{n+1}-1}^{(n)}(x_b)|}{|I_{k-1}^{(n)}(x_b)|}.
$$

From (\ref{estDT1a}), (\ref{estDT2a}) and the  last inequality we get the bounds for $S_n^{(1)}(\beta)$.
Similarly we can derive the following estimate:
$$e^{-3v\beta}|I_{i_0}^{(n)(x_b)}|^{\beta}\sum\limits_{s=0}^{q_n-1}|I_{s}^{(n+1)}(x_b)|^{-\beta}\leq S_n^{(2)}(\beta)\leq e^{3v\beta}|I_{i_0}^{(n)(x_b)}|^{\beta}\sum\limits_{s=0}^{q_n-1}|I_{s}^{(n+1)}(x_b)|^{-\beta}.$$
 Using (\ref{bda}) and  the estimates for $S_n^{(1)}(\beta)$ respectively the ones for $S_n^{(2)}(\beta)$  we get
 the assertion of Lemma \ref{3.1}\newline
\endproof

Next we estimate the sum $S_{n,\beta}(z_0)$.

Let $n>m>0$. Consider the $m$-th renormalization interval $J_m(x_b)=[x_{q_{m+1}},x_{q_m})=I^{(m+1)}_0(x_b)\cup I^{(m)}_0(x_b).$
\\
Define
$$\tilde{\tau}_{n-m}(x_b):=P_n(T,x_b)\cap J_m(x_b),$$
$$\tilde{\tau}_{n-m}^{(m-1)}(x_b):=P_n(T,x_b)\cap [x_{q_{m+1}}, x_b)$$
$$\tilde{\tau}_{n-m}^{(m)}(x_b):=P_n(T,x_b)\cap [x_b, x_{q_{m}})$$

Obviously $\tilde{\tau}_{n-m}(x_b)$ is a partition of $J_m(x_b)$. We will compare $S_{n,\beta}(z_0)$ with the sum of the lengths of the intervals of this partition $\tilde{\tau}_{n-m}(x_b)$.
\begin{equation}\label{uzunten}
\begin{split}S_{n,\beta}(z_0)=
&|I_0^{(n)}( x_b ; z_0)|^{\beta}\left(\sum\limits_{I_{i_k}^{(n)}(x_{b})
\in \tilde{\tau}^{(m+1)}_{n-m}} |I^{(n)}_{i_k}(x_{b})|^{-\beta}+\sum\limits_{I_{j_k}^{(n+1)}(x_{b})\in \tilde{\tau}_{n-m}^{(m+1)}} |I^{(n+1)}_{j_k}(x_{b})|^{-\beta}\right)\\
&+|I_0^{(n)}( x_b ; z_0)|^{\beta} \left(\sum\limits_{I_{i_k}^{(n)}(x_{b})
\in \tilde{\tau}^{(m)}_{n-m}} |I^{(n)}_{i_k}(x_{b})|^{-\beta}+\sum\limits_{I_{j_k}^{(n+1)}(x_{b})\in \tilde{\tau}_{n-m}^{(m)}} |I^{(n+1)}_{j_k}(x_{b})|^{-\beta}\right)\\
&+|I_0^{(n)}( x_b ; z_0)|^{\beta}\sum\limits_{s=1}^{q_{m}-1}\left(\sum\limits_{I_{i_k}^{(n)}(x_{b})
\in \tilde{\tau}^{(m+1)}_{n-m}} |T^s(I^{(n)}_{i_k}(x_{b}))|^{-\beta}+\sum\limits_{I_{j_k}^{(n+1)}(x_{b})\in \tilde{\tau}_{n-m}^{(m+1)}} |T^s(I^{(n+1)}_{j_k}(x_{b}))|^{-\beta}\right)\\
&+|I_0^{(n)}( x_b ; z_0)|^{\beta} \sum\limits_{s=1}^{q_{m+1}-1}\left(\sum\limits_{I_{i_k}^{(n)}(x_{b})
\in \tilde{\tau}^{(m)}_{n-m}} |T^s(I^{(n)}_{i_k}(x_{b}))|^{-\beta}+\sum\limits_{I_{j_k}^{(n+1)}(x_{b})\in \tilde{\tau}_{n-m}^{(m)}} |T^s(I^{(n+1)}_{j_k}(x_{b}))|^{-\beta}\right).
\end{split}
\end{equation}
\\
We denote the first two sums in (\ref{uzunten}) by  $|I_0^{(n)}( x_b ; z_0)|^{\beta} \tilde{D}_{n-m}(x_b,\beta)$ and the remaining part by $|I_0^{(n)}( x_b ; z_0)|^{\beta}\tilde{D}^{(1)}_{n-m}(x_b,\beta).$
Since $0<const\leq T'(x)\leq  Const$ for  $x\in [x_b,x_b+1]$ we have
$$0<const\leq \min\limits_{1\leq s\leq q_m} \inf\limits_{S^1} DT^{s}(x)\leq \max\limits_{1\leq s\leq q_m} \sup\limits_{S^1}DT^s(x)\leq Const.$$
Hence one finds for  some constants $c$ and $C$ depending on $n\, ,\,m$ and $z_0$
\begin{equation}\label{SDbaho}
c\cdot |I_0^{(n)}( x_b ; z_0)|^{\beta} \tilde{D}_{n-m}(x_b,\beta)\leq S_{n,\beta}(z_0) \leq C\cdot |I_0^{(n)}( x_b ; z_0)|^{\beta}\tilde{D}_{n-m}(x_b,\beta).
\end{equation}
To estimate the sum $\tilde{D}_{n-m}(x_b,\beta)$ we make  use of the thermodynamic formalism.
\\
For this consider  the first return map $\pi_m: J_0^{(m)}(x_b)\to J_0^{(m)}(x_b)$ given by
$$\pi_m(x)=\left\{
             \begin{array}{ll}
               T^{q_m}(x), & \hbox{if } x\in I^{(m+1)}_0(x_b), \\
               T^{q_{m+1}}(x), & \hbox{if } x\in I^{(m)}_0(x_b).
             \end{array}
           \right.
$$
We can pass from $[x_{q_{m+1}},x_{q_m})$ to the unit circle $S^1\equiv [0,1)$ by the affine map $z:[x_{q_{m+1}},x_{q_m})\to [0,1)$ given by
$$z=\frac{x-x_{q_{m+1}}}{x_{q_m}-x_{q_{m+1}}},\, x_{q_{m+1}}\leq x< x_{q_m}.$$
We rewrite the maps $T^{q_m}$ and $T^{q_{m+1}}$ in normalized coordinates
$$f_m(z):=\frac{T^{q_m}(x(z))-x_{q_{m+1}}}{x_{q_m}-x_{q_{m+1}}}=\frac{T^{q_m}(x_{q_{m+1}}+z(x_{q_m}-x_{q_{m+1}}))-x_{q_{m+1}}}{x_{q_m}-x_{q_{m+1}}}$$
$$g_m(z):=\frac{T^{q_{m+1}}(x(z))-x_{q_{m+1}}}{x_{q_m}-x_{q_{m+1}}}=\frac{T^{q_{m+1}}(x_{q_{m+1}}+z(x_{q_m}-x_{q_{m+1}}))-x_{q_{m+1}}}{x_{q_m}-x_{q_{m+1}}},$$
where $0\leq z <1$. \newline
With $a_m:=\frac{x_b-x_{q_{m+1}}}{x_{q_m}-x_{q_{m+1}}}$ we define the map
 $$T_m(z):=\left\{
                       \begin{array}{ll}
                         f_m(z), & \hbox{if } 0\leq z<a_m, \\
                         g_m(z), & \hbox{if } a_m\leq z<1.
                       \end{array}
                     \right.
$$
 Since $T_m(0)=f_m(0)=g_m(1)=T_m(1)$ and $T_m(a_m-)=f_m(a_m)=1$ respectively $T_m(a_m+)=g_m(a_m)=0$, the map $T_m$ defines a homeomorphism of the circle $S^1=[0,1)$. Notice that the rotation number $\rho_{T_m}$ of $T_m$ equals the rotation number $\rho_{T}$ which is the "golden mean", i.e.
$$\rho_{T_m}=\frac{\sqrt{5}-1}{2}=[1,1,1,...].$$
Since $DT_m(a_m-)=Df_m(a_m-)=DT^{q_m}(x_b-)$ and $DT_m(a_m+)=DT^{q_{m+1}}(x_b+)$ one finds for the jump ratio $c^2_{T_m}(a_m)=\frac{DT_m(a_m-)}{DT_m(a_m+)}$ of the map $T_m$   at the point
$a_m\in S^1$:  $$c^2_{T_{m}}(a_m)=c^2_T(x_b) \prod\limits_{s=0}^{q_{m-1}-1} DT(x_{q_{m}+s})^{-1}.$$ On the other hand,
since $Df_m(0+)=DT^{q_m}(x_{q_{m+1}}+)$ respectively $Dg_m(1-)=DT^{q_{m+1}}(x_{q_m}-)=DT^{q_m}
(x_{q_{m+1}}-) DT^{q_{m-1}}(x_{q_m}-) $ one finds for the jump ratio $c^2_{T_m}(0)=\frac{DT_m(1-)}
{DT_m(0+)}$ of $T_m$ at the point $z=0\equiv 1$: $$c^2_{T_m}
(0)=\prod\limits_{s=0}^{q_{m-1}-1} DT(x_{q_{m}+s}).$$ Hence we get $c^2_{T_m}(a_m) c^2_{T_m}
(0)=c^2_T(x_b)$, that means the map $T_m:S^1\to S^1$ in general has two break points, namely at
$z=a_m$ and $z=0\equiv 1$.
\\
Hence for the homeomorphism $T_m\in B(G_{br})$ the statements of Theorem \ref{thermo} hold true and and therefore to
$T_m$ corresponds the same potential $U_{\text{br}}$.
\\
Denote by $\tau_{n-m}(a_m)$ the $(n-m)$-th dynamical partition of the circle $S^1$ determined by the point $a_m$ and the map $T_m^{n-m},\,n>m$.
We denote the elements of the partition $\tau_{n-m}(a_m)$  by $\Delta_i^{(n-m+1)}$ and $\Delta_{j}^{(n-m)}$ i.e.
$$\tau_{n-m}(a_m)=\{\Delta_i^{(n-m+1)},\, 0\leq i<q_{n-m};\, \Delta_{j}^{(n-m)}, \, 0\leq j<q_{n-m+1}\}.$$
 The affine map $z:[x_{q_{m+1}},x_{q_{m}}]\to S^1$ induces a $1-1$ correspondence  $\tilde{z}:\tilde{\tau}_{n-m}(x_b)\to \tau_{n-m}(a_m)$  of the intervals $\Delta^{(n-m)}$ of the  partition  $\tilde{\tau}_{n-m}(x_b)$ and the intervals $I^{(n-m)})$ of the dynamical partition $\tau_{n-m}(a_m)$ of $S^1$ such that
$$|\Delta^{(n-m)}|=|\tilde{z }(I^{(n-m)})|=\frac{|I^{(n-m)}|}{|[x_{q_{m-1}}, x_{q_m})|}.$$
Define  $D_{n-m}(\beta)$ as the following sum
$$D_{n-m}(\beta):=\sum\limits_{i=0}^{q_{n-m}-1}|\Delta_i^{(n-m+1)}|^{-\beta}+\sum\limits_{j=0}^{q_{n-m+1}-1}|\Delta_j^{(n-m)}|^{-\beta}.$$
Then obviously
\begin{equation}\label{3*}
D_{n-m}(\beta)=\frac{1}{|[x_{q_{m-1}}, x_{q_m})|}\tilde{D}_{n-m}(x_b, \beta).
\end{equation}
Using the same arguments as for an analogous result in \cite{Dzh2004} (see  Theorem  2.4)
it can be shown that
\begin{equation}\label{4*}
\lim\limits_{n\to \infty}\frac{D_{n-m}(\beta)}{\lambda_{-\beta}^{n-m}}=r(\beta)>0
\end{equation}
Summarizing (\ref{SDbaho}), (\ref{3*}) and (\ref{4*}) we obtain
\begin{equation}
c\cdot r(\beta) |I_0^{(n)}( x_b ; z_0)|^{\beta} \lambda_{-\beta}^{n-m}\leq S_{n,\beta}(z_0)\leq C \cdot r(\beta) |I_0^{(n)}( x_b ; z_0)|^{\beta} \lambda_{-\beta}^{n-m}.
\end{equation}
The last inequalities and Lemma \ref{3.1} imply the statement of Theorem \ref{Lyap1} .
\endproof

\begin{lemm}\label{Eigen.1}
For arbitrary real numbers $\beta$ and $\delta$ with $1\leq \delta< \beta$  the following inequality  holds
$$\lambda_{-\beta}^\delta<\lambda_{-\delta}^\beta,$$
where $\lambda_{t}$ is the leading eigenvalue of the transfer operator $D_{t}$ defined in (\ref{transfer1}).
\end{lemm}
\proof
 Using the bounds  (\ref{ln}) we obtain
\begin{equation} \label{Inq.2}
  \frac{c_{1}}{C_{1}} \left(\frac{\lambda_{-\beta}^{\delta}}{\lambda_{-\delta}^{\beta}}\right)^{n-m}   \leq \frac{(\Lambda_{\beta}(x,q_n))^{\delta}}{(\Lambda_{\delta}(x,q_n))^{\beta}} \leq
 \frac{C_{1}}{c_{1}} \left(\frac{\lambda_{-\beta}^{\delta}}{\lambda_{-\delta}^{\beta}}\right)^{n-m}.
\end{equation}
We will show
$$\lim\limits_{n\to\infty}\frac{(\Lambda_{\beta}(x,q_n))^{\delta}}{(\Lambda_{\delta}(x,q_n))^{\beta}}=0.$$
 This together with (\ref{Inq.2}) proves
  $\lambda_{-\beta}^\delta<\lambda_{-\delta}^\beta$ and hence Lemma \ref{Eigen.1}. \\
Using Lemma \ref{3.1}  we have
\begin{equation} \label{Inq.3}
0 \leq \lim\limits_{n\to\infty}\frac{(\Lambda_{\beta}(x,q_n))^{\delta}}{(\Lambda_{\delta}(x,q_n))^{\beta}}\leq Const
\lim\limits_{n\to\infty}\frac{(S_{n,\beta}(x))^{\delta}}{(S_{n,\delta}(x))^{\beta}}=
\end{equation}
$$=Const \lim\limits_{n\to\infty}
\frac{|I^{(n)}_{0}(x_b;z_0)|^{\beta\delta}\left(\sum\limits_{s=0}^{q_{n+1}-1}|I^{(n)}_{s}(x_{b})|^{-\beta}+\sum\limits_{k=0}^{q_{n}-1}|I^{(n+1)}_{k}(x_{b})|^{-\beta}\right)^{\delta}}
{|I^{(n)}_{0}(x_b;z_0)|^{\beta\delta}\left(
\sum\limits_{s=0}^{q_{n+1}-1}|I^{(n)}_{s}(x_{b})|^{-\delta}+\sum\limits_{k=0}^{q_{n}-1}|I^{(n+1)}_{k}(x_{b})|^{-\delta}\right)^{\beta}}=$$
$$=Const\lim\limits_{n\to\infty}\left(
\frac{\sum\limits_{s=0}^{q_{n+1}-1}|I^{(n)}_{s}(x_{b})|^{-\beta}+\sum\limits_{k=0}^{q_{n}-1}|I^{(n+1)}_{k}(x_{b})|^{-\beta}}
{\left(\sum\limits_{s=0}^{q_{n+1}-1}|I^{(n)}_{s}(x_{b})|^{-\delta}+
\sum\limits_{k=0}^{q_{n}-1}|I^{(n+1)}_{k}(x_{b})|^{-\delta}\right)^{\frac{\beta}{\delta}}}\right)^{\delta}.$$
For any $\beta \in R^{1 }$ we denote $$\tilde{S}_{n,\beta}=\sum\limits_{s=0}^{q_{n+1}-1}|I^{(n)}_{s}(x_{b})|^{-\beta}+\sum\limits_{k=0}^{q_{n}-1}|I^{(n+1)}_{k}(x_{b})|^{-\beta}.$$ We have
\begin{equation}\label{Inq.4}
\lim\limits_{n\to\infty}\left(
\frac{\tilde{S}_{n,\beta}}
{\tilde{S}_{n,\delta}\cdot\tilde{S}_{n,\delta}^{\frac{\beta-\delta}{\delta}}}\right)^{\delta}=
\lim\limits_{n\to\infty}\left(
\frac{\sum\limits_{s=0}^{q_{n+1}-1}|I^{(n)}_{s}(x_{b})|^{-\beta}+\sum\limits_{k=0}^{q_{n}-1}|I^{(n+1)}_{k}(x_{b})|^{-\beta}}
{\sum\limits_{s=0}^{q_{n+1}-1}|I^{(n)}_{s}(x_{b})|^{-\delta}\tilde{S}_{n,\delta}^{\frac{\beta-\delta}{\delta}}+\sum\limits_{k=0}^{q_{n}-1}|I^{(n+1)}_{k}(x_{b})|^{-\delta} \tilde{S}_{n,\delta}^{\frac{\beta-\delta}{\delta}}}\right)^{\delta}.\end{equation}
Next, for every $0\leq s<q_{n+1}$
$$\frac{|I^{(n)}_{s}(x_{b})|^{-\beta}}
{|I^{(n)}_{s}(x_{b})|^{-\delta} \tilde{S}_{n,\delta}^{\frac{\beta-\delta}{\delta}}}=
\left(\frac{|I^{(n)}_{s}(x_{b})|^{-\delta}}
{\tilde{S}_{n,\delta}}\right)^{\frac{\beta-\delta}{\delta}}=
\left(\frac{|I^{(n)}_{s}(x_{b})|^{-\delta}}
{\sum\limits_{s=0}^{q_{n+1}-1}|I^{(n)}_{s}(x_{b})|^{-\delta}+
\sum\limits_{k=0}^{q_{n}-1}|I^{(n+1)}_{k}(x_{b})|^{-\delta}}\right)^{\frac{\beta-\delta}{\delta}}<$$
$$<\left(\frac{|I^{(n)}_{s}(x_{b})|^{-\delta}}
{\sum\limits_{s=0}^{q_{n+1}-1}|I^{(n)}_{s}(x_{b})|^{-\delta}
}\right)^{\frac{\beta-\delta}{\delta}}<
\left(\frac{1}
{n e^{-v\delta}}\right)^{\frac{\beta-\delta}{\delta}}=\frac{e^{v(\beta-\delta)}}{n^{\frac{\beta-\delta}{\delta}}}$$
 Here we used the assertion of  lemma \ref{compare}.\\
 Thus we have
\begin{equation}\label{itemn-1}
|I^{(n)}_{s}(x_{b})|^{-\beta}< \frac{e^{v(\beta-\delta)}}{n^{\frac{\beta-\delta}{\delta}}} |I^{(n)}_{s}(x_{b})|^{-\delta} \tilde{S}_{n,\delta}^{\frac{\beta-\delta}{\delta}}.
\end{equation}
Analogously, we can show , for $0\leq k<q_{n}$
\begin{equation}\label{itemn}
|I^{(n+1)}_{k}(x_{b})|^{-\beta}< \frac{e^{v(\beta-\delta)}}{n^{\frac{\beta-\delta}{\delta}}} |I^{(n+1)}_{k}(x_{b})|^{-\delta} \tilde{S}_{n,\delta}^{\frac{\beta-\delta}{\delta}}.
\end{equation}
Using the bounds  (\ref{itemn-1}) and (\ref{itemn}) we have
$$\left(
\frac{\sum\limits_{s=0}^{q_{n+1}-1}|I^{(n)}_{s}(x_{b})|^{-\beta}+\sum\limits_{k=0}^{q_{n}-1}|I^{(n+1)}_{k}(x_{b})|^{-\beta}}
{\sum\limits_{s=0}^{q_{n+1}-1}|I^{(n)}_{s}(x_{b})|^{-\delta}\tilde{S}_{n,\delta}^{\frac{\beta-\delta}{\delta}}+\sum\limits_{k=0}^{q_{n}-1}|I^{(n+1)}_{k}(x_{b})|^{-\delta} \tilde{S}_{n,\delta}^{\frac{\beta-\delta}{\delta}}}\right)^{\delta}<
\frac{e^{v\delta(\beta-\delta)}}{n^{\beta-\delta}}$$
This bound and  relations  (\ref{Inq.3}) and (\ref{Inq.4}) prove
 $$\lim\limits_{n\to\infty}\frac{(\Lambda_{\beta}(x,q_n))^{\delta}}{(\Lambda_{\delta}(x,q_n))^{\beta}}=0.$$
 and hence also Lemma \ref{Eigen.1}
\endproof

\begin{lemm}\label{hatLambda} For any $z_0\in S^{1}\setminus \{T^{i}(x_b),\,i=0, -1,-2,...,-q_{n}-q_{n+1}+1\},$  and $n\geq1$
\begin{equation}\hat{\Lambda}(z_0,q_{n+1})\leq Const\cdot n\cdot \lambda_{-1}^{n} \end{equation}
where the positive constant $Const$  depends only on  the total variation of  $\log T'$.
\end{lemm}
\proof
One can easily see  (see \cite{Llave} remark 2.1)
\begin{equation}
\Lambda_1(z_0,n+t)=|DT^{t}(z_n)|\Lambda_1(z_0,n)+\Lambda_1(z_n, t).
\end{equation}
Let $1\leq i_n\leq  q_{n+1}$ realize the maximum of
$$\hat{\Lambda}(z_0,q_{n+1})=\max\limits_{1\leq i\leq q_{n+1}-1}\sum\limits_{k=1}^{i}\prod\limits_{j=k}^{i}\left|T'(z_j)\right|=\sum\limits_{k=1}^{i_n-1}\prod\limits_{j=k}^{i_n-1}\left|T'(z_j)\right|=\Lambda_1(z_0,i_n)-1,$$
and decompose it as
 $$i_n=b_{l_m} q_{l_m}+b_{l_{m-1}}q_{l_{m-1}}+...+b_{l_1}q_{l_1},$$ where $1\leq b_j\leq k_j,\, 0\leq l_1<l_2<...<l_m\leq n,$ since $\rho=[k_0,\,k_1,\,...,\, k_n,\,...]$ bounded type.
\begin{equation}
\begin{split}
\Lambda_1(z_0, i_{n})&=\Lambda_1(z_0, b_{l_m} q_{l_m}+b_{l_{m-1}}q_{l_{m-1}}+...+b_{l_1}q_{l_1})\\
&=|DT^{i_n-q_{l_m}}(z_{q_{l_m}})|\Lambda_1(z_0, q_{l_m})+\Lambda_1(z_{q_{l_m}}, i_{n}-q_{l_m})\\
&=\sum\limits_{j=1}^{b_{l_m}}|DT^{i_n-jq_{l_m}}(z_{j\cdot q_{l_m}})|\Lambda_1(z_{(j-1)\cdot q_{l_m}}, q_{l_m})+\Lambda_1(z_{b_{l_m}\cdot q_{l_m}}, i_{n}-b_{l_m}\cdot q_{l_m})\\
&=\sum\limits_{j=1}^{b_{l_m}}|DT^{i_n-jq_{l_m}}(z_{j\cdot q_{l_m}})|\Lambda_1(z_{(j-1)\cdot q_{l_m}}, q_{l_m})+\\
&+\sum\limits_{j=1}^{b_{l_{m-1}}}|DT^{i_n-b_{l_m}q_{l_m}-jq_{l_{m-1}}}(z_{b_{l_m}\cdot q_{l_m}+j\cdot q_{l_{m-1}}})|\Lambda_1(z_{b_{l_m}\cdot q_{l_m}+(j-1)\cdot q_{l_{m-1}}}, q_{l_{m-1}})+\\
&\qquad\ldots\\
&+\sum\limits_{j=1}^{b_{l_1}}|DT^{(b_{l_1}-j)q_{l_1}}(z_{i_n-(b_{l_1}-j})\cdot q_{l_1})|\Lambda_1(z_{i_n-(b_{l_1}-j+1)\cdot q_{l_1}}, q_{l_1})
\end{split}
\end{equation}
Applying Theorem \ref{Lyap1} we obtain:
\begin{equation}\label{hatL}
\begin{split}
\Lambda_1(z_0, i_{n})&\leq C_1\lambda_{-1}^{l_m}\sum\limits_{j=1}^{b_{l_m}}|DT^{i_n-jq_{l_m}}(z_{j\cdot q_{l_m}})||I^{(l_m)}(x_b;z_{(j-1)\cdot q_{l_m}})|+\\
&+C_1\lambda_{-1}^{l_{m-1}}\sum\limits_{j=1}^{b_{l_{m-1}}}|DT^{i_n-b_{l_m}q_{l_m}-jq_{l_{m-1}}}(z_{b_{l_m}\cdot q_{l_m}+j\cdot q_{l_{m-1}}})|I^{(l_{m-1})}(x_b;z_{b_{l_m}\cdot q_{l_m}+(j-1)\cdot q_{l_{m-1}}})|+\\
&\qquad\ldots\\
&+C_1\lambda_{-1}^{l_1}\sum\limits_{j=1}^{b_{l_1}}|DT^{(b_{l_1}-j)q_{l_1}}(z_{i_n-(b_{l_1-j})q_{l_1}})|I^{(l_1)}(x_b; z_{i_n-(b_{l_1}-j+1)q_{l_1}})|
\end{split}
\end{equation}
Next we estimate
$$DT^{i_n-N}(z_{N})|I^{(l_j)}(x_b;z_{N- q_{l_j}})|=\frac{DT^{i_n-N}(z_{N})|I^{(l_j)}(x_b;z_N)|}{\int\limits_{I^{(l_j)}(x_b;z_N)}DT^{i_n-N}(t)dt}\cdot \frac{|I^{(l_j)}(x_b;z_{N-q_{l_j}})|}{|I^{(l_j)}(x_b;z_N)|}\int\limits_{I^{(l_j)}(x_b;z_N)}DT^{i_n-N}(t)dt$$
for $i_n<q_{n+1},\, N<q_{n+1}$ and $1\leq j\leq m$. Using  Finzi's inequality we get
$$\frac{DT^{i_n-N}(z_{N})|I^{(l_j)}(x_b;z_N)|}{\int\limits_{I^{(l_j)}(x_b;z_N)}DT^{i_n-N}(t)dt}=\frac{|I^{(l_j)}(x_b;z_N)|}{\int\limits_{I^{(l_j)}(x_b;z_N)}\frac{DT^{i_n-N}(t)}{DT^{i_n-N}(z_{N})}dt}\leq Const_1.$$
Using Denjoy's inequality we obtain
$$\frac{|I^{(l_j)}(x_b;z_{N-q_{l_j}})|}{|I^{(l_j)}(x_b;z_N)|}\leq Const_2$$
and
$$\int\limits_{I^{(l_j)}(x_b;z_N)}DT^{i_n-N}(t)dt<1.$$
Where $Const_1$ and $Const_2$  only depend on the total variation of $\ln DT.$
Since $\rho$ is of bounded type $b_{l_j}\leq k_{l_j}\leq Const_3.$\\
From these inequalities and (\ref{hatL}) one  obtains
$$\Lambda_1(z_0, i_{n})<\lambda_{-1}^{l_m}C_1 \cdot Const_1\cdot Const_2\sum\limits_{j=1}^{m}b_{l_j}\leq \lambda_{-1}^{l_m}C_1\cdot Const_1\cdot Const_2\cdot Const_3\cdot m\leq Const \cdot n\cdot\lambda_{-1}^{n}.$$
\endproof

\section{ Barycentric coefficients}
 A universal bound for $T^n: S^1\to S^1$ is a constant that does not depend on $n$ and the point $y\in S^1$.\\
Let $[a, b]$ be an interval in $S^1$
and $c\in [a, b]$. The barycentric coefficient of $c$ in $[a, b]$
is the ratio $|[a, c]|/|[a, b]|$, where $|I|$ denotes the length of the interval $I$.

The intervals $[a, b],\,[b,c]$and $[a, c]$ are comparable with each other  if and only if the barycentric coefficient of c in $[a, b]$ is universally bounded in $(0,1)$.

We denote by $I^{(n)}(z)$ the interval of the $n$-th dynamical partition $P_n(T,x_0), \,\,x_{0}\in S^{1}$, containing $z$.

\begin{prop}\label{good}
 Let $T$ be a circle homeomorphism with irrational rotation number $\rho$ of bounded type satisfying  the conditions of Lemma \ref{DEN}. Suppose $x_{0}:=x_{b}\in S^1$ is the unique break point of $T$. Chose a point $z_0\in S^1$ whose orbit $\mathbb{O}_T(z_0)$ is disjoint from the one of  $x_0$.
There exists a subsequence of integers $\{n_{m},\,m=1,2,...\}$ such that
 the barycentric coefficient of the points $z_{k}:=T^{k}{z_{0}}$ in
$I^{(n_{m})}(z_k),\, 0\leq k<q_{n_{m}+1}$ is universally
bounded in $(0,1)$, in other words, there exists a universal  constant $0<K_1 <1$ such that the
length of each of the two connected components of
$I^{(n_{m})}(z_k)\setminus{z_k}$ is at least $K_1|I^{(n_{m})}(z_k)|$.
\end{prop}

Let $T_{\rho}x=x+\rho \mod 1$ be the rotation of $S^1$ by the  irrational angle $\rho \in (0,1)$.
 Suppose that the orbits  $\mathbb{O}_{T_{\rho}}(\tilde{x}_0),\,\,\mathbb{O}_{T_{\rho}}(\tilde{z}_0)$ of $\tilde{x}_{0}:=x_0=x_b$ and $\tilde{z}_0=z_0$
 are disjoint. Consider the dynamical partition  $P_{n}(T_\rho,\tilde{x}_0)$ of $S^1$ determined by $\tilde{x}_0$ and $T_\rho$ and let $\tilde{z}_k=T_{\rho}^k\tilde{z}_0,  \tilde{x}_k=T_{\rho}^k\tilde{x}_0,\, 0\leq k<q_{n+1}$.
 Then $\tilde{z}_0\in \tilde{I}^{(n+1)}_{k_0}(\tilde{x}_0)=T_{\rho}^{k_0}[\tilde{x}_{q_{n+1}},\tilde{x}_0)$ for some $0\leq k_0<q_n$  or
 $\tilde{z}_0\in \tilde{I}^{(n)}_{i_0}(\tilde{x}_0)=T_{\rho}^{i_0}([\tilde{x}_0,\tilde{x}_{-q_{n+1}})\cup[\tilde{x}_{-q_{n+1}},\tilde{x}_{q_n}))$
 for some $0\leq i_0<q_{n+1}$, i.e. $\tilde{z}_0\in T_{\rho}^{i_0}([\tilde{x}_0,\tilde{x}_{-q_{n+1}}))$ or $\tilde{z}_0\in T_{\rho}^{i_0}
 ([\tilde{x}_{-q_{n+1}},\tilde{x}_{q_n}))$.

  In the following we need some lemmas describing the location of the  points $\tilde{z}_k,\, 0\leq k<q_{n+1}$, in intervals of $P_{n}(T_\rho,\tilde{x}_0)$ which we formulate next.

\begin{lemm} \textbf{(Case 1.)}
If  $\tilde{z}_0\in \tilde{I}^{(n+1)}_{k_0}(\tilde{x}_0)$ for some $0\leq k_0<q_n$, then the points $\tilde{z}_k,\, 0\leq k<q_{n+1}$, belong to the following intervals of the dynamical partition $P_{n}(T_\rho, \tilde{x}_0)$:
\begin{itemize}
  \item $\tilde{z}_k\in \tilde{I}^{(n+1)}_{k_0+k}(\tilde{x}_0), \,0\leq k<q_n-k_0$,
  \item $\tilde{z}_k\in \tilde{I}^{(n)}_{k-q_n+k_0}(\tilde{x}_0), \,q_n-k_0\leq k<q_{n+1}$.
\end{itemize}
\end{lemm}

%%%%%%%%%%%%%%%%%%    Picture 3.1 %%%%%%%%%%%%%%%%%%%%%%
%%%%%%%%%%%%%%%%%%%% First two%%%%%%%%%%%%%%%%%%%%%%%%
%%%%%%%%%%%%%%%%%%%%%%%%%%%%%%%%%%%%%%%%%%%%%%%%%
\noindent
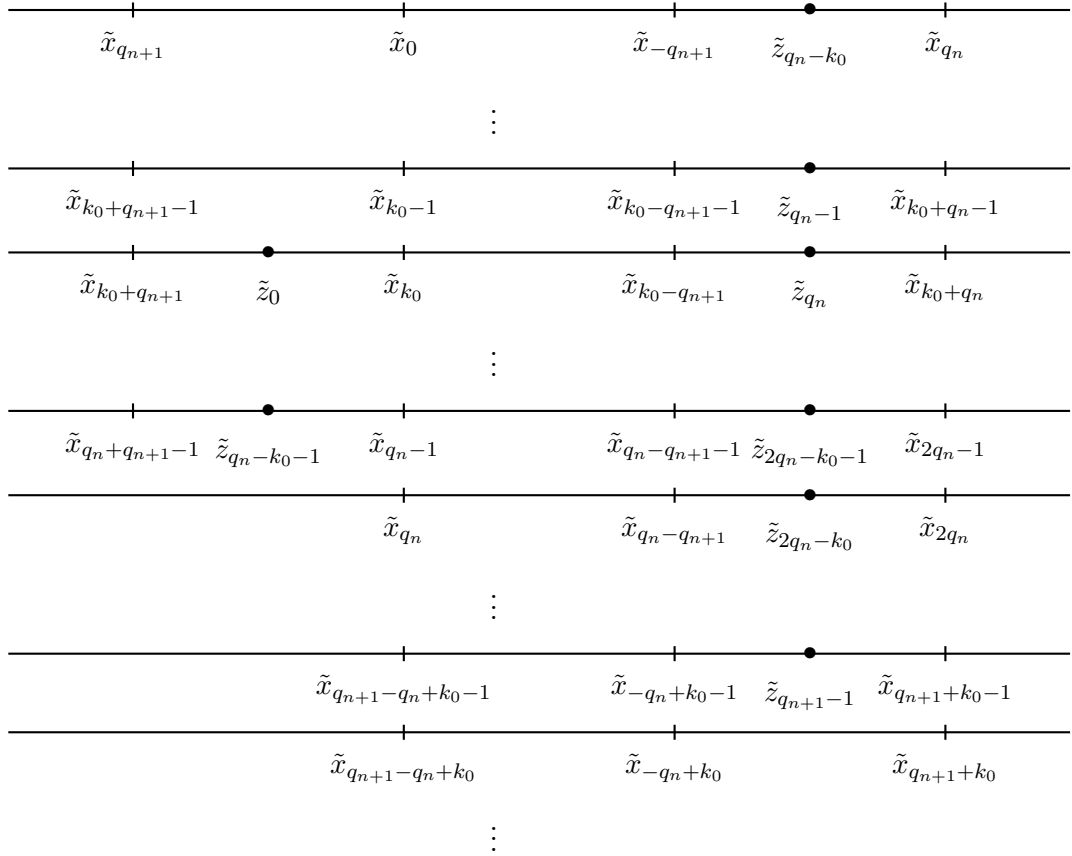
\begin {tikzpicture}[ ,auto ,node distance =1.8cm ,on grid ,
thick ,
state/.style ={ circle ,fill =black,
draw,black , text=white , minimum width =0.1 cm}]
\draw (1.8cm,-.1) -- (1.8cm, .1);   % 1 x 1.8
\draw (5.4cm,-.1) -- (5.4cm, .1);  % 3 x 1.8
\draw (9.0cm,-.1) -- (9.0cm, .1);   % 4 x 1.8
\draw (12.6cm,-.1) -- (12.6cm, .1); % 7 x 1.8
\node[] (A){};
\node[] (B) [right = of A, label=below: $\tilde{x}_{q_{n+1}}$ ] {};
\node[] (C) [right = of B, label=below: ] {};
\node[] (D) [right = of C, label=below: $\tilde{x}_{0}$] {};
\node[] (E) [right = of D, label=below: ] {};
\node[] (F) [right = of E, label=below: $\tilde{x}_{-q_{n+1}}$] {};
%\node[] (G) [right = of F] {};
%\node[] (H) [right = of G] {};
\node[] (G) [right = of F, label=below: $\tilde{z}_{q_{n}-k_0}$ ] {$\bullet$};
\node[] (H) [right = of G, label=below: $\tilde{x}_{q_{n}}$] {};
\node[] (I) [right = of H] {};
\path (A) edge [bend right = 0] (I);
\end{tikzpicture}
\begin{center}
$\vdots\qquad\qquad$
\end{center}
\noindent\begin {tikzpicture}[ ,auto ,node distance =1.8cm ,on grid ,
thick ,
state/.style ={ circle ,fill =black,
draw,black , text=white , minimum width =0.1 cm}]
\draw (1.8cm,-.1) -- (1.8cm, .1);   % 1 x 1.8
\draw (5.4cm,-.1) -- (5.4cm, .1);  % 3 x 1.8
\draw (9.0cm,-.1) -- (9.0cm, .1);   % 4 x 1.8
\draw (12.6cm,-.1) -- (12.6cm, .1); % 7 x 1.8
\node[] (A){};
\node[] (B) [right = of A, label=below: $\tilde{x}_{k_0+q_{n+1}-1}$ ] {};
\node[] (C) [right = of B, label=below: ] {};
\node[] (D) [right = of C, label=below: $\tilde{x}_{k_0-1}$] {};
\node[] (E) [right = of D, label=below: ] {};
\node[] (F) [right = of E, label=below: $\tilde{x}_{k_0-q_{n+1}-1}$] {};
%\node[] (G) [right = of F] {};
%\node[] (H) [right = of G] {};
\node[] (G) [right = of F, label=below: $\tilde{z}_{q_{n}-1}$ ] {$\bullet$};
\node[] (H) [right = of G, label=below: $\tilde{x}_{k_0+q_{n}-1}$] {};
\node[] (I) [right = of H] {};
\path (A) edge [bend right = 0] (I);
\end{tikzpicture}
\\\begin{tikzpicture}[ ,auto ,node distance =1.8cm ,on grid ,
thick ,
state/.style ={ circle ,fill =black,
draw,black , text=white , minimum width =0.1 cm}]
\draw (1.8cm,-.1) -- (1.8cm, .1);   % 1 x 1.8
\draw (5.4cm,-.1) -- (5.4cm, .1);  % 3 x 1.8
\draw (9.0cm,-.1) -- (9.0cm, .1);   % 4 x 1.8
\draw (12.6cm,-.1) -- (12.6cm, .1); % 7 x 1.8
\node[] (A){};
\node[] (B) [right = of A, label=below: $\tilde{x}_{k_0+q_{n+1}}$ ] {};
\node[] (C) [right = of B, label=below: $\tilde{z}_{0}$ ] {$\bullet$};
\node[] (D) [right = of C, label=below: $\tilde{x}_{k_0}$] {};
\node[] (E) [right = of D, label=below: ] {};
\node[] (F) [right = of E, label=below: $\tilde{x}_{k_0-q_{n+1}}$] {};
%\node[] (G) [right = of F] {};
%\node[] (H) [right = of G] {};
\node[] (G) [right = of F, label=below: $\tilde{z}_{q_{n}}$ ] {$\bullet$};
\node[] (H) [right = of G, label=below: $\tilde{x}_{k_0+q_{n}}$] {};
\node[] (I) [right = of H] {};
\path (A) edge [bend right = 0] (I);
\end{tikzpicture}
\begin{center}
$\vdots\qquad\qquad$
\end{center}
%%%%%%%%%%%%%%%%%%%%%%%%%%%%%%%%%%%%%%%%%%%%%%%%%
%%%%%%%%%%%%%%%%%%%% Second two%%%%%%%%%%%%%%%%%%%%%
%%%%%%%%%%%%%%%%%%%%%%%%%%%%%%%%%%%%%%%%%%%%%%%%%
\noindent\begin{tikzpicture}[ ,auto ,node distance =1.8cm ,on grid ,
thick ,
state/.style ={ circle ,fill =black,
draw,black , text=white , minimum width =0.1 cm}]
\draw (1.8cm,-.1) -- (1.8cm, .1);   % 1 x 1.8
\draw (5.4cm,-.1) -- (5.4cm, .1);  % 3 x 1.8
\draw (9.0cm,-.1) -- (9.0cm, .1);   % 4 x 1.8
\draw (12.6cm,-.1) -- (12.6cm, .1); % 7 x 1.8
\node[] (A){};
\node[] (B) [right = of A, label=below: $\tilde{x}_{q_n+q_{n+1}-1}$ ] {};
\node[] (C) [right = of B, label=below: $\tilde{z}_{q_n-k_0-1}$ ] {$\bullet$};
\node[] (D) [right = of C, label=below: $\tilde{x}_{q_n-1}$] {};
\node[] (E) [right = of D, label=below: ] {};
\node[] (F) [right = of E, label=below: $\tilde{x}_{q_n-q_{n+1}-1}$] {};
%\node[] (G) [right = of F] {};
%\node[] (H) [right = of G] {};
\node[] (G) [right = of F, label=below: $\tilde{z}_{2q_{n}-k_0-1}$ ] {$\bullet$};
\node[] (H) [right = of G, label=below: $\tilde{x}_{2q_{n}-1}$] {};
\node[] (I) [right = of H] {};
\path (A) edge [bend right = 0] (I);
\end{tikzpicture}
\\\begin{tikzpicture}[ ,auto ,node distance =1.8cm ,on grid ,
thick ,
state/.style ={ circle ,fill =black,
draw,black , text=white , minimum width =0.1 cm}]
%\draw (1.8cm,-.1) -- (1.8cm, .1);   % 1 x 1.8
\draw (5.4cm,-.1) -- (5.4cm, .1);  % 3 x 1.8
\draw (9.0cm,-.1) -- (9.0cm, .1);   % 4 x 1.8
\draw (12.6cm,-.1) -- (12.6cm, .1); % 7 x 1.8
\node[] (A){};
\node[] (B) [right = of A, label=below:  ] {};
\node[] (C) [right = of B, label=below:  ] {};
\node[] (D) [right = of C, label=below: $\tilde{x}_{q_n}$] {};
\node[] (E) [right = of D, label=below: ] {};
\node[] (F) [right = of E, label=below: $\tilde{x}_{q_n-q_{n+1}}$] {};
%\node[] (G) [right = of F] {};
%\node[] (H) [right = of G] {};
\node[] (G) [right = of F, label=below: $\tilde{z}_{2q_{n}-k_0}$ ] {$\bullet$};
\node[] (H) [right = of G, label=below: $\tilde{x}_{2q_{n}}$] {};
\node[] (I) [right = of H] {};
\path (A) edge [bend right = 0] (I);
\end{tikzpicture}
\begin{center}
$\vdots\qquad\qquad$
\end{center}
%%%%%%%%%%%%%%%%%%%%%%%%%%%%%%%%%%%%%%%%%%%%%%%%%
%%%%%%%%%%%%%%%%%%%% Last  %%%%%%%%%%%%%%%%%%%%%
%%%%%%%%%%%%%%%%%%%%%%%%%%%%%%%%%%%%%%%%%%%%%%%%%
\noindent\begin{tikzpicture}[ ,auto ,node distance =1.8cm ,on grid ,
thick ,
state/.style ={ circle ,fill =black,
draw,black , text=white , minimum width =0.1 cm}]
%\draw (1.8cm,-.1) -- (1.8cm, .1);   % 1 x 1.8
\draw (5.4cm,-.1) -- (5.4cm, .1);  % 3 x 1.8
\draw (9.0cm,-.1) -- (9.0cm, .1);   % 4 x 1.8
\draw (12.6cm,-.1) -- (12.6cm, .1); % 7 x 1.8
\node[] (A){};
\node[] (B) [right = of A, label=below:  ] {};
\node[] (C) [right = of B, label=below:  ] {};
\node[] (D) [right = of C, label=below: $\tilde{x}_{q_{n+1}-q_n+k_0-1}$] {};
\node[] (E) [right = of D, label=below: ] {};
\node[] (F) [right = of E, label=below: $\tilde{x}_{-q_n+k_0-1}$] {};
%\node[] (G) [right = of F] {};
%\node[] (H) [right = of G] {};
\node[] (G) [right = of F, label=below: $\tilde{z}_{q_{n+1}-1}$ ] {$\bullet$};
\node[] (H) [right = of G, label=below: $\tilde{x}_{q_{n+1}+k_0-1}$] {};
\node[] (I) [right = of H] {};
\path (A) edge [bend right = 0] (I);
\end{tikzpicture}
\\\begin{tikzpicture}[ ,auto ,node distance =1.8cm ,on grid ,
thick ,
state/.style ={ circle ,fill =black,
draw,black , text=white , minimum width =0.1 cm}]
%\draw (1.8cm,-.1) -- (1.8cm, .1);   % 1 x 1.8
\draw (5.4cm,-.1) -- (5.4cm, .1);  % 3 x 1.8
\draw (9.0cm,-.1) -- (9.0cm, .1);   % 4 x 1.8
\draw (12.6cm,-.1) -- (12.6cm, .1); % 7 x 1.8
\node[] (A){};
\node[] (B) [right = of A, label=below:  ] {};
\node[] (C) [right = of B, label=below:  ] {};
\node[] (D) [right = of C, label=below: $\tilde{x}_{q_{n+1}-q_n+k_0}$] {};
\node[] (E) [right = of D, label=below: ] {};
\node[] (F) [right = of E, label=below: $\tilde{x}_{-q_n+k_0}$] {};
%\node[] (G) [right = of F] {};
%\node[] (H) [right = of G] {};
\node[] (G) [right = of F, label=below:  ] {};
\node[] (H) [right = of G, label=below: $\tilde{x}_{q_{n+1}+k_0}$] {};
\node[] (I) [right = of H] {};
\path (A) edge [bend right = 0] (I);
\end{tikzpicture}
\begin{center}
$\vdots\qquad\qquad$
\end{center}
%\end{center}
%%%%%%%%%%%%%%%%%%%%%%%%%%%%%%%%%%%%%%%%%%%%%%%%%
%%%%%%%%%%%%%%%%%%%% Finished %%%%%%%%%%%%%%%%%%%%%
%%%%%%%%%%%%%%%%%%%%%%%%%%%%%%%%%%%%%%%%%%%%%%%%%

\begin{lemm}  \textbf{(Case 2.)}
Let   $\tilde{z}_0\in T_{\rho}^{i_0}([\tilde{x}_0,\tilde{x}_{-q_{n+1}}))$ for some $0\leq i_0<q_{n+1}-q_n.$ Then the points  $\tilde{z}_k,\, 0\leq k<q_{n+1}$, belong to the following intervals of the dynamical partition $P_{n}(T_\rho, \tilde{x}_0):$
\begin{itemize}
  \item $\tilde{z}_k\in T_{\rho}^{k+i_0}([\tilde{x}_0,\tilde{x}_{-q_{n+1}}))\subset \tilde{I}^{(n)}_{k+i_0}(\tilde{x}_0), \,0\leq k<q_{n+1}-i_0;$
  \item $\tilde{z}_k\in \tilde{I}^{(n+1)}_{k-q_{n+1}+i_0}(\tilde{x}_0), \,q_{n+1}-i_0\leq k<q_{n+1}.$
\end{itemize}
\end{lemm}

\begin{lemm}  \textbf{(Case 2'.)}
If  $\tilde{z}_0\in T_{\rho}^{i_0}([\tilde{x}_0,\tilde{x}_{-q_{n+1}}))$
for some $q_{n+1}-q_n\leq i_0<q_{n+1}$, then the points $\tilde{z}_k,\, 0\leq k<q_{n+1}$, belong to the following intervals of the dynamical partition $P_{n}(T\rho, \tilde{x}_0)$:
\begin{itemize}
 \item $\tilde{z}_k\in T_{\rho}^{k+i_0}([x_0,\tilde{x}_{-q_{n+1}}))\subset \tilde{I}^{(n)}_{k+i_0}(\tilde{x}_0), \,0\leq k<q_{n+1}-i_0$,
 \item $\tilde{z}_k\in \tilde{I}^{(n+1)}_{k-q_{n+1}+i_0}(\tilde{x}_0), \,q_{n+1}-i_0\leq k<q_{n+1}+q_n-i_0$,
\item $\tilde{z}_k\in \tilde{I}^{(n)}_{k-q_{n+1}-q_n+i_0}(\tilde{x}_0), \,q_{n+1}+q_n-i_0\leq k<q_{n+1}$.
\end{itemize}
\end{lemm}
\begin{lemm}  \textbf{(Case 3.)}
If   $\tilde{z}_0\in T_{\rho}^{i_0}([\tilde{x}_{-q_{n+1},\tilde{x}_{q_n}}))$ for some $0\leq i_0<q_{n+1}$, then the points of  $\tilde{z}_k,\, 0\leq k<q_{n+1}$, belong to the following intervals of the dynamical partition $P_{n}(T_\rho, \tilde{x}_0)$:
\begin{itemize}
 \item $\tilde{z}_k\in T_{\rho}^{k+i_0}([\tilde{x}_{-q_{n+1},\tilde{x}_{q_n}}))\subset \tilde{I}^{(n)}_{k+i_0}(\tilde{x}_0), \,0\leq k<q_{n+1}-i_0$,
 \item $\tilde{z}_k\in T_{\rho}^{k+i_0}([\tilde{x}_{-q_{n+1}},\tilde{x}_{q_n}))\subset \tilde{I}^{(n)}_{k-q_{n+1}+i_0}(\tilde{x}_0), \,q_{n+1}-i_0\leq k<q_{n+1}$.
\end{itemize}
\end{lemm}

Denote by $\tilde{I}^{(n)}(\tilde{z}_k)$ the interval of $P_n(T_\rho,\tilde{x}_0)$ which  contains $\tilde{z}_k$.
\begin{lemm}\label{ajral}
Let $T_{\rho}x=x+\rho \mod 1$ be the rotation of the circle with  $\rho$ irrational of bounded type. Chose a point $\tilde{z}_0\in S^1$ such that the orbits  $\mathbb{O}_{T_{\rho}}(\tilde{x}_0)$ and $\mathbb{O}_{T_{\rho}}(\tilde{z}_0)$   to be disjoint. Then there exists a subsequence $\{n_m\}$ of integers  such that the barycentric coefficient of the points $\tilde{z}_k$ in
$\tilde{I}^{(n_m)}(\tilde{z}_k)\in P_{n_m}(T_\rho,\tilde{x}_0) ,\,0\leq k< q_{n_m+1}$ is universally
bounded in $(0,1)$.
\end{lemm}

\proof
Denote by  $\tilde{P}_n(T_\rho,\tilde{x}_0)$ the partition generated by the points $\{\tilde{x}_{-q_{n+1}},\,\tilde{x}_{-q_{n+1}+1},\,\ldots,\,\tilde{x}_0,$ $\,\ldots,\, \tilde{x}_{q_n+q_{n+1}-1}\}$. First, we prove that there exists a subsequence of positive integers $\{ n_m\}$,  such that the barycentric coefficient
of the point $\tilde{z}_0$ in the interval $\tilde{\Delta}^{(n_m)}(\tilde{z}_0)\in\tilde{P}_{n_m}(T_\rho,\tilde{x}_0)$ is universally bounded in $(0,1)$.\\
Suppose on the contrary that there is no such subsequence. \\
{\bf{Fact 1.}}\emph{For sufficiently large $n$ the point $\tilde{z}_0$ is always at least in one of the two  intervals $\tilde{\Delta}^{(n+4)}_-$ and $\tilde{\Delta}^{(n+4)}_+$ of $\tilde{P}_{(n+4)}(T_\rho,\tilde{x}_0)$ contained and intersecting the boundary of $\tilde{\Delta}^{(n)}(\tilde{z}_0)$. } But if $\tilde{z}_0$ does not belong to  both these intervals, the point $\tilde{z}_0$ splits the interval $\tilde{\Delta}^{(n)}(\tilde{z}_0)$ into two intervals $\tilde{\Delta}^{(n)}_+(\tilde{z}_0)$ respectively $\tilde{\Delta}^{(n)}_-(\tilde{z}_0)$, the length of both these intervals being larger  than $|\tilde{\Delta}^{(n+4)}_-|$ respectively $|\tilde{\Delta}^{(n+4)}_+|$. Hence $\frac{|\tilde{\Delta}^{(n+4)}_\pm|}{|\tilde{\Delta}^{(n)}(\tilde{z}_0)|}\leq \frac{|\tilde{\Delta}^{(n)}_\pm(\tilde{z}_0)|}{|\tilde{\Delta}^{(n)}(\tilde{z}_0)|}\leq 1-\frac{|\tilde{\Delta}^{(n+4)}_\pm|}{|\tilde{\Delta}^{(n)}(\tilde{z}_0)|}$ and therefore $0<c<\frac{|\tilde{\Delta}^{(n+4)}_\pm|}{|\tilde{\Delta}^{(n)}(\tilde{z}_0)|}<C<1$ where these constants $c$ and $C$ do not depend on $n$. This shows that the barycentric coefficient of any $\tilde{y}_0\in \tilde{\Delta}^{(n)}(\tilde{z}_0)$ is universally bounded in $(0,1)$.\\

{\bf{Fact 2.}}\emph{ The point $\tilde{z}_0$ cannot be contained for all $n$ in the same  of the two intervals $\tilde{\Delta}^{(n+4)}_\pm$ of the partition $\tilde{P}_{n+4}(T_\rho,\tilde{x}_0)$  contained and intersecting the boundary of $\tilde{\Delta}^{(n)}(\tilde{z}_0)$}. Otherwise the point $\tilde{z}_0$ must converge for $n\to \infty$ to the limit boundary point of the intervals $\tilde{\Delta}^{(n)}(\tilde{z}_0)$, which however belongs to the orbit $\mathbb{O}_{T_\rho}(\tilde{x}_0)$ in contradiction to our assumption the orbits of $\tilde{z}_0$ and $\tilde{x}_0$ being disjoint.
Hence there exists a subsequence $n_m$ such that $\tilde{z}_0$ belongs alternatively to an interval in  $\tilde{P}^{(n+4)}(\rho,\tilde{x}_0)$ contained in and intersecting  the left respectively the right  boundary of $\tilde{\Delta}^{(n)}(\tilde{z}_0)$. Therefore the point $\tilde{z_0}$ must be separated from both the boundary points of $\tilde{\Delta}^{(n)}(\tilde{z}_0)$ by at least one interval in  $\tilde{P}^{(n+6)}(\rho,\tilde{x}_0)$ contained in and intersecting the boundary of $\tilde{\Delta}^{(n)}(\tilde{z}_0)$. Otherwise there would be another subsequence $n_l$ with $\tilde{z}_{n_l}$ approaching the limit boundary point of the intervals $\tilde{\Delta}^{(n_l)}(\tilde{z_0})$ which however is on $\mathbb{O}_{T_\rho}(\tilde{x_0})$, contrary to our assumption. That shows that the barycentric coefficient of the point $\tilde{z}_0$ in $\tilde{\Delta}^{(n)}(\tilde{z}_0)$ is universally bounded in $(0,1)$. This holds true also for the interval $\tilde{\Delta}^{(n)}(\tilde{z}_0)$ replaced by the interval $\tilde{I}^{(n)}\in P_n(T_\rho,\tilde{x}_0)$ with $\tilde{x}_0=x_0$.
\\
Next we proof the barycentric coefficient of $\tilde{z}_k$ in $\tilde{I}^{(n)}(\tilde{z}_k)\in {P}_{n}(T_\rho,\tilde{x}_0)$ for all $0\leq k\leq q_n-1$  to be universally bounded in $(0,1)$.\\
For the point $\tilde{z}_0$  the following three cases are possible:
\begin{enumerate}
  \item $\tilde{z}_0\in \tilde{I}^{(n+1)}_{k_0}(\tilde{x}_0),\,0\leq k_0\leq q_{n}-1 ;$
  \item $\tilde{z}_0\in [T_{\rho}^{k_0}({\tilde{x}_0}),T_{\rho}^{k_0-q_{n+1}}({\tilde{x}_0}))\subset \tilde{I}^{(n)}_{k_0}(\tilde{x}_0),\, 0\leq k_0\leq q_{n+1}-1;$
  \item $\tilde{z}_0\in [T_{\rho}^{k_0-q_{n+1}}(\tilde{x}_0),T_{\rho}^{k_0+q_{n}}(\tilde{x_0}))\subset \tilde{I}^{(n)}_{k_0}(\tilde{x_0}),\, 0\leq k_0\leq q_{n+1}-1$.
\end{enumerate}
 \textbf{Case 1.} $\tilde{z}_0\in \tilde{I}^{(n+1)}_{k_0}(\tilde{x_0}),\,0\leq k_0\leq q_{n}-1$ and $\tilde{z}_0$ is separated from both boundary points of $\tilde{I}^{(n+1)}_{k_0}(\tilde{x_0})$. Then for $0\leq k\leq q_n-k_0-1$, $\tilde{z}_k$ is also separated form the boundary points of $\tilde{I}^{(n+1)}_{k+k_0}(\tilde{x_0})$.
For $q_n-k_0 \leq k\leq q_{n+1}-1$, $\tilde{z}_k$ is clearly separated form the boundary points of $T^{k}_{\rho}\tilde{I}^{(n+1)}_{q_n-1}(\tilde{x_0})$. But since  $T^{k}_{\rho}\tilde{I}^{(n+1)}_{q_n-1}(\tilde{x_0})\subset \tilde{I}^{(n)_{k-q_n-k_0}}(\tilde{x_0})$, the barycentric coefficient of $\tilde{z}_k$ in $\tilde{I}^{(n)_{k-q_n-k_0}}(\tilde{x_0})$ is universally bounded in $(0,1)$.\\
\textbf{Case 2.} $\tilde{z}_0\in [T_{\rho}^{k_0}(\tilde{x_0}),T_{\rho}^{k_0-q_{n+1}}(\tilde{x_0}))\subset \tilde{I}^{(n)}_{k_0}(\tilde{x_0}),\, 0\leq k_0\leq q_{n+1}-1$ and $\tilde{z}_0$ is separated from the points  $T_{\rho}^{k_0}(\tilde{x_0})$ and $T_{\rho}^{k_0-q_{n+1}}(\tilde{x_0})$. Then for $0\leq k\leq q_{n+1}-k_0-1$, $z_k$ is separated from the points  $T_{\rho}^{k+k_0}(\tilde{x_0})$ and $T_{\rho}^{k+k_0-q_{n+1}}(\tilde{x_0})$. Since $[T_{\rho}^{k+k_0}(\tilde{x_0}),T_{\rho}^{k+k_0-q_{n+1}}(\tilde{x_0}))\subset \tilde{I}^{(n)}_{k+k_0}(\tilde{x_0})$, the barycentric coefficient of $\tilde{z}_k$ in $\tilde{I}^{(n)}_{k+k_0}(\tilde{x_0})$ is universally bounded. For $q_{n+1}-k_0\leq k\leq q_{n+1}-1$, $\tilde{z}_k$ is separated from the points  $T_{\rho}^{k+k_0}(\tilde{x_0})$ and $T_{\rho}^{k+k_0-q_{n+1}}(\tilde{x_0})$, but $[T_{\rho}^{k+k_0}(\tilde{x_0}),T_{\rho}^{k+k_0-q_{n+1}}(\tilde{x_0}))=\tilde{I}^{(n+1)}_{k+k_0-q_{n+1}}(\tilde{x_0})$. Then the barycentric coefficient of $\tilde{z}_k$ in $\tilde{I}^{(n+1)}_{k+k_0-q_{n+1}}(\tilde{x_0})$ is universally bounded in $(0,1)$.\\
\textbf{Case 3.} $\tilde{z}_0\in [T_{\rho}^{k_0-q_{n+1}}(\tilde{x_0}),T_{\rho}^{k_0+q_{n}}(\tilde{x_0}))\subset \tilde{I}^{(n)}_{k_0}(\tilde{x_0}),\, 0\leq k_0\leq q_{n+1}-1$ and $\tilde{z}_0$ is separated from the points  $T_{\rho}^{k_0-q_{n+1}}(\tilde{x_0})$ and $T_{\rho}^{k_0+q_{n}}(\tilde{x_0})$. For $0\leq k\leq q_{n+1}-k_0-1$, $\tilde{z}_k$ is separated from the points  $T_{\rho}^{k+k_0-q_{n+1}}(\tilde{x_0})$ and $T_{\rho}^{k+q_{n}}(\tilde{x_0})$. Since $[T_{\rho}^{k+k_0-q_{n+1}}(\tilde{x_0}),T_{\rho}^{k+k_0-q_{n+1}}(\tilde{x_0}))\subset \tilde{I}^{(n)}_{k+k_0}(\tilde{x_0})$, the barycentric coefficient of $\tilde{z}_k$ in $\tilde{I}^{(n)}_{k+k_0}(\tilde{x_0})$ is universally bounded. For $q_{n+1}-k_0\leq k\leq q_{n+1}-1$, $\tilde{z}_k$ is separated from the points  $T_{\rho}^{k+k_0-q_{n+1}}(\tilde{x_0})$ and $T_{\rho}^{k+q_{n}}(\tilde{x_0})$, but $[T_{\rho}^{k+k_0-q_{n+1}}(\tilde{x_0}),T_{\rho}^{k+q_{n}}(\tilde{x_0}))=\tilde{I}^{(n)}_{k+k_0-q_{n+1}}(\tilde{x_0})$. Then the barycentric coefficient of $\tilde{z}_k$ in $\tilde{I}^{(n+1)}_{k+k_0-q_{n+1}}(\tilde{x_0})$ is universally bounded in $(0,1)$.
\endproof

%%%%%%%%%%%%%%%%%%%% letter%%%%%%%%%%%%%%%%%%%%%%%

We proved in Lemma 4.6 that there exists a subsequence of positive integers $\{ n_m\}$,  such that the barycentric coefficient
of the point $\tilde{z}_0$ in the interval $\tilde{\Delta}^{(n_m)}(\tilde{z}_0)\in\tilde{P}_{n_m}(T_\rho,\tilde{x}_0)$ is universally bounded in $(0,1)$.\\
That means there exists a positive integer $l:=l(z_{0}),$ such that the point
 $\tilde{z}_0$ cannot  be in  the two  intervals $\tilde{\Delta}^{(n_m+l)}_-$ and $\tilde{\Delta}^{(n_{m}+l)}_+$ of $\tilde{P}_{(n_{m}+l)}(T_\rho,\tilde{x}_0)$ contained and intersecting the boundary of $\tilde{\Delta}^{(n_{m})}(\tilde{z}_0)$.
There are then three possible cases:\\
\textbf{Case 1.}
$\tilde{z}_0\in \tilde{\Delta}^{(n_m+1)}_{k_0}(\tilde{z_0}):=T_\rho^{k_0}[\tilde{x}_{q_{n_m+1}},\tilde{x}_0]\in \tilde{P}_{n_m}(T_\rho,\tilde{x}_0),\,\,\text{for some } k_{0},\,\,0\leq k_0\leq q_{n_m}-1$\\
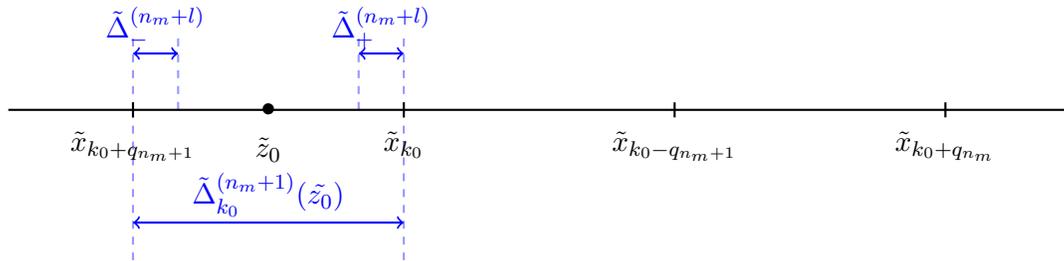
\begin{tikzpicture}[ ,auto ,node distance =1.8cm ,on grid ,
thick ,
state/.style ={ circle ,fill =black,
draw,black , text=white , minimum width =0.1 cm}]
\lineann[-2]{1.8}{5.4}{$\tilde{\Delta}^{(n_m+1)}_{k_0}(\tilde{z_0})$};
\lineann[1]{1.8}{2.4}{$\tilde{\Delta}^{(n_m+l)}_-$};
\lineann[1]{4.8}{5.4}{$\tilde{\Delta}^{(n_m+l)}_+$};
\draw (1.8cm,-.1) -- (1.8cm, .1);   % 1 x 1.8
\draw (5.4cm,-.1) -- (5.4cm, .1);  % 3 x 1.8
\draw (9.0cm,-.1) -- (9.0cm, .1);   % 4 x 1.8
\draw (12.6cm,-.1) -- (12.6cm, .1); % 7 x 1.8
\node[] (A){};
\node[] (B) [right = of A, label=below: $\tilde{x}_{k_0+q_{n_m+1}}$ ] {};
\node[] (C) [right = of B, label=below: $\tilde{z}_{0}$ ] {$\bullet$};
\node[] (D) [right = of C, label=below: $\tilde{x}_{k_0}$] {};
\node[] (E) [right = of D, label=below: ] {};
\node[] (F) [right = of E, label=below: $\tilde{x}_{k_0-q_{n_m+1}}$] {};
%\node[] (G) [right = of F] {};
%\node[] (H) [right = of G] {};
\node[] (G) [right = of F] {};
\node[] (H) [right = of G, label=below: $\tilde{x}_{k_0+q_{n_m}}$] {};
\node[] (I) [right = of H] {};
\path (A) edge [bend right = 0] (I);
\end{tikzpicture}\\
\textbf{Case 2.}
 $\tilde{z}_0\in \tilde{\Delta}_{k_0}^{-(n_m+1)}:= T_\rho^{k_0}[\tilde{x}_0,\tilde{x}_{-q_{n_m+1}}]\in \tilde{P}_{n_m}(T_\rho,\tilde{x}_0)$ for some $k_0, 0\leq k_0\leq q_{n_m+1}-1$. Obviously $\tilde{\Delta}_{k_0}^{-(n_m+1)}\subset \tilde{\Delta}_{k_0}^{(n_m)}$.\\
\begin{tikzpicture}[ ,auto ,node distance =1.8cm ,on grid ,
thick ,
state/.style ={ circle ,fill =black,
draw,black , text=white , minimum width =0.1 cm}]
\lineann[-2]{5.4}{9.0}{$\tilde{\Delta}^{-(n_m+1)}_{k_0}(\tilde{z_0})$};
\lineann[1]{5.4}{6}{$\tilde{\Delta}^{(n_m+l)}_-$};
\lineann[1]{8.4}{9}{$\tilde{\Delta}^{(n_m+l)}_+$};
\draw (1.8cm,-.1) -- (1.8cm, .1);   % 1 x 1.8
\draw (5.4cm,-.1) -- (5.4cm, .1);  % 3 x 1.8
\draw (9.0cm,-.1) -- (9.0cm, .1);   % 4 x 1.8
\draw (12.6cm,-.1) -- (12.6cm, .1); % 7 x 1.8
\node[] (A){};
\node[] (B) [right = of A, label=below: $\tilde{x}_{k_0+q_{n_m+1}}$ ] {};
\node[] (C) [right = of B] {};
\node[] (D) [right = of C, label=below: $\tilde{x}_{k_0}$] {};
\node[] (E) [right = of D, label=below: $\tilde{z}_{0}$ ] {$\bullet$};
\node[] (F) [right = of E, label=below: $\tilde{x}_{k_0-q_{n_m+1}}$] {};
%\node[] (G) [right = of F] {};
%\node[] (H) [right = of G] {};
\node[] (G) [right = of F] {};
\node[] (H) [right = of G, label=below: $\tilde{x}_{k_0+q_{n_m}}$] {};
\node[] (I) [right = of H] {};
\path (A) edge [bend right = 0] (I);
\end{tikzpicture}\\
\textbf{Case 3.}
$\tilde{z}_0\in \tilde{\Delta}_{k_0}^{-(n_m+1),n_m}:= T_\rho^{k_0}[\tilde{x}_{-q_{n_m+1}},\tilde{x}_{q_{n_m}}]\in \tilde{P}_{n_m}(T_\rho,\tilde{x_0}$ for some $k_0, 0\leq k_0\leq q_{n_m+1}-1$. Obviously $ \tilde{\Delta}_{k_0}^{-(n_m+1),n_m}\subset \tilde{\Delta}_{k_0}^{(n_m)}$.\\
\begin{tikzpicture}[ ,auto ,node distance =1.8cm ,on grid ,
thick ,
state/.style ={ circle ,fill =black,
draw,black , text=white , minimum width =0.1 cm}]
\lineann[-2]{9}{12.6}{$\tilde{\Delta}^{(-(n_m+1),n_m)}_{k_0}(\tilde{z_0})$};
\lineann[1]{9}{9.6}{$\tilde{\Delta}^{(n_m+l)}_-$};
\lineann[1]{12}{12.6}{$\tilde{\Delta}^{(n_m+l)}_+$};
\draw (1.8cm,-.1) -- (1.8cm, .1);   % 1 x 1.8
\draw (5.4cm,-.1) -- (5.4cm, .1);  % 3 x 1.8
\draw (9.0cm,-.1) -- (9.0cm, .1);   % 4 x 1.8
\draw (12.6cm,-.1) -- (12.6cm, .1); % 7 x 1.8
\node[] (A){};
\node[] (B) [right = of A, label=below: $\tilde{x}_{k_0+q_{n_m+1}}$ ] {};
\node[] (C) [right = of B] {};
\node[] (D) [right = of C, label=below: $\tilde{x}_{k_0}$] {};
\node[] (E) [right = of D ] {};
\node[] (F) [right = of E, label=below: $\tilde{x}_{k_0-q_{n_m+1}}$] {};
%\node[] (G) [right = of F] {};
%\node[] (H) [right = of G] {};
\node[] (G) [right = of F, label=below: $\tilde{z}_{0}$] {$\bullet$};
\node[] (H) [right = of G, label=below: $\tilde{x}_{k_0+q_{n_m}}$] {};
\node[] (I) [right = of H] {};
\path (A) edge [bend right = 0] (I);
\end{tikzpicture}\\
One also can say, there exists at least one
interval of the partition $\tilde{P}_{(n_m+l)}(T_\rho,\tilde{x}_0)$ between the points $\tilde{z_k}$ and
$\tilde{x}_{-j}, 0\leq j\leq q_{n_m+1}-1$.

%%%%%%%%%%%%%%%%%%%% letter%%%%%%%%%%%%%%%%%%%%%%%

\begin{lemm}\label{bound} (see \cite{DzLi}). Let $T$ be a circle homeomorphism with irrational rotation number of bounded type satisfying the conditions of Lemma \ref{DEN}.  Put $\theta_{\pm}:=(1+e^{\pm v})^{-\frac{1}{2}}$ and let $n$ and $l\geq 2$ be positive integers. Then there exist universal positive constants $C_{1},\,C_{2},\,\,C_{1}<C_{2}$  such that for arbitrary ${I}^{(n)} \in P_{n}(T, x_0)$ and $I^{(n+l)}\in P_{n+l}(T, x_0)$ with $I^{(n+l)}\subset I^{(n)}$ the following bounds hold
$$C_{1} \theta_{+}^{l} \leq \frac{I^{(n+l)}}{I^{(n)}}\leq C_{2} \theta_{-}^l.$$
\end{lemm}

\textbf{Proof of Proposition \ref{good}.}

\proof
Consider the dynamical partition $P_n(T,x_0)$   and the finite part   $ \mathbb{O}_T^{_{n+1}}(z_0)= \{z_k,\,0\leq k< q_{n+1}\} $  of the orbit  of $z_0$ under the map $T$. $I^{(n)}(z_k)$ denotes again the interval in $P_{n}(T,x_0)$ which contains $z_k$.
Since $T$ and $T_{\rho}$ are topologically conjugate, Lemma \ref{ajral} shows,  that there exists for $l\geq 4$ a subsequence $\{n_m\}$ in $\mathbb{N}$ such that
the point $z_k$ cannot be in the two boundary intervals of the partition $P_{n_m+l}(T,x_0)$  contained in $I^{(n_m)}(z_k)$. According to Lemma \ref{bound}  the  barycentric coefficient of $z_k$ in $I^{(n_m)}(z_k)$   hence is universally bounded in $(0,1)$.
\endproof

\section{Proof of Theorem \ref{main}}

Take any point $z_{0}$ not belonging to  the orbit $\mathbb{O}_T(x_b)$ of the break point $x_b$.
Before giving the main steps in  the proof of our main theorem we emphasize an important point:
a key ingredient in the proof of this theorem  is
the Taylor expansion of the process $\bar{z}_{q_{n+1}}(z_0,\sigma)$
\begin{equation} \label{Tayl1}
\bar{z}_{q_{n+1}}(z_0, \sigma)=T^{q_{n+1}}(z_0)+\sigma L_{q_{n+1}}(z_0)+\sigma^2 Q_{q_{n+1}}(z_0,\sigma),
\end{equation}
where the linear term $L_{q_{n+1}}$ is the sum of independent random variables as
defined in (\ref{Ln}). However we can use  Taylor's formula (\ref{Tayl1}) only  in case
the neighborhood $U ^ {(n)}$  of the point $z_{q_{n+1}}=T^{q_{n+1}}(z_0) $ does not contain any break point of  $T^{q_{n+1}}$.\newline
Let us briefly sketch the main steps of the proof of Theorem \ref{main}, which will overcome this difficulty. \newline
\textbf{I.}  By Lemma 4.1 there exists an increasing  sequence of natural numbers $\{n_{m},m=1,2,...\}$, such that
for $0 \leq k< q_{n_{m}+1}$  the interval $I^{(n_{m})}(z_ {k})$ of
the partition $ P^{(n_{m})} (T, x_{b})$ contains only the point $z_{k}:=T^{k}(z_{0})$ from the finite orbit $\{z_i=T^i z_0,\, 0\leq i <q_{n_m+1}\}$. Moreover the barycentric coordinates of the points $z_{k}$ are strictly separated
from 0 and 1 by universal constants, i.e. the points $z_k$ are uniformly separated from the boundaries $\partial I^{(n_m)}(z_k)$.\newline
For every $m\geq 1$ we will construct  a series of neighborhoods $A^{n_{m}}_k\subset I^{(n_{m})}(z_{k}),\,\, 0\leq k < q_{n_{m}+1}$,
 of the points $z_{k}$, such that every interval   $A^{n_{m}}_k$ does not contain any break point of $T^{q_{n_m+1}}$.\newline
\textbf{II.}
For the stochastic sequence   $$\bar{z}_k=T(\bar{z}_{k-1})+\sigma_{q_{n_m+1}}\xi_k,\,\,\,\bar{z}_0=z_0$$
we show that the probabilities of the events
$$B_{n_m}:=\{\bar{z}_1\in A^{n_m}_{1}, \,\, \bar{z}_2\in A^{n_m}_{2},\,...,\,
\bar{z}_{q_{n_m+1}-1}\in A^{n_m}_{q_{n_m+1}-1}\}$$
 tend to $1$ as $m\to\infty$.\newline
 \textbf{III.}
For fixed $z_0\in S^1\setminus\{x_b\}$    the Taylor expansion of the process $\bar{z}_{q_{n_m+1}}(z_0,\sigma_{q_{n_m+1}-1})$
in the variables $\xi_1,\,\xi_2,\,...,\, \xi_{q_{n_m+1}-1}$   allows to decompose it under the condition $\bar{z}_j\in A_j^{n_m},\,1\leq j\leq q_{n_m+1}-1$, as
\begin{equation}
\bar{z}_{q_{n_m+1}-1}(z_0, \sigma_{q_{n_m+1}-1})=
\end{equation}
$$T^{q_{n_m+1}-1}(z_0)+\sigma_{q_{n_m+1}-1}L_{q_{n_m+1}-1}(z_0)+\sigma_{q_{n_m+1}-1}^2 Q_{q_{n_m+1}-1}(z_0,\sigma_{q_{n_m+1}-1}),
$$
where the linear term $L_{q_{n_m+1}}$ is the sum of independent random variables defined by (\ref{Ln}).\\
\textbf{VI.}
We prove the CLT for this linear part $ L_ {q_ {n_m + 1}-1} $ which finally leads to the proof of  Theorem \ref{main}.\\

To achieve this program we  formulate and prove in a first step several lemmas.\\

 For this take a point $z_0\in S^1\setminus \mathbb{O}_T(x_b)$ and
  the sequence of increasing natural numbers $n_{m},\,m=1,2,...$ determined
 by Proposition 4.1. Consider the two partitions $\tilde{P}_{n_m}(T,x_b)$ (generated by the points $\{x_{-q_{n_m+1}},\,\ldots,\,x_0,\,\ldots,\,x_{q_{n_m}+q_{n_m+1}-1}\}$) respectively $\tilde{P}_{n_m+l}(T,x_b) $.
Each interval $\Delta^{(n_m)}(z_k)\in \tilde{P}_{n_m}(T,x_b), 0\leq k<q_{n_m+1}$ contains at least $q_l\geq q_4$ and hence at least three intervals of the partition $\tilde{P}_{n_m+l}(T,x_b)$.
By Proposition 4.1 the point $ z_k $ cannot be located
in the intervals $\Delta_\pm^{n_m+l}$ of the partition $ \tilde{P}_ {n_m + l} (T, x_t) $ contained in and intersecting the interval $ \Delta^{(n_m)} (z_k) $. \\
 We introduce  next  certain connected intervals $A^{(n_m)}_k$ composed of intervals of the partition $ \tilde{P}_{2n_m+l+1}(x_b)$ and containing the point $z_k$.
For this let $\Delta^{(2n_m+l+1)}(z_0)$ be the element of the partition $\tilde{P}_{2n_m+l+1}(x_b)$ which contains $z_0$.
 As above there are again three possible cases for this interval: $\Delta^{(2n_m+l+1)}(z_0)=\Delta_{{t}_0}^{2n_m+l+2}=[x_{{t}_0+q_{2n_m+l+2}}, x_{{t}_0}]$,  $\Delta^{(2n_m+l+1)}(z_0)=\Delta_{{t}_0}^{-(2n_m+l+2)}=[x_{{t}_0},x_{{t}_0-q_{2n_m+l+2}}]$ respectively $\Delta^{(2n_m+l+1)}(z_0)=\Delta_{{t}_0}^{-(2n_m+l+2),2n_m+l+1}=[x_{{t}_0-q_{2n_m+l+2}},x_{{t}_0+q_{2n_m+l+1}}]$. Thereby $t_0=t_0(k_0,l)$ is choosen such that $\Delta^{(2n_m+l+1)}(z_0)\subset \Delta_{k_0}^{n_m+1}(z_0)\setminus (\Delta_-^{n_m+l}\cup\Delta_+^{n_m+l})$.
Then we define
\begin{itemize}
  \item  $A^{(n_m)}_0:=\Delta^{(2n_m+l+1)}(z_0)\in \tilde{P}_{2n_m+l+1}(x_b)$,
  \item for every  $k,\,\, 1\leq k<q_{n_m+1}$, we set
  $$A^{(n_m)}_k:=A^{-}_{k}(n_{m})\cup T(A^{(n_m)}_{k-1})\cup A^{+}_{k}(n_{m}),$$
\end{itemize}
where $A^{-}_{k}(n_{m})$ and $A^{+}_{k}(n_{m})$  are the left  and right neighbours of  $T(A^{(n_m)}_{k-1})$ in the partition $ \tilde{P}_{2n_m+l+1}(x_b)$  respectively.\\
\begin{lemm}
 $A^{(n_m)}_k\subset \Delta^{(n_m)}(z_k),$ for all $0\leq k<q_{n_m+1}.$
\end{lemm}
\proof
Let the sequence $\left\{ n_{m},\,\,m=1,2,...\right\}$ be determined by Lemma 4.1.
 Consider the partitions
$\tilde{P}_{n_m}(T,x_b),\,\tilde{P}_{n_m+l}(T,x_b)$ and $\tilde{P}_{2n_m+l+1}(T,x_b)$.
 Recall that every interval of $\tilde{P}_{n_m}(T,x_b)$ contains
 at least three intervals of $\tilde{P}_{n_m+l}(T,x_b)$.\\
 Furthermore every interval of $\tilde{P}_{n_m+l}(T,x_b)$ contains  at least $q_{n_m+1}$ intervals of $\tilde{P}_{2n_m+l+1}(T,x_b)$.
 This implies that
 $$A^{(n_m)}_{k}\subset \Delta^{(n_m)}(z_{k})\in P_{n_m}(T,x_b),\,0\leq k<q_{n_m+1}.$$
\endproof
\textbf{The intervals $A_{k}^{(n_{m})}$ in red in  case I.}\\
%%%%%%%%%%%%%%%%%%    Picture 3.1 %%%%%%%%%%%%%%%%%%%%%%
%%%%%%%%%%%%%%%%%%%% First two%%%%%%%%%%%%%%%%%%%%%%%%
%%%%%%%%%%%%%%%%%%%%%%%%%%%%%%%%%%%%%%%%%%%%%%%%%
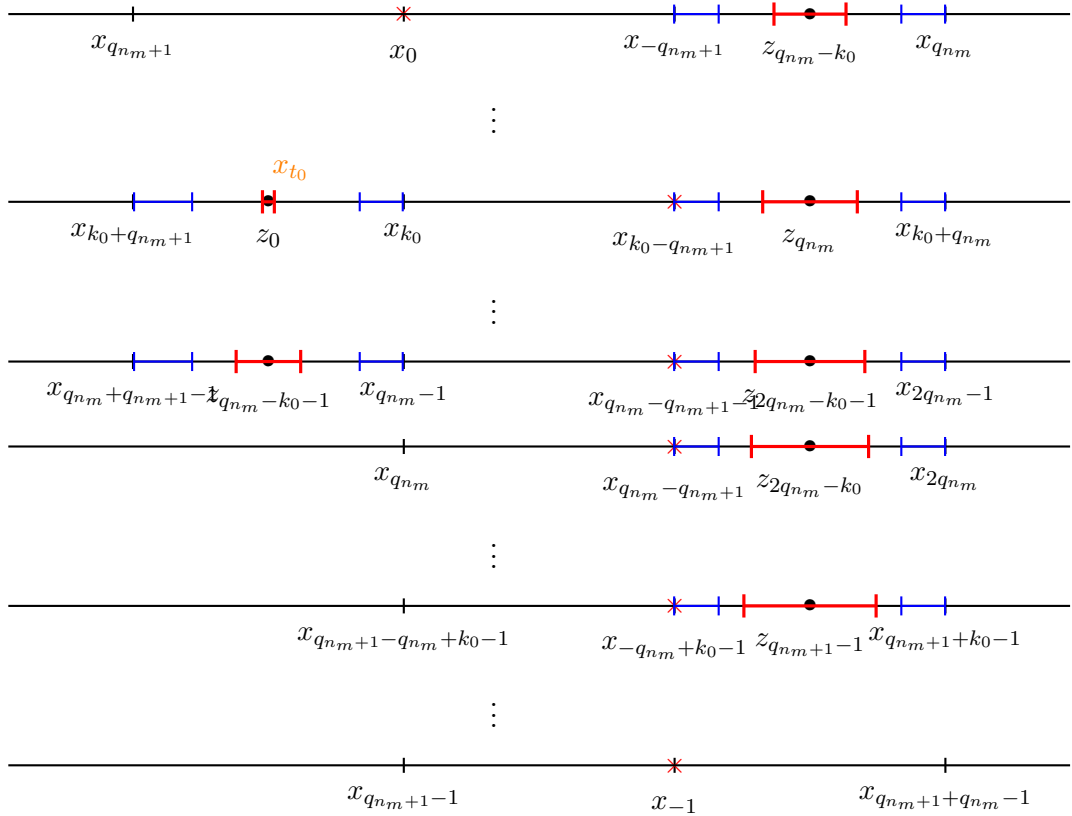
\noindent\begin {tikzpicture}[ ,auto ,node distance =1.8cm ,on grid ,
thick ,
state/.style ={ circle ,fill =black,
draw,black , text=white , minimum width =0.1 cm}]
\draw (1.8cm,-.1) -- (1.8cm, .1);   % 1 x 1.8
\draw (5.4cm,-.1) -- (5.4cm, .1);  % 3 x 1.8
\draw (9.0cm,-.1) -- (9.0cm, .1);   % 4 x 1.8
\draw (12.6cm,-.1) -- (12.6cm, .1); % 7 x 1.8
\node[] (A){};
\node[] (B) [right = of A, label=below: $x_{q_{n_m+1}}$ ] {};
\node[] (C) [right = of B, label=below: ] {};
\node[] (D) [right = of C, label=below: $x_{0}$] {\textcolor{red}{$\times$}};
\node[] (E) [right = of D, label=below: ] {};
\node[] (F) [right = of E, label=below: $x_{-q_{n_m+1}}$] {};
%\node[] (G) [right = of F] {};
%\node[] (H) [right = of G] {};
\node[] (G) [right = of F, label=below: $z_{q_{n_m}-k_0}$ ] {$\bullet$};
\node[] (H) [right = of G, label=below: $x_{q_{n_m}}$] {};
\node[] (I) [right = of H] {};
\path (A) edge [bend right = 0] (I);
\draw[|-|,blue] (8.98,0) -- (9.6,0);
\draw[|-|,blue] (12,0) -- (12.61,0);
\Anm[0.5]{10.8}
\end{tikzpicture}
\begin{center}
$\vdots\qquad\qquad$
\end{center}
\noindent\begin{tikzpicture}[ ,auto ,node distance =1.8cm ,on grid ,
thick ,
state/.style ={ circle ,fill =black,
draw,black , text=white , minimum width =0.1 cm}]
\draw (1.8cm,-.1) -- (1.8cm, .1);   % 1 x 1.8
\draw (5.4cm,-.1) -- (5.4cm, .1);  % 3 x 1.8
\draw (9.0cm,-.1) -- (9.0cm, .1);   % 4 x 1.8
\draw (12.6cm,-.1) -- (12.6cm, .1); % 7 x 1.8
\filldraw (3.9,0) circle (0pt) node[align=left,   label=above: $\textcolor{orange}{x_{t_0}}$] {};
\node[] (A){};
\node[] (B) [right = of A, label=below: $x_{k_0+q_{n_m+1}}$ ] {};
\node[] (C) [right = of B, label=below: $z_{0}$ ] {$\bullet$};
\node[] (D) [right = of C, label=below: $x_{k_0}$] {};
\node[] (E) [right = of D, label=below: ] {};
\node[] (F) [right = of E, label=below: $x_{k_0-q_{n_m+1}}$] {\textcolor{red}{$\times$}};
%\node[] (G) [right = of F] {};
%\node[] (H) [right = of G] {};
\node[] (G) [right = of F, label=below: $z_{q_{n_m}}$ ] {$\bullet$};
\node[] (H) [right = of G, label=below: $x_{k_0+q_{n_m}}$] {};
\node[] (I) [right = of H] {};
\path (A) edge [bend right = 0] (I);
\draw[|-|,blue] (1.8,0) -- (2.6,0);
\draw[|-|,blue] (4.8,0) -- (5.4,0);
\draw[|-|,blue] (8.98,0) -- (9.6,0);
\draw[|-|,blue] (12,0) -- (12.61,0);
\Anm[0.65]{10.8}
\Anm[0.1]{3.6}
\end{tikzpicture}
\begin{center}
$\vdots\qquad\qquad$
\end{center}
%%%%%%%%%%%%%%%%%%%%%%%%%%%%%%%%%%%%%%%%%%%%%%%%%
%%%%%%%%%%%%%%%%%%%% Second two%%%%%%%%%%%%%%%%%%%%%
%%%%%%%%%%%%%%%%%%%%%%%%%%%%%%%%%%%%%%%%%%%%%%%%%
\noindent\begin{tikzpicture}[ ,auto ,node distance =1.8cm ,on grid ,
thick ,
state/.style ={ circle ,fill =black,
draw,black , text=white , minimum width =0.1 cm}]
\draw (1.8cm,-.1) -- (1.8cm, .1);   % 1 x 1.8
\draw (5.4cm,-.1) -- (5.4cm, .1);  % 3 x 1.8
\draw (9.0cm,-.1) -- (9.0cm, .1);   % 4 x 1.8
\draw (12.6cm,-.1) -- (12.6cm, .1); % 7 x 1.8
\node[] (A){};
\node[] (B) [right = of A, label=below: $x_{q_{n_m}+q_{n_m+1}-1}$ ] {};
\node[] (C) [right = of B, label=below: $z_{q_{n_m}-k_0-1}$ ] {$\bullet$};
\node[] (D) [right = of C, label=below: $x_{q_{n_m}-1}$] {};
\node[] (E) [right = of D, label=below: ] {};
\node[] (F) [right = of E, label=below: $x_{q_{n_m}-q_{n_m+1}-1}$] {\textcolor{red}{$\times$}};
%\node[] (G) [right = of F] {};
%\node[] (H) [right = of G] {};
\node[] (G) [right = of F, label=below: $z_{2q_{n_m}-k_0-1}$ ] {$\bullet$};
\node[] (H) [right = of G, label=below: $x_{2q_{n_m}-1}$] {};
\node[] (I) [right = of H] {};
\path (A) edge [bend right = 0] (I);
\draw[|-|,blue] (8.98,0) -- (9.6,0);
\draw[|-|,blue] (12,0) -- (12.61,0);
\draw[|-|,blue] (1.8,0) -- (2.6,0);
\draw[|-|,blue] (4.8,0) -- (5.4,0);
\Anm[0.75]{10.8}
\Anm[0.45]{3.6}
\end{tikzpicture}
\\\begin{tikzpicture}[ ,auto ,node distance =1.8cm ,on grid ,
thick ,
state/.style ={ circle ,fill =black,
draw,black , text=white , minimum width =0.1 cm}]
%\draw (1.8cm,-.1) -- (1.8cm, .1);   % 1 x 1.8
\draw (5.4cm,-.1) -- (5.4cm, .1);  % 3 x 1.8
\draw (9.0cm,-.1) -- (9.0cm, .1);   % 4 x 1.8
\draw (12.6cm,-.1) -- (12.6cm, .1); % 7 x 1.8
\node[] (A){};
\node[] (B) [right = of A, label=below:  ] {};
\node[] (C) [right = of B, label=below:  ] {};
\node[] (D) [right = of C, label=below: $x_{q_{n_m}}$] {};
\node[] (E) [right = of D, label=below: ] {};
\node[] (F) [right = of E, label=below: $x_{q_{n_m}-q_{n_m+1}}$] {\textcolor{red}{$\times$}};
%\node[] (G) [right = of F] {};
%\node[] (H) [right = of G] {};
\node[] (G) [right = of F, label=below: $z_{2q_{n_m}-k_0}$ ] {$\bullet$};
\node[] (H) [right = of G, label=below: $x_{2q_{n_m}}$] {};
\node[] (I) [right = of H] {};
\path (A) edge [bend right = 0] (I);
\draw[|-|,blue] (8.98,0) -- (9.6,0);
\draw[|-|,blue] (12,0) -- (12.61,0);
\Anm[0.8]{10.8}
\end{tikzpicture}
\begin{center}
$\vdots\qquad\qquad$
\end{center}
%%%%%%%%%%%%%%%%%%%%%%%%%%%%%%%%%%%%%%%%%%%%%%%%%
%%%%%%%%%%%%%%%%%%%% Last  %%%%%%%%%%%%%%%%%%%%%
%%%%%%%%%%%%%%%%%%%%%%%%%%%%%%%%%%%%%%%%%%%%%%%%%
\noindent\begin{tikzpicture}[ ,auto ,node distance =1.8cm ,on grid ,
thick ,
state/.style ={ circle ,fill =black,
draw,black , text=white , minimum width =0.1 cm}]
%\draw (1.8cm,-.1) -- (1.8cm, .1);   % 1 x 1.8
\draw (5.4cm,-.1) -- (5.4cm, .1);  % 3 x 1.8
\draw (9.0cm,-.1) -- (9.0cm, .1);   % 4 x 1.8
\draw (12.6cm,-.1) -- (12.6cm, .1); % 7 x 1.8
\node[] (A){};
\node[] (B) [right = of A, label=below:  ] {};
\node[] (C) [right = of B, label=below:  ] {};
\node[] (D) [right = of C, label=below: $x_{q_{n_m+1}-q_{n_m}+k_0-1}$] {};
\node[] (E) [right = of D, label=below: ] {};
\node[] (F) [right = of E, label=below: $x_{-q_{n_m}+k_0-1}$] {\textcolor{red}{$\times$}};
%\node[] (G) [right = of F] {};
%\node[] (H) [right = of G] {};
\node[] (G) [right = of F, label=below: $z_{q_{{n_m}+1}-1}$ ] {$\bullet$};
\node[] (H) [right = of G, label=below: $x_{q_{n_m+1}+k_0-1}$] {};
\node[] (I) [right = of H] {};
\path (A) edge [bend right = 0] (I);
\draw[|-|,blue] (8.98,0) -- (9.6,0);
\draw[|-|,blue] (12,0) -- (12.61,0);
\Anm[0.9]{10.8}
\end{tikzpicture}
\begin{center}
$\vdots\qquad\qquad$
\end{center}
\noindent\begin{tikzpicture}[ ,auto ,node distance =1.8cm ,on grid ,
thick ,
state/.style ={ circle ,fill =black,
draw,black , text=white , minimum width =0.1 cm}]
%\draw (1.8cm,-.1) -- (1.8cm, .1);   % 1 x 1.8
\draw (5.4cm,-.1) -- (5.4cm, .1);  % 3 x 1.8
\draw (9.0cm,-.1) -- (9.0cm, .1);   % 4 x 1.8
\draw (12.6cm,-.1) -- (12.6cm, .1); % 7 x 1.8
\node[] (A){};
\node[] (B) [right = of A, label=below:  ] {};
\node[] (C) [right = of B, label=below:  ] {};
\node[] (D) [right = of C, label=below: $x_{q_{n_m+1}-1}$] {};
\node[] (E) [right = of D, label=below: ] {};
\node[] (F) [right = of E, label=below: $x_{-1}$] {\textcolor{red}{$\times$}};
%\node[] (G) [right = of F] {};
%\node[] (H) [right = of G] {};
\node[] (G) [right = of F, label=below:  ] {};
\node[] (H) [right = of G, label=below: $x_{q_{n_m+1}+q_{n_m}-1}$] {};
\node[] (I) [right = of H] {};
\path (A) edge [bend right = 0] (I);
\end{tikzpicture}\\
%\end{center}
%%%%%%%%%%%%%%%%%%%%%%%%%%%%%%%%%%%%%%%%%%%%%%%%%
%%%%%%%%%%%%%%%%%%%% Finished %%%%%%%%%%%%%%%%%%%%%
%%%%%%%%%%%%%%%%%%%%%%%%%%%%%%%%%%%%%%%%%%%%%%%%%
\textbf{The intervals $A_{k}^{(n_{m})}$ in red in case  II.}\\
%%%%%%%%%%%%%%%%%%    Picture 3.1 %%%%%%%%%%%%%%%%%%%%%%
%%%%%%%%%%%%%%%%%%%% First two%%%%%%%%%%%%%%%%%%%%%%%%
%%%%%%%%%%%%%%%%%%%%%%%%%%%%%%%%%%%%%%%%%%%%%%%%%
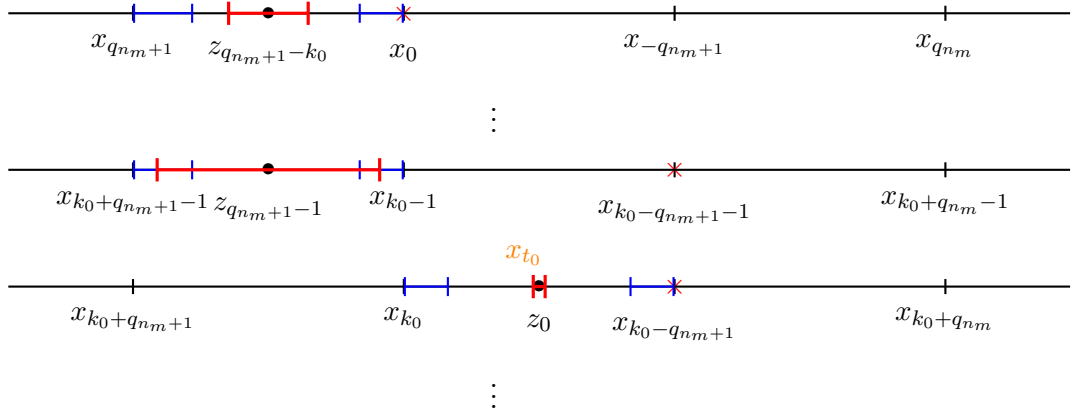
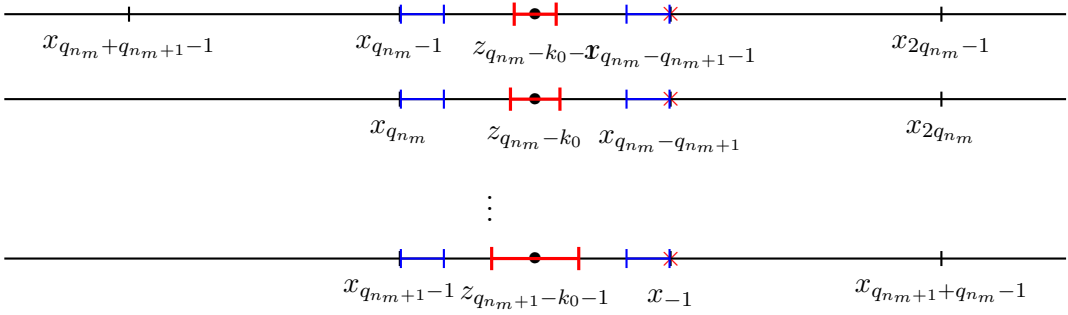
\noindent\begin {tikzpicture}[ ,auto ,node distance =1.8cm ,on grid ,
thick ,
state/.style ={ circle ,fill =black,
draw,black , text=white , minimum width =0.1 cm}]
\draw (1.8cm,-.1) -- (1.8cm, .1);   % 1 x 1.8
\draw (5.4cm,-.1) -- (5.4cm, .1);  % 3 x 1.8
\draw (9.0cm,-.1) -- (9.0cm, .1);   % 4 x 1.8
\draw (12.6cm,-.1) -- (12.6cm, .1); % 7 x 1.8
\node[] (A){};
\node[] (B) [right = of A, label=below: $x_{q_{n_m+1}}$ ] {};
\node[] (C) [right = of B, label=below:$z_{q_{n_m+1}-k_0}$ ] {$\bullet$};
\node[] (D) [right = of C, label=below: $x_{0}$] {\textcolor{red}{$\times$}};
\node[] (E) [right = of D, label=below: ] {};
\node[] (F) [right = of E, label=below: $x_{-q_{n_m+1}}$] {};
%\node[] (G) [right = of F] {};
%\node[] (H) [right = of G] {};
\node[] (G) [right = of F, label=below: ] {};
\node[] (H) [right = of G, label=below: $x_{q_{n_m}}$] {};
\node[] (I) [right = of H] {};
\path (A) edge [bend right = 0] (I);
\draw[|-|,blue] (1.8,0) -- (2.6,0);
\draw[|-|,blue] (4.8,0) -- (5.4,0);
\Anm[0.55]{3.6}
\end{tikzpicture}
\begin{center}
$\vdots\qquad\qquad$
\end{center}
\noindent\begin {tikzpicture}[ ,auto ,node distance =1.8cm ,on grid ,
thick ,
state/.style ={ circle ,fill =black,
draw,black , text=white , minimum width =0.1 cm}]
\draw (1.8cm,-.1) -- (1.8cm, .1);   % 1 x 1.8
\draw (5.4cm,-.1) -- (5.4cm, .1);  % 3 x 1.8
\draw (9.0cm,-.1) -- (9.0cm, .1);   % 4 x 1.8
\draw (12.6cm,-.1) -- (12.6cm, .1); % 7 x 1.8
\node[] (A){};
\node[] (B) [right = of A, label=below: $x_{k_0+q_{n_m+1}-1}$ ] {};
\node[] (C) [right = of B, label=below:$z_{q_{n_m+1}-1}$  ] {$\bullet$};
\node[] (D) [right = of C, label=below: $x_{k_0-1}$] {};
\node[] (E) [right = of D, label=below: ] {};
\node[] (F) [right = of E, label=below: $x_{k_0-q_{n_m+1}-1}$] {\textcolor{red}{$\times$}};
%\node[] (G) [right = of F] {};
%\node[] (H) [right = of G] {};
\node[] (G) [right = of F, label=below:  ] {};
\node[] (H) [right = of G, label=below: $x_{k_0+q_{n_m}-1}$] {};
\node[] (I) [right = of H] {};
\path (A) edge [bend right = 0] (I);
\draw[|-|,blue] (1.8,0) -- (2.6,0);
\draw[|-|,blue] (4.8,0) -- (5.4,0);
\Anm[1.5]{3.6}
\end{tikzpicture}
\\\begin{tikzpicture}[ ,auto ,node distance =1.8cm ,on grid ,
thick ,
state/.style ={ circle ,fill =black,
draw,black , text=white , minimum width =0.1 cm}]
\draw (1.8cm,-.1) -- (1.8cm, .1);   % 1 x 1.8
\draw (5.4cm,-.1) -- (5.4cm, .1);  % 3 x 1.8
\draw (9.0cm,-.1) -- (9.0cm, .1);   % 4 x 1.8
\draw (12.6cm,-.1) -- (12.6cm, .1); % 7 x 1.8
\filldraw (7.0,0) circle (0pt) node[align=left,   label=above: $\textcolor{orange}{x_{t_0}}$] {};
\node[] (A){};
\node[] (B) [right = of A, label=below: $x_{k_0+q_{n_m+1}}$ ] {};
\node[] (C) [right = of B, label=below:  ] {};
\node[] (D) [right = of C, label=below: $x_{k_0}$] {};
\node[] (E) [right = of D, label=below:$z_{0}$ ] {$\bullet$};
\node[] (F) [right = of E, label=below: $x_{k_0-q_{n_m+1}}$] {\textcolor{red}{$\times$}};
%\node[] (G) [right = of F] {};
%\node[] (H) [right = of G] {};
\node[] (G) [right = of F, label=below: ] {};
\node[] (H) [right = of G, label=below: $x_{k_0+q_{n_m}}$] {};
\node[] (I) [right = of H] {};
\path (A) edge [bend right = 0] (I);
\draw[|-|,blue] (5.4,0) -- (6,0);
\draw[|-|,blue] (8.4,0) -- (9,0);
\Anm[0.1]{7.2}
\end{tikzpicture}
\begin{center}
$\vdots\qquad\qquad$
\end{center}
%%%%%%%%%%%%%%%%%%%%%%%%%%%%%%%%%%%%%%%%%%%%%%%%%
%%%%%%%%%%%%%%%%%%%% Second two%%%%%%%%%%%%%%%%%%%%%
%%%%%%%%%%%%%%%%%%%%%%%%%%%%%%%%%%%%%%%%%%%%%%%%%
\noindent\begin{tikzpicture}[ ,auto ,node distance =1.8cm ,on grid ,
thick ,
state/.style ={ circle ,fill =black,
draw,black , text=white , minimum width =0.1 cm}]
\draw (1.8cm,-.1) -- (1.8cm, .1);   % 1 x 1.8
\draw (5.4cm,-.1) -- (5.4cm, .1);  % 3 x 1.8
\draw (9.0cm,-.1) -- (9.0cm, .1);   % 4 x 1.8
\draw (12.6cm,-.1) -- (12.6cm, .1); % 7 x 1.8
\node[] (A){};
\node[] (B) [right = of A, label=below: $x_{q_{n_m}+q_{n_m+1}-1}$ ] {};
\node[] (C) [right = of B, label=below: ] {};
\node[] (D) [right = of C, label=below: $x_{q_{n_m}-1}$] {};
\node[] (E) [right = of D, label=below: $z_{q_{n_m}-k_0-1}$] {$\bullet$};
\node[] (F) [right = of E, label=below: $x_{q_{n_m}-q_{n_m+1}-1}$] {\textcolor{red}{$\times$}};
%\node[] (G) [right = of F] {};
%\node[] (H) [right = of G] {};
\node[] (G) [right = of F, label=below:  ] {};
\node[] (H) [right = of G, label=below: $x_{2q_{n_m}-1}$] {};
\node[] (I) [right = of H] {};
\path (A) edge [bend right = 0] (I);
\draw[|-|,blue] (5.4,0) -- (6,0);
\draw[|-|,blue] (8.4,0) -- (9,0);
\Anm[0.3]{7.2}
\end{tikzpicture}
\\\begin{tikzpicture}[ ,auto ,node distance =1.8cm ,on grid ,
thick ,
state/.style ={ circle ,fill =black,
draw,black , text=white , minimum width =0.1 cm}]
%\draw (1.8cm,-.1) -- (1.8cm, .1);   % 1 x 1.8
\draw (5.4cm,-.1) -- (5.4cm, .1);  % 3 x 1.8
\draw (9.0cm,-.1) -- (9.0cm, .1);   % 4 x 1.8
\draw (12.6cm,-.1) -- (12.6cm, .1); % 7 x 1.8
\node[] (A){};
\node[] (B) [right = of A, label=below:  ] {};
\node[] (C) [right = of B, label=below:  ] {};
\node[] (D) [right = of C, label=below: $x_{q_{n_m}}$] {};
\node[] (E) [right = of D, label=below: $z_{q_{n_m}-k_0}$] {$\bullet$};
\node[] (F) [right = of E, label=below: $x_{q_{n_m}-q_{n_m+1}}$] {\textcolor{red}{$\times$}};
%\node[] (G) [right = of F] {};
%\node[] (H) [right = of G] {};
\node[] (G) [right = of F, label=below:  ] {};
\node[] (H) [right = of G, label=below: $x_{2q_{n_m}}$] {};
\node[] (I) [right = of H] {};
\path (A) edge [bend right = 0] (I);
\draw[|-|,blue] (5.4,0) -- (6,0);
\draw[|-|,blue] (8.4,0) -- (9,0);
\Anm[0.35]{7.2}
\end{tikzpicture}
\begin{center}
$\vdots\qquad\qquad$
\end{center}
%%%%%%%%%%%%%%%%%%%%%%%%%%%%%%%%%%%%%%%%%%%%%%%%%
%%%%%%%%%%%%%%%%%%%% Last  %%%%%%%%%%%%%%%%%%%%%
%%%%%%%%%%%%%%%%%%%%%%%%%%%%%%%%%%%%%%%%%%%%%%%%%
\noindent\begin{tikzpicture}[ ,auto ,node distance =1.8cm ,on grid ,
thick ,
state/.style ={ circle ,fill =black,
draw,black , text=white , minimum width =0.1 cm}]
%\draw (1.8cm,-.1) -- (1.8cm, .1);   % 1 x 1.8
\draw (5.4cm,-.1) -- (5.4cm, .1);  % 3 x 1.8
\draw (9.0cm,-.1) -- (9.0cm, .1);   % 4 x 1.8
\draw (12.6cm,-.1) -- (12.6cm, .1); % 7 x 1.8
\node[] (A){};
\node[] (B) [right = of A, label=below:  ] {};
\node[] (C) [right = of B, label=below:  ] {};
\node[] (D) [right = of C, label=below: $x_{q_{n_m+1}-1}$] {};
\node[] (E) [right = of D, label=below: $z_{q_{n_m+1}-k_0-1}$] {$\bullet$};
\node[] (F) [right = of E, label=below: $x_{-1}$] {\textcolor{red}{$\times$}};
%\node[] (G) [right = of F] {};
%\node[] (H) [right = of G] {};
\node[] (G) [right = of F, label=below: ] {};
\node[] (H) [right = of G, label=below: $x_{q_{n_m+1}+q_{n_m}-1}$] {};
\node[] (I) [right = of H] {};
\path (A) edge [bend right = 0] (I);
\draw[|-|,blue] (5.4,0) -- (6,0);
\draw[|-|,blue] (8.4,0) -- (9,0);
\Anm[0.6]{7.2}
\end{tikzpicture}\\
\textbf{The intervals $A_{k}^{n_{m}}$ in red in case III.}\\
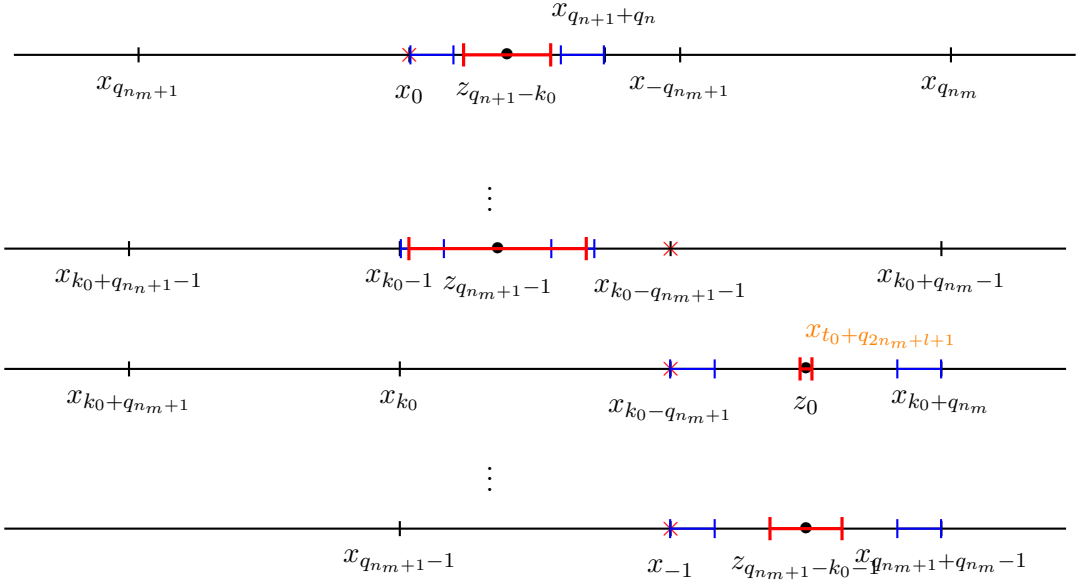
\vspace{0.5cm}
\noindent\begin {tikzpicture}[ ,auto ,node distance =1.8cm ,on grid ,
thick ,
state/.style ={ circle ,fill =black,
draw,black , text=white , minimum width =0.1 cm}]
\draw (1.8cm,-.1) -- (1.8cm, .1);   % 1 x 1.8
\draw (5.4cm,-.1) -- (5.4cm, .1);  % 3 x 1.8
\draw (9.0cm,-.1) -- (9.0cm, .1);   % 4 x 1.8
\draw (12.6cm,-.1) -- (12.6cm, .1); % 7 x 1.8
\draw (8cm,-.1) -- (8cm, .1) node [ label=above: $x_{q_{n+1}+q_n}$] {}; % 7 x 1.8
\filldraw (6.7,0) circle (1pt) node[align=left,   label=below: $z_{q_{n+1}-k_0}$] {$\bullet$};
\node[] (A){};
\node[] (B) [right = of A, label=below: $x_{q_{n_m+1}}$ ] {};
\node[] (C) [right = of B, label=below: ] {};
\node[] (D) [right = of C, label=below: $x_{0}$] {\textcolor{red}{$\times$}};
\node[] (E) [right = of D, label=below: ] {};
\node[] (F) [right = of E, label=below: $x_{-q_{n_m+1}}$] {};
%\node[] (G) [right = of F] {};
%\node[] (H) [right = of G] {};
\node[] (G) [right = of F, label=below: ] {};
\node[] (H) [right = of G, label=below: $x_{q_{n_m}}$] {};
\node[] (I) [right = of H] {};
\path (A) edge [bend right = 0] (I);
\draw[|-|,blue] (5.4,0) -- (6,0);
\draw[|-|,blue] (7.4,0) -- (8,0);
\Anm[0.6]{6.7}
\end{tikzpicture}
\begin{center}
$\vdots\qquad\qquad$
\end{center}
\noindent\begin {tikzpicture}[ ,auto ,node distance =1.8cm ,on grid ,
thick ,
state/.style ={ circle ,fill =black,
draw,black , text=white , minimum width =0.1 cm}]
\draw (1.8cm,-.1) -- (1.8cm, .1);   % 1 x 1.8
\draw (5.4cm,-.1) -- (5.4cm, .1);  % 3 x 1.8
\draw (9.0cm,-.1) -- (9.0cm, .1);   % 4 x 1.8
\draw (12.6cm,-.1) -- (12.6cm, .1); % 7 x 1.8
\filldraw (6.7,0) circle (1pt) node[align=left,   label=below: $z_{q_{n_m+1}-1}$  ] {$\bullet$};
\node[] (A){};
\node[] (B) [right = of A, label=below: $x_{k_0+q_{n_n+1}-1}$ ] {};
\node[] (C) [right = of B, label=below:] {};
\node[] (D) [right = of C, label=below: $x_{k_0-1}$] {};
\node[] (E) [right = of D, label=below: ] {};
\node[] (F) [right = of E, label=below: $x_{k_0-q_{n_m+1}-1}$] {\textcolor{red}{$\times$}};
%\node[] (G) [right = of F] {};
%\node[] (H) [right = of G] {};
\node[] (G) [right = of F, label=below:  ] {};
\node[] (H) [right = of G, label=below: $x_{k_0+q_{n_m}-1}$] {};
\node[] (I) [right = of H] {};
\path (A) edge [bend right = 0] (I);
\draw[|-|,blue] (5.4,0) -- (6,0);
\draw[|-|,blue] (7.4,0) -- (8,0);
\Anm[1.2]{6.7}
\end{tikzpicture}
\\\begin{tikzpicture}[ ,auto ,node distance =1.8cm ,on grid ,
thick ,
state/.style ={ circle ,fill =black,
draw,black , text=white , minimum width =0.1 cm}]
\draw (1.8cm,-.1) -- (1.8cm, .1);   % 1 x 1.8
\draw (5.4cm,-.1) -- (5.4cm, .1);  % 3 x 1.8
\draw (9.0cm,-.1) -- (9.0cm, .1);   % 4 x 1.8
\draw (12.6cm,-.1) -- (12.6cm, .1); % 7 x 1.8
\filldraw (11.8,0) circle (0pt) node[align=left,   label=above: $\textcolor{orange}{x_{t_0+q_{2n_m+l+1}}}$] {};
\node[] (A){};
\node[] (B) [right = of A, label=below: $x_{k_0+q_{n_m+1}}$ ] {};
\node[] (C) [right = of B, label=below:  ] {};
\node[] (D) [right = of C, label=below: $x_{k_0}$] {};
\node[] (E) [right = of D, label=below:] {};
\node[] (F) [right = of E, label=below: $x_{k_0-q_{n_m+1}}$] {\textcolor{red}{$\times$}};
%\node[] (G) [right = of F] {};
%\node[] (H) [right = of G] {};
\node[] (G) [right = of F, label=below: $z_{0}$ ] {$\bullet$};
\node[] (H) [right = of G, label=below: $x_{k_0+q_{n_m}}$] {};
\node[] (I) [right = of H] {};
\path (A) edge [bend right = 0] (I);
\draw[|-|,blue] (8.98,0) -- (9.6,0);
\draw[|-|,blue] (12,0) -- (12.61,0);
\Anm[0.1]{10.8}
\end{tikzpicture}
\begin{center}
$\vdots\qquad\qquad$
\end{center}
%%%%%%%%%%%%%%%%%%%%%%%%%%%%%%%%%%%%%%%%%%%%%%%%%
%%%%%%%%%%%%%%%%%%%% Last  %%%%%%%%%%%%%%%%%%%%%
%%%%%%%%%%%%%%%%%%%%%%%%%%%%%%%%%%%%%%%%%%%%%%%%%
\noindent\begin{tikzpicture}[ ,auto ,node distance =1.8cm ,on grid ,
thick ,
state/.style ={ circle ,fill =black,
draw,black , text=white , minimum width =0.1 cm}]
%\draw (1.8cm,-.1) -- (1.8cm, .1);   % 1 x 1.8
\draw (5.4cm,-.1) -- (5.4cm, .1);  % 3 x 1.8
\draw (9.0cm,-.1) -- (9.0cm, .1);   % 4 x 1.8
\draw (12.6cm,-.1) -- (12.6cm, .1); % 7 x 1.8
\node[] (A){};
\node[] (B) [right = of A, label=below:  ] {};
\node[] (C) [right = of B, label=below:  ] {};
\node[] (D) [right = of C, label=below: $x_{q_{n_m+1}-1}$] {};
\node[] (E) [right = of D, label=below:] {};
\node[] (F) [right = of E, label=below: $x_{-1}$] {\textcolor{red}{$\times$}};
%\node[] (G) [right = of F] {};
%\node[] (H) [right = of G] {};
\node[] (G) [right = of F, label=below: $z_{q_{n_m+1}-k_0-1}$] {$\bullet$};
\node[] (H) [right = of G, label=below: $x_{q_{n_m+1}+q_{n_m}-1}$] {};
\node[] (I) [right = of H] {};
\path (A) edge [bend right = 0] (I);
\draw[|-|,blue] (8.98,0) -- (9.6,0);
\draw[|-|,blue] (12,0) -- (12.61,0);
\Anm[0.5]{10.8}
\end{tikzpicture}

\begin{lemm}\label{Bnm}
Let $T$ be a  circle homeomorphism, $(\xi_n)_{n=1}^\infty$ be a sequence of independent
 random variables  satisfying the conditions of Theorem \ref{main} and let $\theta_{+}=(1+e^{v})^{-\frac{1}{2}}<1.$\\
Suppose that the sequence $\left\{\sigma^2_{q_{n+1}},\,n\geq1 \right\} $  fulfills the condition
$$\lim\limits_{n\to \infty}\frac{q_{n+1}\sigma^2_{q_{n+1}}}{\theta_{+}^{2n+4}}=0.$$
 Then the probability of the event
\begin{equation}\label{B_n_m} B_{n_m}:=\{\bar{z}_1\in A^{(n_m)}_{1}, \bar{z}_2\in A^{(n_m)}_{2},\,...,\,\bar{z}_{q_{n_m+1}-1}\in A^{(n_m)}_{q_{n_m+1}-1}\}
\end{equation}
 tends to $1$ as $m$ tends to infinity, where
  $$\bar{z}_k=T(\bar{z}_{k-1})+\sigma_{q_{n_m+1}}\xi_k,\,1\leq k\leq q_{n_m+1}-1, \bar{z}_0=z_0.$$
\end{lemm}
\proof It is clear that
\begin{equation}\begin{split}\label{prob0}
\mathbf{P}(B_{n_m})&=  \mathbf{P}\left(\bigcap\limits_{k=1}^{q_{n_m+1}-1}\{\bar{z}_k\in A^{n_m}_k \}\right)\\
&=\mathbf{P}\left(\bar{z}_1\in A^{n_m}_1\right)\prod\limits_{k=2}^{q_{n_m+1}-1}\mathbf{P}\left(\{\bar{z}_k\in A^{n_m}_k \}|\bigcap\limits_{i=1}^{k-1}\{\bar{z}_i\in A^{n_m}_i \}\right).\end{split}\end{equation}
Next we estimate the factors of the last product.\\
Applying Chebeshev's inequality we obtain
$$\mathbf{P}(\bar{z}_1\in A^{n_m}_1)=\mathbf{P}(T(z_0)+\sigma_{q_{n_m+1}}\xi_1\in A^{-}_1\cup TA^{n_m}_0\cup A^{+}_1)\geq$$
$$\geq \mathbf{P}(|\sigma_{q_{n_m+1}}\xi_1|\leq \min\{|A^{-}_1|,|A^{+}_1|\})\geq 1-\frac{\sigma_{q_{n_m+1}}^2\mathrm{Var}\xi_1}{\min\{|A^{-}_1|^2,|A^{+}_1|^2\}}.$$
Since $A^{-}_1$ and $A^{+}_1$ are elements  of $P_{2n_m+\textcolor{blue}{7}}$, we have $\min\{|A^{-}_1|^2,|A^{+}_1|^2\}\geq C \theta_{+}^{2n_m+\textcolor{blue}{7}}.$
Using the last bound we get
\begin{equation}\label{Prob1}\mathbf{P}(\bar{z}_1\in A^{n_m}_1)\geq 1-\frac{C\sigma^2_{q_{n_m+1}}}{\theta_{+}^{2n_m+\textcolor{blue}{7}}}.
\end{equation}
For the other  $k,\,\,2\leq k<q_{n_m+1}$, analogyously  we have
$$\mathbf{P}(\{\bar{z}_k\in A^{n_m}_k \}/\bigcap\limits_{i=1}^{k-1}\{\bar{z}_i\in A^{n_m}_i \})=$$
$$=\mathbf{P}\left(\{T(\bar{z}_k)+\sigma_{q_{n_m+1}}\xi_k\in A^{-}_k\cup TA^{n_m}_{k-1}\cup A^{+}_k\}/\{\bar{z}_{k-1}\in A^{n_m}_{k-1}\} \right)\geq$$
$$\geq \mathbf{P}(|\sigma_{q_{n_m+1}}\xi_k|\leq \min\{|A^{-}_k|,|A^{+}_k|\})\geq 1-\frac{\sigma_{q_{n_m+1}}^2\mathrm{Var}\;\xi_k}{\min\{|A^{-}_k|^2,|A^{+}_k|^2\}},$$
and therefore as before
\begin{equation}\label{Prob2}\mathbf{P}(\{\bar{z}_k\in A^{n_m}_k \}/\bigcap\limits_{i=1}^{k-1}\{\bar{z}_i\in A^{n_m}_i \})\geq 1-\frac{C\sigma^2_{q_{n_m+1}}}{\theta_{+}^{2n_m+\textcolor{blue}{7}}}.\end{equation}

Summarizing  (\ref{prob0}), (\ref{Prob1}) and (\ref{Prob2}) we get:
\begin{equation} \label{bound1}
\mathbf{P}(B_{n_m})=\mathbf{P}(\bigcap\limits_{k=1}^{q_{n_m+1}-1}\{\bar{z}_k\in A_{k}^{(n_m)}(z_k) \})\geq \left(1-\frac{C\sigma_{q_{n_m+1}}^2}{\theta_{+}^{2n_m+\textcolor{blue}{7}}}\right)^{q_{n_m+1}}.
\end{equation}
By assumption
$$\lim\limits_{m\to \infty}\frac{q_{n_m+1}\sigma^2_{q_{n_m+1}}}{\theta_{+}^{2n_m+\textcolor{blue}{7}}}=0.$$

This together with (\ref{Prob2})
implies that
$$\mathbf{P}(B_{n_m})\to 1.$$
as ${m\to\infty}.$
\endproof
 We fix a point $z_0\in S^1\setminus\{x_b\} $.   Assuming $\bar{z}_j\in A_j^{n_m},\,1\leq j\leq q_{n_m+1}-1$
we can use the Taylor expansion for the random process $\bar{z}_{q_{n_m+1}}(z_0,\sigma_{q_{n_m+1}-1})$
\begin{equation}\label{Proc.}
\bar{z}_{q_{n_m+1}-1}(z_0, \sigma_{q_{n_m+1}-1})=T^{q_{n_m+1}-1}(z_0)+\sigma L_{q_{n_m+1}-1}(z_0)+\sigma_{q_{n_m+1}-1}^2 Q_{q_{n_m+1}-1}(z_0,\sigma_{q_{n_m+1}-1}),
\end{equation}
where the linear term $L_{q_{n_m+1}}$ is the sum of independent random variables
defined by (\ref{Ln}).\\
 To the  linear  process
\begin{equation}\label{y} y_n(z_0,\sigma)=T^n(z_0)+\sigma L_n(z_0)\end{equation}
 we can apply the  straightforward extension of the central limit theorem
 as proved in \cite{Llave} (see  Lemma 3.1 ).
\begin{lemm} \label{5.3} (see  \cite{Llave}). Let $T\in C^2(S^1\setminus\{x_b\})$ be a circle homeomorphism with a break point $x_b$
 and let $\{\xi_n\}$  be a sequence of independent random variables with $p>2$ moments satisfying
the  conditions (\ref{Cond1}) and (\ref{Cond2}). Assume condition (\ref{maincon}) holds for some point
 $z_{0}\in S^1\setminus \{\mathbb{O}_T(x_b)\}$ and some increasing sequence $\left\{n_k\right\}$ of positive integers,
 then
$$l_{n_k}(z_{0})=\frac{L_{n_k}(z_{0})}{\sqrt{\mathrm{var}\,L_{n_k}(z_{0})}}$$
converges in distribution to the standard Gaussian as $k\to \infty$. Moreover,
there is a universal constant $C$ such that
\begin{equation}\label{BerreEsseyn}\sup\limits_{t\in R} |P(l_{n_k}(z_{0})<t)-\Phi(t)|\leq C \frac{\Lambda_{\min(p,3)}(x,n_k)}{(\Lambda_2(x,n_k))^{\min(p,3)/2}}.\end{equation}
\end{lemm}

Condition  (\ref{maincon}) holds indeed in our case, namely
using (\ref{ln}) , Lemma \ref{Eigen.1} and $p>2$ we get
\begin{equation}
\lim\limits_{n\to\infty}\frac{\Lambda_p(z_{0},q_{n})}{(\Lambda_2(z_{0},q_{n}))^{p/2}}=
\lim\limits_{n\to\infty}\frac{|I^{(n)}_{0}({x_{b}};z_{0})|^p\lambda_{-p}^n}{(|I^{(n)}_{0}({x_{b}};z_{0})|^2\lambda_{-2}^n)^{p/2}}=0.
\end{equation}

Next we treat the nonlinear part of the process (\ref{Proc.}). For this consider the process

$$\omega_{q_{n_m+1}-1}=\frac{\bar{z}_{q_{n_m+1}-1}-T^{q_{n_m+1}-1}(z_0)}{\sigma_{q_{n_m+1}-1}\sqrt{\mathrm{var}(L_{q_{n_m+1}-1}(z_0))}}=
\frac{L_{q_{n_m+1}-1}(z_0)}{\sqrt{\mathrm{var}(L_{q_{n_m+1}-1}(z_0))}}+\sigma_{q_{n_m+1}-1}
\frac{Q_{q_{n_m+1}-1}(z_0,\sigma_{q_{n_m+1}-1})}{\sqrt{\mathrm{var}(L_{q_{n_m+1}-1}(z_0))}}$$

Lemma \ref{5.3} implies that
 $\frac{L_{q_{n_m+1}-1}(z_0)}{\sqrt{\mathrm{var}(L_{q_{n_m+1}-1}(z_0))}}$ converges weekly (in distribution) to the standard Gaussian.

Next we will  show that the random process   $\sigma_{q_{n_m+1}-1}
\frac{Q_{q_{n_m+1}-1}(z_0,\sigma_{q_{n_m+1}-1})}{\sqrt{\mathrm{var}(L_{q_{n_m+1}-1}(z_0))}}$ converges to $0$ in probability. For this we introduce the following constants
 $$K_1=\sup\limits_{x\in S^1}\frac{1}{|T'(x)|},\,\,
K_2=\exp\left(K_1\sup\limits_{x\in S^1}|T_{+}''(x)|\right),\,\,
 K=K_1\cdot K_2\sup\limits_{x\in S^1}|T_{+}''(x)|.$$
\begin{lemm}\label{analogously} Suppose a circle map $T$ satisfies the conditions of Theorem \ref{main} and
the sequence $\sigma_{n}$  satisfies relation (15).
 Let $D_{n_{m}}$ be the event
\begin{equation} \label{CondQn} D_{n_m}=\big\{K\cdot\sigma_{q_{n_m+1}}\max\limits_{1\leq i\leq q_{n_m+1}}|\xi_i| \left(\hat{\Lambda}(z_0,q_{n_m+1})\right)^2<\frac{1}{2}\big\}.
\end{equation}
Under the condition  of the event $B_{n_m}$ in  (\ref{B_n_m}) the following inequality holds
$$\mathbf{P}\left(\left|\frac{\sigma_{q_{n_m+1}} Q_{q_{n_m+1}}(z_0, \sigma_{q_{n_m+1}})}{\sqrt{\mathrm{var} L_{q_{n_m+1}}(z_0)}}\mathbf{1}_{D_{n_m}}\right|>\varepsilon\right)\leq \left(\frac{2K\sigma_{q_{n_m+1}}\left( \hat{\Lambda}(x_0,q_{n+1})\right)^3 \left(E(\max\limits_{1\leq i\leq q_{n_m+1}}|\xi_i|^p)\right)^{\frac{2}{p}}}{\varepsilon\sqrt{\Lambda_{2}(z_0, q_{n_m+1})}}\right)^{\frac{p}{2}}.$$
\end{lemm}
\proof
In a first step we estimate  $|\sigma^2_{q_{n_m+1}} Q_{q_{n_m+1}}(z_0, \sigma_{q_{n_m+1}})|$.\\
Using (\ref{Ln}) we obtain  the recurrence  relation
\begin{equation}\label{Ln1}L_{k+1}(z_0)=L_{k}(z_0) T'(z_k)+ \xi_{k+1}.\end{equation}
Let $\bar{z}_k\in A^{n_m}_k,\,0\leq k<q_{n_m+1}$,  then we have
\begin{equation}\label{barx}
\begin{split}\bar{z}_{k+1}(z_0,\sigma_{q_{n_m+1}})&=T(\bar{z}_k)+\sigma_{q_{n_m+1}} \xi_{k+1}\\
&=T(z_k+\sigma_{q_{n_m+1}}L_k(z_0)+\sigma_{q_{n_m+1}}^2Q_{k}(z_0))+\sigma_{q_{n_m+1}} \xi_{k+1}\\
&=z_{k+1}+T'(\hat{z}_k)(\sigma_{q_{n_m+1}}L_k(z_0)+\sigma_{q_{n_m+1}}^2Q_{k}(z_0))+\sigma_{q_{n_m+1}} \xi_{k+1},
\end{split}
\end{equation}
where $|\hat{z}_k-z_k|\leq |\sigma_{q_{n_m+1}}L_k(z_0)+\sigma_{q_{n_m+1}}^2Q_{k}(z_0)|=|\bar{z}_k-z_k|\leq |A^{n_m}_k|.$\\
On the other hand
$$\bar{z}_{k+1}(z_0,\sigma_{q_{n_m+1}})=z_{k+1}+\sigma_{q_{n_m+1}}L_{k+1}(z_0)+\sigma_{q_{n_m+1}}^2Q_{k+1}(z_0).$$
 The last relation together with (\ref{barx}) and (\ref{Ln1}) implies that
$$\sigma_{q_{n_m+1}}^2Q_{k+1}(z_0)=\sigma_{q_{n_m+1}}L_{k}(z_0)(T'(\hat{z}_k)-T'(z_k))+ \sigma_{q_{n_m+1}}^2Q_k(z_0)T'(\hat{z}_k).$$
Iterating the last recurrence relation we obtain
\begin{equation}\label{Q_n}
\left|\sigma_{q_{n_m+1}}^2Q_{k+1}(z_0)\right|=\sigma_{q_{n_m+1}}\sum\limits_{i=1}^{k}\frac{|T'(\hat{z}_i)-T'(z_i)|}{T'(\hat{z}_i)}|L_i(z_0)|\prod\limits_{s=i}^{k}T'(\hat{z}_s),
\end{equation}
 Next we estimate the right hand side of (\ref{Q_n}).\\
For $1\leq i\leq k$ one finds
\begin{equation}\label{supT''}|T'(\hat{z}_i)-T'(z_i)|\leq \sup\limits_{x\in S^1}|T_{+}''(z)||\hat{z}_i-z_i|\leq \sup\limits_{z\in S^1}|T_{+}''(z)||\sigma_{q_{n_m+1}}L_i(z_0)+\sigma_{q_{n_m+1}}^2Q_{i}(z_0)|. \end{equation}
But
\begin{equation}\frac{1}{|T'(\hat{z}_i)|}\leq K_1.\end{equation}
Therefore
\begin{equation}\label{prod}
\begin{split}\prod\limits_{s=i}^{k}T'(\hat{z}_s)&=\prod\limits_{s=i}^{k}T'(z_s)\prod\limits_{s=i}^{k}\frac{T'(\hat{z}_s)}{T'(z_s)}\\
&=\prod\limits_{s=i}^{k}T'(z_s)\prod\limits_{s=i}^{k}\left(1+\frac{T'(\hat{z}_s)-T'(z_s)}{T'(z_s)}\right)\\
&\leq \prod\limits_{s=i}^{k}T'(z_s)\prod\limits_{s=i}^{k}\left(1+K_1\sup\limits_{z\in S^1}T_{+}''(z)|\hat{z}_i-z_i|\right)\\
&\leq\exp\left(\sum\limits_{s=i}^{k}K_1\sup\limits_{z\in S^1}T_{+}''(z)|\hat{z}_i-z_i|\right)\prod\limits_{s=i}^{k}T'(z_s)\\
&\leq\exp\left(K_1\sup\limits_{z\in S^1}|T_{+}''(z)|\right)\prod\limits_{s=i}^{k}T'(z_s)\\
&= K_2 \prod\limits_{s=i}^{k}T'(z_s).
\end{split}
\end{equation}
Using (\ref{Q_n}) and (\ref{supT''})-(\ref{prod}) we have for $0\leq k<q_{n_m+1}$
\begin{equation}
\begin{split}\left|\sigma^2_{q_{n_m+1}}Q_{k+1}(z_0)\right|&\leq K_1K_2\sup\limits_{z\in S^1}|T_{+}''(z)| \sigma_{q_{n_m+1}} \sum\limits_{i=1}^{k}|L_i(z_0)||\sigma_{q_{n_m+1}}L_i(z_0)+\sigma_{q_{n_m+1}}^2Q_{i}(z_0)|\prod\limits_{s=i}^{k}T'(z_s)\\
&\leq K\sigma_{q_{n_m+1}}^2 \sum\limits_{i=1}^{k}|L_i(z_0)|^2\prod\limits_{s=i}^{k}T'(z_s)
+K\sigma_{q_{n_m+1}} \sum\limits_{i=1}^{k}|\sigma_{q_{n_m+1}}^2Q_{i}(z_0)||L_i(z_0)|\prod\limits_{s=i}^{k}T'(z_s)\\
&\leq K\sigma_{q_{n_m+1}}^2\left(\max\limits_{1\leq i\leq k}|\xi_i| \hat{\Lambda}(z_0,k)\right)^2\sum\limits_{i=1}^{k}\prod\limits_{s=i}^{k}T'(z_s)+\\
&+K\sigma_{q_{n_m+1}}\max\limits_{1\leq i\leq k}|\sigma_{q_{n_m+1}}^2Q_{i}(z_0)|\max\limits_{1\leq i\leq k}|\xi_i| \hat{\Lambda}(z_0,k)\sum\limits_{i=1}^{k}\prod\limits_{s=i}^{k}T'(z_s).
\end{split}
\end{equation}
It is then clear that
\begin{equation}.
\begin{split}\max\limits_{1\leq i\leq k+1}|\sigma_{q_{n_m+1}}^2Q_{i}(z_0)|&\leq K\sigma_{q_{n_m+1}}^2\max\limits_{1\leq i\leq k}|\xi_i|^2\left( \hat{\Lambda}(z_0,k)\right)^3+\\
&+\max\limits_{1\leq i\leq k}|\sigma_{q_{n_m+1}}^2Q_{i}(z_0)|K\sigma_{q_{n_m+1}}\max\limits_{1\leq i\leq k}|\xi_i| \left(\hat{\Lambda}(z_0,k)\right)^2\\
&\leq K\sigma_{q_{n_m+1}}^2\max\limits_{1\leq i\leq k}|\xi_i|^2\left( \hat{\Lambda}(z_0,k)\right)^3+\\
&+\max\limits_{1\leq i\leq k+1}|\sigma_{q_{n_m+1}}^2Q_{i}(z_0)|K\sigma_{q_{n_m+1}}\max\limits_{1\leq i\leq k}|\xi_i| \left(\hat{\Lambda}(z_0,k)\right)^2.
\end{split}
\end{equation}
Hence we have the following bound
\begin{equation}\label{Max}
\begin{split}
\max\limits_{1\leq i\leq k+1}|\sigma_{q_{n_m+1}}^2Q_{i}(z_0)|&\leq K\sigma_{q_{n_m+1}}^2\max\limits_{1\leq i\leq k}|\xi_i|^2\left( \hat{\Lambda}(z_0,k)\right)^3+\\
&+\max\limits_{1\leq i\leq k+1}|\sigma_{q_{n_m+1}}^2Q_{i}(z_0)|K\sigma_{q_{n_m+1}}\max\limits_{1\leq i\leq k}|\xi_i| \left(\hat{\Lambda}(z_0,k)\right)^2,
\end{split}
\end{equation}
 Now we can prove  Lemma \ref{analogously}.\\
Using (\ref{CondQn}) and (\ref{Max}) we get
\begin{equation}
\begin{split}
|\sigma_{q_{n_m+1}}^2Q_{q_{n_m+1}}(z_0) \mathbf{1}_{D_{n_m}}|&\leq \max\limits_{1\leq i\leq q_{n_m+1}}|\sigma_{q_{n_m+1}}^2Q_{i}(z_0)\mathbf{1}_{D_{n_m}}|\\
&\leq 2K\sigma_{q_{n_m+1}}^2\max\limits_{1\leq i\leq q_{n_m+1}}|\xi_i|^2\left( \hat{\Lambda}(z_0,q_{n_m+1})\right)^3.
\end{split}
\end{equation}
Consequently
$$\left(E|\sigma_{q_{n_m+1}}^2Q_{q_{n_m+1}}(z_0) \mathbf{1}_{D_{n_m}}|^{p/2}\right)^{\frac{2}{p}}\leq 2K\sigma_{q_{n_m+1}}^2\left( \hat{\Lambda}(z_0,q_{n_m+1})\right)^3 \left(E(\max\limits_{1\leq i\leq q_{n_m+1}}|\xi_i|^p)\right)^{\frac{2}{p}}.$$
 Let $\varepsilon>0$ and using Chebeshev's inequality we obtain
$$\mathbf{P}\left(\left|\frac{\sigma_{q_{n_m+1}} Q_{q_{n_m+1}}(z_0, \sigma_{q_{n_m+1}})}{\sqrt{\mathrm{var} L_{q_{n_m+1}}(z_0)}}\mathbf{1}_{D_{n_m}}\right|>\varepsilon\right)\leq \left(\frac{2K\sigma_{q_{n_m+1}}\left( \hat{\Lambda}(z_0,q_{n_m+1})\right)^3 \left(E(\max\limits_{1\leq i\leq q_{n_m+1}}|\xi_i|^p)\right)^{\frac{2}{p}}}{\varepsilon\sqrt{\mathrm{var} L_{q_{n_m+1}}(z_0)}}\right)^{\frac{p}{2}}.$$
Conditions (\ref{Cond1}) and (\ref{Cond2}) imply
$$const \,\Lambda_2(z_0,q_{n_m+1})\leq \mathrm{var}\, L_{q_{n+1}}(z_0) \leq Const\, \Lambda_2(z_0,q_{n_m+1})$$
this and the previous bounds imply the  assertion of Lemma \ref{analogously}.
\endproof

\begin{lemm}\label{Dnm}  Suppose a circle map $T$ satisfies the conditions of Theorem \ref{main} and the sequence $\sigma_n$ satisfies the relation
\begin{equation}
\lim\limits_{n\to\infty}\sigma^p_{q_{n+1}} q_{n+1}\lambda_{-1}^{2np} n^{2p}=0.
\end{equation}
Then the probabilities of the events $D_{n_m}$ defined in (\ref{CondQn}) tend to 1 as  $m\to \infty$.
\end{lemm}
\proof
Consider the probability
\begin{equation}\label{CondQn1}
\begin{split}\mathbf{P}(\Omega\setminus D_{n_m})&=\mathbf{P}\left\{C\sigma_{q_{n_m+1}}\max\limits_{1\leq i\leq q_{n_m+1}}|\xi_i| \left(\hat{\Lambda}(z_0,q_{n_m+1})\right)^2\geq\frac{1}{2}\right\}\\
&=\mathbf{P}\left\{\max\limits_{1\leq i\leq q_{n_m+1}}|\xi_i|\geq\frac{1}{2 C\sigma_{q_{n_m+1}}\left(\hat{\Lambda}(z_0,q_{n_m+1})\right)^2}\right\}\\
 &\leq E(\max\limits_{1\leq i\leq q_{n_m+1}}|\xi_i|^p) \left(2 C\sigma_{q_{n_m+1}}\left(\hat{\Lambda}(z_0,q_{n_m+1})\right)^2\right)^p
\end{split}
\end{equation}
Condition (\ref{Cond2})  shows that
$$const\leq E(\max\limits_{1\leq i\leq q_{n+1}}|\xi_i|^p)\leq q_{n+1} Const$$
and   Lemma  \ref{hatLambda}   implies $\hat{\Lambda}(z_0,q_{n_m+1})\leq Const\cdot n_{m}\cdot \lambda_{-1}^{n_m} $.\\
Using these  bounds we get the following estimate for $\mathbf{P}(\Omega\setminus D_{n_m}):$
$$\mathbf{P}(\Omega\setminus D_n)\leq Const_1 \cdot q_{n_m+1}\cdot n_m^{2p}\cdot\lambda_{-1}^{2pn_m}\cdot\sigma_{q_{n_m+1}}^p,$$
which implies Lemma \ref{Dnm}.
\endproof
\begin{lemm}\label{lim} Assume a circle map $T$ satisfies the conditions of Theorem \ref{main} and the sequence $\sigma_n$ satisfies the relation: $$\lim\limits_{n \to\infty}\sigma_{q_{n+1}}\cdot n^3\cdot q_{n+1}^{2/p}\cdot\lambda_{-1}^{5n/2}\cdot \theta_+^{-2n}=0.$$ Then
$$\lim\limits_{m \to\infty}\frac{\sigma_{q_{n_m+1}}\left( \hat{\Lambda}(x_0,q_{n_m+1})\right)^3 \left(E(\max\limits_{1\leq i\leq q_{n_m+1}}|\xi_i|^p)\right)^{\frac{2}{p}}}{\sqrt{\Lambda_{2}(z_0, q_{n_m+1})}}=0.$$
\end{lemm}
\proof
In the Proof of  Lemma \ref{Dnm}
we showed
$$const\leq E(\max\limits_{1\leq i\leq q_{n+1}}|\xi_i|^p)\leq q_{n+1} \,Const$$
respectively
$$\hat{\Lambda}(z_0,q_{n_m+1})\leq C\cdot n_{m}\cdot \lambda_{-1}^{n_m}.$$
Since $\lambda_{-1}\leq\lambda_{-2} $ Theorem \ref{Lyap1} shows
$$c_{1}|I^{(n_m+1)}_{0}({x_{b}};z_{0})|^2 \lambda_{-2}^{n_m}\leq \Lambda_2(z_0,q_{n_m+1})\leq C_{1} |I^{(n_m+1)}_{0}({x_{b}};z_{0})|^2 \lambda_{-2}^{n_m} .$$
 Lemma \ref{bound} on the other hand implies
$$|I^{(n_m+1)}_{0}({x_{b}};z_{0})|\geq C_1 \theta_{+}^{n_m+1}.$$
Using these inequalities we find
$$\lim\limits_{m \to\infty}\frac{\sigma_{q_{n_m+1}}\left( \hat{\Lambda}(x_0,q_{n_m+1})\right)^3 \left(E(\max\limits_{1\leq i\leq q_{n_m+1}}|\xi_i|^p)\right)^{\frac{2}{p}}}{\sqrt{\Lambda_{2}(z_0, q_{n_m+1})}} \leq Const\frac{\sigma_{q_{n_m+1}} \cdot n_{m}^3\cdot \lambda_{-1}^{3n_m} q^{\frac{2}{p}}_{n+1} }{\sqrt{\theta_{+}^{2n_m+2} \lambda_{-2}^{n_m}}}\leq$$
$$\leq Const\frac{\sigma_{q_{n_m+1}} \cdot n_{m}^3\cdot \lambda_{-1}^{3n_m} q^{\frac{2}{p}}_{n+1} }{\theta_{+}^{2n_m+2} \lambda_{-1}^{n_m/2}}=Const\,\sigma_{q_{n_m+1}} \cdot n_{m}^3\cdot \lambda_{-1}^{5n_m/2} q^{\frac{2}{p}}_{n_m+1} \theta_{+}^{-2n_m-2} $$
\endproof
For the proof of Theorem \ref{main} we make use of
$$\lim\limits_{n\to\infty}\frac{\ln q_{n}}{n}=\ln\rho^{-1}_T,$$ where $\rho_T=[k_{1},k_{2},..,k_{m},1,1,...],\,\,m\geq1$.
Thus for arbitrary $0<\varepsilon<\frac{1}{2}$ there exists $N\in\mathbb{N}$ such that  for all $n>N$ $$\rho_T^{-n/2}<\rho_T^{-n(1-\varepsilon)}<q_n<\rho_T^{-n(1+\varepsilon)}<\rho_T^{-3n/2}.$$
Define  $\gamma:=\max\{\frac{2}{p}-\frac{5\ln\lambda_{-1}+4\ln\theta^{-1}_{+}+6}{\ln\rho_T},\frac{1}{2}+\frac{2\ln\theta_{+}}{\ln\rho_T}\}$. In order to prove the first part of Theorem  \ref{main}, it is enough to verify the conditions of Lemmas \ref{Bnm}, \ref{Dnm} and \ref{lim}.\\
\textbf{1. Condition of Lemma \ref{Bnm}.} For sufficiently large $n$
  $$\frac{q_{n+1}\sigma^2_{q_{n+1}}}{\theta_{+}^{2n+4}}=q_{n+1}\sigma^2_{q_{n+1}}\theta_{+}^{-2n-4}= q_{n+1}\sigma^2_{q_{n+1}}\rho_T^{\frac{\ln\theta_{+}^{-1}}{\ln\rho_T}(2n-4)}<$$
  $$<\left(q_{n+1}^{\frac{1}{2}-\frac{\ln\theta_{+}^{-1}}{\ln\rho_T}\frac{n-2}{(n+1)(1-\varepsilon)}}\sigma_{q_{n+1}}\right)^2<
   \left(q_{n+1}^{\gamma}\sigma_{q_{n+1}}\right)^2.$$
\textbf{2. Condition of Lemma \ref{Dnm}.}
$$\sigma^p_{q_{n+1}} q_{n+1}\lambda_{-1}^{2np} n^{2p}=\left(\sigma_{q_{n+1}} \cdot q_{n+1}^{1/p}\cdot\lambda_{-1}^{2n}\cdot n^{2}\right)^p< \left(\sigma_{q_{n+1}}\cdot q_{n+1}^{2/p}\cdot\lambda_{-1}^{5n/2}\cdot n^3\cdot \theta_+^{-2n}\right)^p.$$
\textbf{3. Condition of Lemma \ref{lim}.} For sufficiently large $n$
$$\sigma_{q_{n+1}}\cdot q_{n+1}^{2/p}\cdot\lambda_{-1}^{5n/2}\cdot n^3\cdot \theta_+^{-2n}= \sigma_{q_{n+1}}\cdot q_{n+1}^{2/p}\cdot\rho_T^{\frac{\ln\lambda_{-1}}{\ln\rho_T}\cdot \frac{5n}{2}}\cdot \rho_T^{\frac{3\ln n}{\ln\rho_T}}\cdot \rho_T^{2n\cdot\frac{\ln\theta^{-1}_+}{\ln\rho_T}}<$$
$$<\sigma_{q_{n+1}}\cdot q_{n+1}^{2/p}\cdot q_{n+1}^{-\frac{\ln\lambda_{-1}}{\ln\rho_T}\cdot \frac{5n}{2(n+1)(1-\varepsilon)}}\cdot q_{n+1}^{-\frac{3\ln n}{\ln\rho_T(n+1)(1-\varepsilon)}}\cdot q_{n+1}^{-\frac{2n}{(n+1)(1-\varepsilon)}\cdot\frac{\ln\theta^{-1}_+}{\ln\rho_T}}<$$
$$<\sigma_{q_{n+1}}\cdot q_{n+1}^{2/p-\frac{\ln\lambda_{-1}}{\ln\rho_T}\cdot \frac{5n}{(n+1)}-\frac{4n}{(n+1)}\cdot\frac{\ln\theta^{-1}_+}{\ln\rho_T}-\frac{6\ln n}{(n+1)\ln\rho_T}}<\sigma_{q_{n+1}}\cdot q_{n+1}^{\gamma}.$$
Thus
\begin{equation}\label{cond5.3}
\sigma_{q_{n+1}} \cdot n^3\cdot \lambda_{-1}^{5n/2} q^{\frac{2}{p}}_{n+1} \theta_{+}^{-2n-2}<\sigma_{q_{n+1}}\cdot q_{n+1}^{\gamma}.
\end{equation}
For the second part of Theorem \ref{main} it is enough to show that the rate of convergence to $0$ of $\sigma_{q_{n_m+1}-1}
\frac{Q_{q_{n_m+1}-1}(z_0,\sigma_{q_{n_m+1}-1})}{\sqrt{\mathrm{var}(L_{q_{n_m+1}-1}(z_0))}}$  is bounded by the right hand side in the estimate in (\ref{BerreEsseyn}), that means
$$\sigma_{q_{n+1}} \cdot n^3\cdot \lambda_{-1}^{5n/2} q^{\frac{2}{p}}_{n+1} \theta_{+}^{-2n-2}<\frac{\Lambda_{\min(p,3)}(x,q_n)}{(\Lambda_2(x,q_n))^{\min(p,3)/2}} .$$
To prove this inequality we use (\ref{cond5.3}) and  Theorem \ref{Lyap1} to show
$$const\cdot\left(\frac{\lambda_{-\min(p,3)}^2}{\lambda_{-2}^{\min(p,3)}}\right)^{\frac{n}{2}}\leq \frac{\Lambda_{\min(p,3)}(x,q_n)}{(\Lambda_2(x,q_n))^{\min(p,3)/2}}.$$
 It is therefore enough to get
\begin{equation}\label{last}\sigma_{q_{n+1}}\cdot q_{n+1}^{\gamma}\leq const\cdot\left(\frac{\lambda_{-\min(p,3)}^2}{\lambda_{-2}^{\min(p,3)}}\right)^{\frac{n}{2}},\end{equation}
when $\sigma_{q_{n+1}}\leq C_1 \cdot q_{n+1}^{-\tau}$.
%$$\sigma_{q_{n+1}} \cdot n^3\cdot \lambda_{-1}^{5n/2} q^{\frac{2}{p}}_{n+1} \theta_{+}^{-2n-2}<\sigma_{q_{n+1}}\cdot q_{n+1}^{\gamma}\leq const\cdot\left(\frac{\lambda_{-\min(p,3)}^2}{\lambda_{-2}^{\min(p,3)}}\right)^{\frac{n}{2}}\leq \frac{\Lambda_{\min(p,3)}(x,q_n)}{(\Lambda_2(x,q_n))^{\min(p,3)/2}} $$
Let us choose $s=\min(p,3)$ and  $\tau\geq\gamma {+}\frac{2\ln\lambda_{-s}-s\ln\lambda_{-2}}{3\ln\rho_T}$.
Then
$$\sigma_{q_{n+1}}\leq C_1 \cdot q_{n+1}^{-\tau}\leq C_1\cdot q_{n+1}^{-\gamma{-}\frac{2\ln\lambda_{-s}-s\ln\lambda_{-2}}{3\ln\rho_T}}\leq C_1\cdot q_{n+1}^{-\gamma}\cdot \rho_T^{\frac{2\ln\lambda_{-s}-s\ln\lambda_{-2}}{3\ln\rho_T}\cdot\frac{3(n+1)}{2}},$$
from which (\ref{last}) follows.\hfill$\Box$\newline

To finish the proof we can next apply Lemma 3.2 in \cite{Llave} by noting that there the quantity $||f''||_{C_0}$ has to be replaced in our case 
%by $\sup\limits_{z\in [x_b,x_b+1]} T^''(z)$ 
by $\sup\limits_{z\in[x_b,x_b+1]}|T''(z)|$
which exists for $T\in C^{2+\epsilon}(S^1\setminus\{x_b\})$. Following then the arguments of 
the proof of the CLT in Section 3.3 of \cite{Llave} leads finally to the convergence of the process $\omega_{n_k}(z_0,\sigma_{n_k})$ to the 
standard Gaussian distribution.

{\bf{Acknowledgement}} The authors thank a referee for several very helpful remarks.

\label{lastpage}

\end{document}